\documentclass[a4paper] {article}

	\usepackage {amsmath}
	\usepackage {amscd}
	\usepackage {amssymb}
	\usepackage {amsthm}
	
	\usepackage {dsfont}

	\usepackage[T1]{fontenc}
	\usepackage[utf8]{inputenc}

	\usepackage {enumerate}

	\usepackage{titling}
\newcommand{\subtitle}[1]{%
  \posttitle{%
    \par\end{center}
    \begin{center}\large#1\end{center}
    \vskip0.5em}%
}

	\newtheorem{TheoremABC}{Theorem}

	\newtheorem {Theorem}{Theorem}[section]
	\newtheorem {Cor}[Theorem]{Corollary}
	\newtheorem {Lemma}[Theorem]{Lemma}
	\newtheorem {Proposition}[Theorem]{Proposition}
	\newtheorem {Definition}[Theorem]{Definition}

	\theoremstyle{definition}
	\newtheorem {Remark}[Theorem]{Remark}

\begin {document}

\title{A survey of the GIT picture for the Yang-Mills equation over Riemann surfaces}

\author{Samuel Trautwein\thanks{Partially supported by the Swiss National Science Foundation Grant 156000.}}

\maketitle

\abstract{The purpose of this paper is to give a self-contained exposition of the Atiyah-Bott picture \cite{AtBott:YangMillsEq} for the Yang-Mills equation over Riemann surfaces with an emphasis on the analogy to finite dimensional geometric invariant theory. The main motivation is to provide a careful study of the semistable and unstable orbits: This includes the analogue of the Ness uniqueness theorem for Yang-Mills connections, the Kempf-Ness theorem, the Hilbert-Mumford criterion and a new proof of the moment-weight inequality following an approach outlined by Donaldson \cite{Donaldson:2005}. A central ingredient in our discussion is the Yang-Mills flow for which we assume longtime existence and convergence (see \cite{Rade1992}).
}

\newpage

\tableofcontents

\newpage


\section{Introduction}
The purpose of this paper is threefold: The first goal is to provide a self-contained and essentially complete exposition of the geometric invariant theory for the Yang-Mills equation over Riemann surfaces from the differential geometric point of view. We follow closely the line of arguments of finite dimensional GIT and emphasize this analogy throughout (see e.g. the exposition \cite{RobSaGeo}). There is no claim of originality. The main motivation has been that the underlying ideas are spread over the literature and it seems worthwhile to provide an unified exposition. 

The second goal is to include a careful study of the semistable and unstable orbits. This is in contrast to most of the developments after the landmark paper \cite{AtBott:YangMillsEq} of Atiyah and Bott, which deal with the characterization of stable objects in more general moduli problems, i.e. the analogue of the Narasimhan-Seshadri-Ramanathan theorem. In the unitary case Daskalopoulos \cite{Dask:1992} established the Morse theoretic picture of Atiyah and Bott. A direct corollary of this stratification is the analogue of the Ness uniqueness theorem and the moment limit theorem (see Theorem A below). We give an alternative proof of this result following the ideas from finite dimensional GIT. This argument does not depend on the Harder-Narasimhan filtration or on other aspects from the holomorphic point of view and works for general structure groups. Following an approach outlined by Donaldson \cite{Donaldson:2005}, we also carry out a new proof of the moment-weight inequality which is essentially contained in the work of Atiyah and Bott.

The third goal is to provide a transparent exposition of the central ideas used in gauge theoretical moduli problems. While several results are known in greater generality, the key ideas are still immanent in our treatment. We hope that this enables non experts to explore the beauty of this subject without having to worry about the technical difficulties which come along with more general moduli problems. On the other hand, we expect that some of the well-known results in finite dimensional GIT quickly lead to difficult and open questions in higher dimensional moduli problems.
\\

Let $G$ be a compact connected Lie group and let $\Sigma$ be a closed Riemann surface. Fix a volume form on $\Sigma$, compatible with the orientation, and let $P \rightarrow \Sigma$ be a principal $G$ bundle. Atiyah and Bott \cite{AtBott:YangMillsEq} observed that the curvature
					$$\mu(A) := *F_A \in \Omega^0(\Sigma,\text{ad}(P))$$
defines a moment map for the action of the gauge group $\mathcal{G}(P)$ on the space of connections $\mathcal{A}(P)$. For any constant central section $\tau$, the symplectic quotient
				$$\mathcal{A}(P)/\!/\mathcal{G}(P) := \mu^{-1}(\tau)/\mathcal{G}(P)$$
yields the moduli space of projectively flat connections on $P$ with constant central curvature $\tau$.

There is a general procedure relating symplectic quotients with quotients in algebraic geometry (see \cite{RobSaGeo}). For this let $G^c$ be the complexification of $G$ and $P^c := P\times_G G^c$ the associated principal $G^c$ bundle. The complexification of the gauge group is $\mathcal{G}^c(P) := \mathcal{G}(P^c)$. The space $\mathcal{A}(P)$ can naturally be identified with the space $\mathcal{J}(P^c)$ of holomorphic structures on $P^c$ (see Lemma \ref{PreLem2}) and the complexified gauge group $\mathcal{G}^c(P)$ acts naturally on this space. The corresponding GIT quotient
											$$\mathcal{A}^{ss}(P)/\!/\mathcal{G}^c(P)$$
of $\mathcal{A}(P) $ by $\mathcal{G}^c(P)$ is obtained in two steps. First, one defines a dense and open subset $\mathcal{A}^{ss}(P) \subset \mathcal{A}(P)$ of semistable connections or holomorphic structures on $E$ and second, one identifies two semistable orbits in the quotient if they cannot be separated in $\mathcal{A}(P)$. The restriction to semistable orbits is necessary to obtain a \textit{good quotient} in the sense of algebraic geometry. There are essentially two approaches to define semistable objects. In the symplectic approach, one uses a moment map for the gauge action on $\mathcal{A}(P)$ to define semistable objects. In the classical algebraic geometric approach, one defines a notion of semistability $\mathcal{J}^{ss}(P^c) \subset \mathcal{J}(P^c)$ on the space of holomorphic structures on $P^c$. A classical result due to Narasimhan and Seshadri \cite{NarasimhanSeshadri:1965} in the case $G = U(n)$ and due to Ramanathan \cite{Ramanathan:1975} for general $G$ shows that both of these notions essentially agree if one restricts to further open subsets of stable objects.\\

The Yang-Mills picture introduced by Atiyah and Bott \cite{AtBott:YangMillsEq} shed new light on this result and inspired Donaldson \cite{Donaldson:1983NS} to an analytic proof of the Narasimhan-Seshadri theorem. The Yang-Mills functional is given by the formula 
		$$\mathcal{YM}: \mathcal{A}(P) \rightarrow \mathbb{R}, \qquad \mathcal{YM}(A) := \frac{1}{2}\int_{\Sigma} ||F_A||^2 \, \text{dvol}_{\Sigma}.$$
Standard arguments from Chern-Weil theory show that there exists an unique central element $\tau \in Z(\mathfrak{g})$ such that 
				\begin{align} \label{Inteq1} \mathcal{YM}(A) = \inf_{B \in \mathcal{A}(P)} \mathcal{YM}(B) \qquad \Longleftrightarrow \qquad *F_A = \tau. \end{align}
We call $\tau$ the central type of $P$. We shall consider in the following connections of Sobolev class $W^{1,2}$ and gauge transformations of Sobolev class $W^{2,2}$. Rade \cite{Rade1992} showed in this setting that for every initial data $A_0 \in \mathcal{A}(P)$ the negative gradient flow of the Yang-Mills functional
				\begin{align} \label{Inteq2} \partial_t A(t) = - \nabla \mathcal{YM}(A(t)) = - d_{A(t)}^* F_{A(t)}, \qquad A(0) = A_0 \end{align}
has an unique (weak) solution which exists for all time. Moreover, this solution remains in a single complexified $\mathcal{G}^c(P)$-orbit and converges in the $W^{1,2}$-topology to a Yang-Mills connection $A_{\infty} \in \overline{\mathcal{G}^c(A_0)}$. The following is the analogue of the Ness Uniqueness theorem in finite dimensional GIT. 

\begin{TheoremABC}[\textbf{Uniqueness of Yang-Mills connections}]
	Let $A_0 \in \mathcal{A}(P)$ and $A_{\infty}$ be the limit of the Yang-Mills flow (\ref{Inteq2}) starting at $A_0$. Then
		\begin{enumerate}
				\item $\mathcal{YM}(A_{\infty}) = \inf_{g \in \mathcal{G}^c(P)} \mathcal{YM}(gA)$.
				
				\item If $B \in \overline{\mathcal{G}^c(A_0)}$ is contained in the $W^{1,2}$-closure of $\mathcal{G}^c(A_0)$ and
					$$\mathcal{YM}(B) = \inf_{g \in \mathcal{G}^c(P)} \mathcal{YM}(gA)$$
			then $\mathcal{G}(B) = \mathcal{G}(A_{\infty})$.
		\end{enumerate}
\end{TheoremABC}

In the case of a vector bundle $E \rightarrow \Sigma$ (i.e. $G = U(n)$) Daskalopoulos \cite{Dask:1992} established the convergence of the Yang-Mills flow over Riemann surfaces by different methods. He proves a suitable slice theorm near Yang-Mills connections and shows that the limiting Yang-Mills connection $A_{\infty}$ is determined up to an unitary gauge transformation by the isomorphism class of the Harder-Narasimhan filtration of $(E, \bar{\partial}_{A_0})$. This proves Theorem A in the unitary case and it should be possible to deduce the general result from this using the methods in \cite{BiswanWilkin:2010}. We present a different proof of Theorem A in Theorem \ref{MLT} and Theorem \ref{ThmNUgen} by following the line of arguments from finite dimensional GIT (\cite{RobSaGeo}, Chapter 6). These arguments were originally given by Calabi-Chen \cite{CalabiChen:2002} and Chen-Sun \cite{ChenSun:2010} in the context of extremal Kähler metrics. 

A connection $A \in \mathcal{A}(P)$ is called $\mu_{\tau}$-semistable resp. $\mu_{\tau}$-unstable if
					$$\inf_{g \in \mathcal{G}^c(P)} ||*F_{gA} - \tau||_{L^2} = 0 \qquad \text{resp.} \qquad \inf_{g \in \mathcal{G}^c(P)} ||*F_{gA} - \tau||_{L^2} > 0$$
where $\tau$ is defined by (\ref{Inteq1}). Moreover, $A$ is called $\mu_{\tau}$-polystable if there exists $g \in \mathcal{G}^c(P)$ with $*F_{gA} = \tau$ and it is called $\mu_{\tau}$-stable if $gA$ is in addition irreducible. Then Theorem A implies that the map which sends $A_0 \in \mathcal{A}(P)$ to the limit $A_{\infty}$ of the Yang-Mills flow starting at $A_0$ yields the identifications
					$$\mathcal{A}^{ss}(P)/\!/\mathcal{G}^c(P) \cong \mathcal{A}^{ps}(P)/\mathcal{G}^c(P) \cong \mu^{-1}(\tau)/\mathcal{G}(P).$$
Conversely, the $\mu_{\tau}$-unstable orbits converge to higher critical points of the Yang-Mills functional. More details on this correspondence are given in Theorem \ref{ThmYMStab}.\\

The theory has greatly evolved since the paper \cite{AtBott:YangMillsEq} of Atiyah and Bott. The main goal in those developments has been the characterization of stable objects in more general moduli problems (e.g. \cite{Donaldson:ASD4}, \cite{Donaldson:1987}, \cite{UYau:1986}, \cite{Hitchin:1987}, \cite{Simpson:1987}, \cite{Bradlow:1991}). The characterization of unstable orbits is in general much more difficult as it refers to higher critical points of the Yang-Mills functional. Given a connection $A \in \mathcal{A}(P)$ and $\xi \in \Omega^0(\Sigma, \text{ad}(P))$ the weight $w_{\tau}(A,\xi)$ is defined by
				$$w_{\tau}(A,\xi) := \lim_{t\rightarrow \infty} \langle *F_{e^{\textbf{i}t\xi}A} - \tau, \xi \rangle \in \mathbb{R}\cup\{\infty\}.$$
The first part of the following theorem is the analogue of the moment-weight inequality and the last two claims are the analogue of the Kempf existence and uniqueness theorem in finite dimensional GIT.

\begin{TheoremABC}[\textbf{Atiyah-Bott}]
Let $A \in \mathcal{A}(P)$ and let $\tau \in Z(\mathfrak{g})$ be defined by (\ref{Inteq1}). Then
		\begin{enumerate}
				\item For all $0 \neq \xi \in \Omega^0(\Sigma, \text{ad}(P))$ holds
								$$- \frac{w_{\tau}(A,\xi)}{||\xi||} \leq \inf_{g \in \mathcal{G}^c(P)} ||*F_A - \tau||^2.$$
				\item If the right-hand-side is positive, then there exists up to scaling an unique $0\neq \xi_0 \in \Omega^0(\Sigma,\text{ad}(P))$ such that
								$$- \frac{w_{\tau}(A,\xi_0)}{||\xi_0||} = \inf_{g \in \mathcal{G}^c(P)} ||*F_A - \tau||^2.$$
					Moreover, $\xi_0$ is rational in the sense that it generates a closed one parameter subgroup of $\mathcal{G}(P)$.
				\item Let $A_{\infty}$ be the limit of the Yang-Mills flow starting at $A_0$. Then there exists $u \in \mathcal{G}(P)$ such that $\xi_0$ agrees up to scaling with $u(*F_{A_{\infty}} - \tau)u^{-1}$.
		\end{enumerate}
\end{TheoremABC}

This is essentially contained in the work of Atiyah and Bott (\cite{AtBott:YangMillsEq}, Prop. 8.13 and Prop. 10.13). A connection $A \in \mathcal{A}(P)$ induces a holomorphic structure on the complexified bundle $P^c := P \times_G G^c$ and its Lie algebra bundle $\text{ad}(P^c)$. Atiyah and Bott explicitly determine the infimum of the Yang-Mills functional over $\mathcal{G}^c(A)$ in terms of the Harder-Narasimhan filtration of $\text{ad}(P^c)$. The analogous result has been shown by Calabi, Chen, Donaldson and Sun \cite{CalabiChen:2002, Donaldson:2005, Chen:2008, Chen:2009, ChenSun:2010} in the context of extremal Kähler metrics. Donaldson \cite{Donaldson:2005} compares the Atiyah-Bott picture in the vector bundle case $G=U(n)$ with their results on the Calabi functional and mentions that their methods should lead to a new proof of the moment-weight inequality in the Atiyah-Bott case. We carry out this proof in Theorem \ref{thmMWI}. The same arguments are used in \cite{RobSaGeo} to establish the finite dimensional version of the moment-weight inequality. We reformulate and prove the last two claims in Theorem \ref{ThmDomWeight}. The case $G=U(n)$ follows along the line of arguments of Atiyah and Bott from the Harder-Narasimhan filtration and the Narasimhan-Seshadri theorem. The general case can be reduced to this by the use of Theorem A. For this, choose a faithful representation $G \hookrightarrow U(n)$. Then any $G$-connection $A$ can be considered as $U(n)$-connection and Theorem A implies
					$$\inf_{g\in \mathcal{G}^c(E)} \mathcal{YM}(gA) = \inf_{g \in \text{GL}(E)} \mathcal{YM}(gA).$$
It now remains to compare the weights for the gauge action with respect to the two structure groups $G$ and $U(n)$ to conclude the proof. We would also like to mention the work of Bruasse and Teleman \cite{BruasseTeleman:2003, Bruasse:2006}. They prove for more general gauge theoretical moduli problems that whenever the supremum over the normalized weights is positive, then it is attained in an unique direction corresponding to the Harder-Narasimhan filtration.\\

There is a classical algebraic geometric notion of stability for holomorphic principal bundles (see Definition \ref{DefnStab2}). In the vector bundle case $G = U(n)$ this corresponds to the notion of (slope-)stable holomorphic vector bundles, which are easier to define: A holomorphic vector bundle $E$ is called stable (semistable) if
					$$\frac{c_1(F)}{\text{rk}(F)} < \frac{c_1(E)}{\text{rk}(E)} \qquad \left( \frac{c_1(F)}{\text{rk}(F)} \leq \frac{c_1(E)}{\text{rk}(E)} \right)$$
holds for every proper holomorphic subbundle $0 \neq F \subset E$. Moreover, $E$ is called polystable if it decomposes as the direct sum of stable vector bundles all having the same slope and $E$ is called unstable if it is not semistable.

\begin{TheoremABC}[\textbf{Generalized Narasimhan-Seshadri-Ramanathan theorem}]
Let $A \in \mathcal{A}(P)$ and define $\tau$ by (\ref{Inteq1}). Then $A$ induces a holomorphic structure $J_A$ on the complexified bundle $P^c := P \times_G G^c$ and the following holds:
		\begin{enumerate}
				\item $(P^c,J_A)$ is stable if and only if $A$ is $\mu_{\tau}$-polystable and the kernel of the infinitesimal action
										$L_A: \Omega^0(\Sigma, \text{ad}(P^c)) \rightarrow \Omega^1(\Sigma, \text{ad}(P))$
														$$L_A(\xi + \textbf{i}\eta) := - d_A \xi - * d_A \eta$$
					contains only constant central sections.					
				\item $(P^c, J_A)$ is $\mu_{\tau}	$-polystable if and only if $A$ is $\mu_{\tau}$-polystable.
				\item $(P^c, J_A)$ is $\mu_{\tau}	$-semistable if and only if $A$ is $\mu_{\tau}$-semistable.
				\item $(P^c, J_A)$ is $\mu_{\tau}	$-unstable if and only if $A$ is $\mu_{\tau}$-unstable.

		\end{enumerate}
\end{TheoremABC}

The first part of this Theorem is the Narasimhan-Seshadri-Ramanathan theorem. We present an analytic proof of this result in Theorem \ref{ThmNSR} which was originally given by Bradlow \cite{Bradlow:1991} and Mundet \cite{Mundet:2000} for more general moduli problems. The main step in their proof is to establish the analogue of the Kempf-Ness theorem (see Theorem \ref{ThmKNgen}) in the stable case. The polystable case is deduced from the stable case by induction on the dimension of $G$. We characterize in Proposition \ref{PropStabChar} the stability of $(P^c, J_A)$ in terms of the weights $w_{\tau}(A,\xi)$. In particular, $(P^c, J_A)$ is unstable if and only if there exists a negative weight $w_{\tau}(A,\xi) < 0$. Therefore, the unstable and semistable case of Theorem C follow directly from Theorem B. We reformulate and prove Theorem C in Theorem \ref{ThmNSRgen}.

\subsubsection*{Outline}

In Section 2 we review the necessary preliminaries. The first part deals with the relevant background on gauge theory. Besides fixing notation, the main goals are to provide an explicit description of the complexified gauge action in both the vector bundle and principal bundle case and to describe the moment map picture of Atiyah and Bott. We show that this picture remains valid if one considers connections and gauge transformations in suitable Sobolev completions. The second part discusses parabolic subgroups of complex reductive Lie groups. These play a crucial role in the algebraic geometric definition of stability and the geometric description of the weights.

In Section 3 we discuss the algebraic and symplectic definitions of stability. The main result in this section is the generalized Narasimhan-Seshadri-Ramanathan theorem (Theorem \ref{ThmNSRgen}) which states that these definitions are essentially equivalent. The proof of this theorem is based on the whole remainder of the exposition. 

In Section 4 we review the analytical properties of the Yang-Mills flow which Rade \cite{Rade1992} established in his thesis. Following the line of arguments in (\cite{RobSaGeo}, Chapter 6), we prove Theorem A in Theorem \ref{MLT} and Theorem \ref{ThmNUgen}. We close this section with Theorem \ref{ThmYMStab} which characterizes the $\mu_{\tau}$-stability of a connection $A \in \mathcal{A}(P)$ in terms of the limit $A_{\infty}$ of the Yang-Mills flow starting at $A$.

In Section 5 we introduce the weights $w_{\tau}(A,\xi)$ and show that they are closely related to holomorphic parabolic reductions of the complexified bundle $(P^c, J_A)$. Proposition \ref{PropStabChar} shows that the weights provide an alternative describes of the algebraic notion of stability. We close this section with the proof of the moment weight inequality (Theorem \ref{thmMWI}) following the approach outlined by Donaldson \cite{Donaldson:2005}.

In Section 6 we describe a general procedure which associates to a given connection $A \in \mathcal{A}(P)$ a $\mathcal{G}(P)$-invariant functional $\Phi_A : \mathcal{G}^c(P) \rightarrow \mathbb{R}$. We call this the Kempf-Ness functional of $A$. The slope of this functional at infinity agrees with the weights discussed in Chapter 5 and hence relates to the algebraic notion of stability by Proposition \ref{PropStabChar}. The analogue of the Kempf-Ness theorem (see Theorem \ref{ThmKNgen}) relates the global behavior of $\Phi_A$ to the symplectic $\mu_{\tau}$-stability of $A$. This provides a link between the algebraic and symplectic notions of stability and leads to an analytic proof of the Narasimhan-Seshadri-Ramanathan theorem in Theorem \ref{ThmNSR}. These arguments are given by Bradlow \cite{Bradlow:1991} and Mundet \cite{Mundet:2000} in more general settings. 

In Section 7 we establish the analogue of the Kempf existence and uniqueness theorem (see Theorem \ref{ThmDomWeight}). We include a self-contained account on the Harder-Narasimhan filtration for the convenience of the reader. This is the main new ingredient needed for the proof.

\subsubsection*{Higher dimensional base manifolds}
We restrict our discussion to the case where $\Sigma$ is a Riemann surface, although several results remain valid in greater generality. The main reason for this is to simplify the presentation. Let us indicate in the following to which degree the discussion could be generalized.

Replace $\Sigma$ by a closed Kähler manifold $(X,J,\omega)$ and denote by
			$$\Lambda : \Omega^{1,1}(X) \rightarrow \Omega^0(X)$$
the adjoint operator of $f \mapsto f\omega$. The Hermitian Einstein equation is given by
					$$\Lambda F_A = \tau$$
for some constant central element $\tau \in \Omega^0(X,\text{ad}(P))$. Denote by $\mathcal{A}^{1,1}$ the space of connections on $P$ whose curvature $F_A$ is of type $(1,1)$. This space can be given a Kähler structure and $\mu(A) = \Lambda F_A$ yields a moment map for the gauge action. In the vector bundle case, the Narasimhan-Seshadri theorem has been generalized to this setting by Donaldson \cite{Donaldson:ASD4, Donaldson:1987} in the algebraic framework and by Uhlenbeck and Yau \cite{UYau:1986} in the analytic framework over arbitrary Kähler manifolds. We would like to point out an observation by Anouche and Biswan \cite{AnchoucheBiswan:2001}. They show that a holomorphic principal bundle $P^c$ is polystable (resp. semistable), if and only if the associated holomorphic vector bundle $\text{ad}(P^c)$ is polystable (resp. semistable). Further generalizations involving more complicated moduli problems have been studied by Hitchin \cite{Hitchin:1987}, Simpson \cite{Simpson:1987} and Bradlow \cite{Bradlow:1991}. In his thesis \cite{Mundet:2000}, Mundet generalizes this correspondence to a very general moduli problem.

Our discussion of the Yang-Mills flow in Chapter 4 relies heavily on the fact that $\Sigma$ is a Riemann surface. In particular, the group of $W^{2,2}$ gauge transformations acts no more continuously on the space of $W^{1,2}$ connections for higher dimensional base manifolds. To avoid this issue, one could consider the flow directly on the space of smooth connections. Donaldson showed in \cite{Donaldson:ASD4} that the Yang-Mills flow starting at smooth $\mathcal{A}^{1,1}$ connections admits a smooth solution which exists for all time. In the stable case, Donaldson used this flow to prove his extension of the Narasimhan-Seshadri theorem. See \cite{Siu:1987} for a survey on this approach. The main issue is the complicated limiting behavior of solutions which yields profound technical difficulties. Bando and Siu (\cite{BandoSiu:1994}, Theorem 4) showed that the limit "breaks up" into Hermitian-Einstein sheaves in the unstable case and conjectured that the limit corresponds essentially to the Harder-Narasimhan filtration. This is very similar to our discussion in Chapter 7. The Bando-Siu conjecture has been confirmed by Daskalopoulos-Wentworth \cite{DaskWentworth:2004} in the case of Kähler surfaces and by Sibling \cite{Sibley:2015} and Jacob \cite{Jacob:2013, Jacob:2013b} for general Kähler manfiolds. This yields the analogue of Theorem C for vector bundles over Kähler manifolds.

Our calculation of the weights in Chapter 5 remains valid over arbitrary Kähler manifold. However, the weakly holomorphic filtration yields in this case only a filtration by torsion-free subsheaves. The proof of the moment-weight inequality generalizes ad verbatim to this case. The proof which we present for the Narasimhan-Seshadri-Ramanathan theorem remains valid in this setting as well (see \cite{Mundet:2000}).

The Harder-Narasimhan filtration is well defined for holomorphic vector bundles over Kähler manifolds, but consists of torsion-free subsheaves instead of holomorphic subbundles. It corresponds again to the supremum over the normalized weights. This is shown by Bruasse \cite{Bruasse:2006} and we present part of his argument in Chapter 7. It is a nontrivial result that the infimum of $||\Lambda F_{gA}||$ over the (smooth) complexified gauge orbit yields the same value and follows from the Bando-Siu conjecture. Bruasse gives an alternative and direct argument to prove that the supremum is in fact attained.

\subsubsection*{General assumptions}
Let $G$ be a compact connected (real) Lie group, $\Sigma$ a closed Riemann surface and $P \rightarrow \Sigma$ a principal $G$ bundle. We fix a volume form $dvol_{\Sigma}$ on $\Sigma$ and assume for convenience that the volume form is scaled such that
			$$\text{vol}(\Sigma) = 1.$$
Note that the volume form also induces a fixed Riemannian metric on $\Sigma$.

All Lie groups which we encounter in our discussion are assumed to be connected. When $G$ is a compact connected Lie group, then its complexification $G^c$, its parabolic subgroups $Q(\zeta)$ and their Levi subgroups $L(\zeta)$ are automatically connected (see Lemma \ref{pblLem1}).

As a general rule, we consider connections of Sobolev class $W^{1,2}$ and gauge transformations of Sobolev class $W^{2,2}$. The gauge action extends smoothly over these Sobolev spaces, since the base manifold is a Riemann surface. These regularity assumptions do not affect the overall picture and we shall discuss them in more detail in the preliminaries below.

\subsubsection*{Acknowledgement}
I would like thank my supervisor D.A Salamon for many helpful discussions throughout the process of writing this paper.


\section{Preliminaries}
First, we review the underlying notions from gauge theory and set up our notation. The main goal is to describe the complexification of the gauge action and the moment map picture of Atiyah and Bott. We also discuss the regularity assumptions which are crucial for our further analytic discussion. In the second subsection, we describe parabolic subgroups of complex reductive Lie groups. We also include a brief discussion of the root space decomposition of semisimple Lie algebras for the sake of completeness. 

\subsection{Gauge theory}

We consider throughout this section fiber bundles over a closed connected Riemann surface $\Sigma$, although this is not really necessary for most of our discussion.
 
\subsubsection{Basic gauge theory}

We start with the general framework of fiber bundles and specialize our discussion afterwards to the cases of vector bundles and principal bundles.

\paragraph{Fiber bundles.}
Let $E$, $F$ and $B$ be smooth manifolds. The manifold $E$ together with a projection map $\pi: E \rightarrow B$ is called a fiber bundle over $B$ with fiber $F$, if for every $x \in B$ there exists a neighborhood $x\in U \subset B$ and a diffeomorphisms
					$$\psi: \pi^{-1}(U) \rightarrow U\times F$$
such that $\text{pr}_1 \circ \psi = \pi|_U$. Here $\text{pr}_1: U\times F \rightarrow U$ denotes the projection onto the first factor. The map $\psi$ is called a local trivialization of the fiber bundle $E$. Suppose $\psi_{\alpha}$ and $\psi_{\beta}$ are local trivializations over $U_{\alpha}$ and $U_{\beta}$. Then there exists a unique map $g_{\beta \alpha} : U_{\alpha}\cap U_{\beta} \rightarrow \text{Diff}(F)$ satisfying
					$$\psi_{\beta \alpha}(x,u) := (\psi_{\beta}\circ\psi_{\alpha}^{-1})(x, u) = (x, g_{\beta \alpha}(x) u)$$
for all $x \in U_{\alpha} \cap U_{\beta}$	and $u \in F$. A reduction of the structure group of $E$ to a subgroup $G \subset \text{Diff}(F)$ consists of an open cover $\{ U_{\alpha} \}$ of $B$ together with local trivializations $\psi_{\alpha}$ such that all transition maps $g_{\beta \alpha}$ take values in $G$. The bundle $E$ together with a fixed choice of such trivialization is called as fiber bundle with structure group $G$.

The tangent bundle $TE$ contains a canonical vertical subbundle $V := \ker{d\pi}$. A connection on $E$ is a splitting of the exact sequence
				$$0 \rightarrow V \rightarrow TE \rightarrow TE/V \rightarrow 0 $$
and corresponds to a horizontal distribution $H \subset TE$ satisfying $TE = H \oplus V$. Identifying $H$ with the projection of $TE$ onto $V$, we can describe a connection by a $V$-valued $1$-form $A \in \Omega^1(E,V)$. The curvature of a connection is the $2$-Form $F_A \in \Omega^2(E,V)$ defined by
				$$F_A (x; v,w) := [ v - A_x(v), w - A_x(w) ] = [ v^{hor}, w^{hor}]^{vert}.$$
It measures the integrability of the horizontal distribution $H_A \subset TE$.\\

\textbf{Affine connections and vector bundles.}
A vector bundle is a fiber bundle $E$ whose fiber $F = V$ is a vector space and whose structure group $G \subset \text{GL}(V)$ is linear. In this case every fiber $E_z := \pi^{-1}(z)$ has a canonical structure of a vector space and we have well-defined maps 
		$$ \forall \lambda \in \mathbb{C}: \quad S_{\lambda}: E \rightarrow E, \qquad x \mapsto \lambda x  $$
		$$a : E \oplus E \rightarrow E, \qquad (x,y) \mapsto x + y.$$
A connection on $E$ is a connection $A \in \Omega^1(E,V)$ of the underlying  fiber bundle which is compatible with the linear structure on the fibers: Denote by $H_A \subset TE$ the horizontal distribution corresponding to $A$ and by $\tilde{H}_A \subset T(E\oplus E)$ the induced horizontal distribution consisting of pairs $(v,w) \in H\oplus H$ satisfying $d\pi(v) = d\pi(w)$. Then one requires
			\begin{align} \label{Preeq1} dS_{\lambda}(H) \subset H \quad \forall \lambda \in \mathbb{C} \qquad \text{and} \qquad da(\tilde{H}) \subset H. \end{align}
Alternatively, one can think of a connection as a covariant derivation
				$$d_A: \Omega^0(\Sigma, E) \stackrel{d}{\longrightarrow} \Omega^1(\Sigma, TE) \stackrel{A}{\longrightarrow} \Omega^1 (\Sigma, V) \cong \Omega^1(\Sigma, E)$$
where the last map comes from the canonical identification of the vertical bundle with the vector bundle itself. The linearity condition (\ref{Preeq1}) says precisely that this defines an affine connection.

\begin{Definition}
Let $E \rightarrow \Sigma$ be a complex vector bundle. An \textbf{affine connection} on $E$ is a linear operator $D: \Omega^0(\Sigma, E) \rightarrow \Omega^{1}(\Sigma, E)$ which satisfies the Leibniz rule
				$$D(fs) = df \otimes s + f \otimes D f$$
for all $f: \Sigma \rightarrow \mathbb{C}$ and $s \in \Omega^0(\Sigma,E)$.
\end{Definition}

We denote by $\mathcal{A}(E)$ the space of affine connections on $E$. Let $\psi_{\alpha}: E|_{U_{\alpha}} \rightarrow U_{\alpha} \times V$ be a local trivialization and denote for a local section $s : U_{\alpha} \rightarrow E$ with respect to this trivialization $s_{\alpha} := \text{pr}_2\circ \psi_{\alpha}$. Then an affine connection $D$ has the shape
					$$(Ds)_{\alpha} = d s_{\alpha} + A_{\alpha} s_{\alpha}$$
for some $A_{\alpha} \in \Omega^1(U_{\alpha}, \text{End}(V))$. These $A_{\alpha}$ are called connection potentials for the affine connection $D$. If all connection potentials take values in the Lie algebra $\mathfrak{g} \subset \text{End}(V)$ of the structure group $G \subset \text{GL}(V)$, then the affine connection $D$ is called a $G$-connection. We denote by $\mathcal{A}_G(E)$ the space of all $G$-connections on $E$.

An affine connection $D$ induces higher covariant derivations by the formula
						$$D: \Omega^k(\Sigma, E) \rightarrow \Omega^{k+1}(\Sigma, E), \qquad D(\tau\otimes s) = d \tau \otimes s + (-1)^{k} \tau \wedge Ds$$
for $\tau \in \Omega^k(\Sigma)$ and $s \in \Omega^0(\Sigma, E)$. The curvature $F_D \in \Omega^2(\Sigma, \text{End}(E))$ is the unique tensor satisfying
						$$ (D\circ D)s = F_D \cdot s$$
for all $s \in \Omega^0(\Sigma, E)$. It is the obstruction to $D^2 = 0$ and not directly related to the curvature of the horizontal distribution defined by $D$. It rather corresponds to curvature of the induced horizontal distribution in the frame bundle of $E$ as we shall see below.

\paragraph{Connections on principal bundles.} 
Let $G$ be a Lie group with Lie algebra $\mathfrak{g}$. A principal $G$ bundle over $\Sigma$ is a fiber bundle $\pi: P \rightarrow \Sigma$ together with a fiber preserving right action $P\times G \rightarrow P$ which is free and transitive on the fibers. In particular, the fibers are isomorphic to $G$ and using the right action we can always construct equivariant local trivializations of $P$. For $p \in P$ and $\xi \in \mathfrak{g}$ the infinitesimal action of $\xi$ is defined by
				$$p\xi := \left.\frac{d}{dt}\right|_{t=0} p\exp(t\xi) \in T_p P.$$
The collection of these tangent vectors defines the vertical subbundle 
			$$V = \ker{d\pi} = \{ p\xi \,| \, p \in P, \xi \in \mathfrak{g}\} \subset TP.$$ 
A connection on $P$ is an equivariant connection of the underlying fiber bundle and corresponds to an equivariant horizontal distribution $H \subset TP$ satisfying $TP = V \oplus H$.  Identifying $H$ with the projection $\Pi: TP = V\oplus H\rightarrow V$, we can describe such a connection by a $\mathfrak{g}$-valued $1$-form $A \in \Omega^1(P,\mathfrak{g})$ via the relation $\Pi_p(\hat{p}) = p A_p(\hat{p})$ for all $p \in P$ and $\hat{p} \in T_pP$. The connection $1$-Form $A$ satisfies the conditions
		\begin{align} \label{gaugeeq1} \quad A_p(p\xi) = \xi \quad\text{and} \quad  A_{pg} (\hat{p}g) = g^{-1} A_p(\hat{p}) g \end{align}
for all $g\in G$, $\xi \in \mathfrak{g}$, $p \in P$ and $\hat{p} \in T_pP$. Conversely, the kernel of any $A \in \Omega^1(P,\mathfrak{g})$ satisfying (\ref{gaugeeq1}) gives rise to an equivariant horizontal distribution $H \subset TP$ and we define by
				$$\mathcal{A}(P) := \{ A \in \Omega^1(P,\mathfrak{g})\,|\, \text{$A$ satisfies (\ref{gaugeeq1})} \}$$
the space of connections on $P$.

The curvature of a connection $A \in \mathcal{A}(P)$ is defined as
					$$F_A := dA + \frac{1}{2}[A \wedge A] \in \Omega^2(P, \mathfrak{g})$$				
where $[A \wedge A]$ is given by the usual formula for the exterior product with multiplication replaced by the Lie bracket. This curvature is linked to the curvature of the corresponding horizontal distribution by the relation
				$$ [X, Y]^{vert} =  [X^{hor}, Y^{hor}]^{vert} = pF_A(p; X,Y) $$
for $p \in P$ and $X,Y \in T_p P$.

\paragraph{Associated bundles.}
Let $P \rightarrow \Sigma$ be a principal $G$ bundle as above. A smooth manifold $F$ together with a representation $\rho: G \rightarrow \text{Diff}(F)$ gives rise to the associated fiber bundle $P\times_{\rho} F$ with fiber $F$ which is defined by
				$$P \times_{\rho} F := (P\times F) / G$$
where $G$ acts diagonally by $g(p,x) = (pg, \rho(g)^{-1}x)$. We denote the orbit of $(p,x) \in P\times F$ under this action by $[p,x]$. A connection $A \in \mathcal{A}(P)$ induces a connection on the fiber bundle $P\times_{\rho} F$, which is given by the image of the horizontal distribution under  $TP \subset TP\times TF \rightarrow T (P\times_{\rho} F)$.

Important examples arise from the action of $G$ on itself by inner automorphism and from the adjoint action of $G$ on its Lie algebra. We denote the associated bundles for these actions by
				$$\text{Ad}(P) := P \times_G G \qquad \text{and} \qquad \text{ad}(P) := P \times_{ad} \mathfrak{g}.$$
Note that the bundle $\text{Ad}(P)$ is a fiber bundle with fiber $G$ but not a principal bundle. The fibers of $\text{ad}(P)$ inherit from $\mathfrak{g}$ a well-defined Lie algebra structure. 

The difference $a := A_1 - A_2$ of two connection $1$-forms $A_1, A_2 \in \mathcal{A}(P)$ satisfies
			$$a_p(p\xi) = 0 \qquad \text{and} \qquad a_{pg}(\hat{p}g) = g^{-1}a_p(\hat{p})g$$ 
for all $p \in P$, $\hat{p} \in T_p P$, $\xi \in \mathfrak{g}$ and $g \in G$. Hence $a$ corresponds to a $\text{ad}(P)$-valued $1$-form $\bar{a}$ on $\Sigma$ by the formula $\bar{a}(\pi(p); d\pi(p)\hat{p}) = [p, a(p; \hat{p})]$. This describes $\mathcal{A}(P)$ as an affine space with underlying linear space $\Omega^1(\Sigma, \text{ad}(P))$ and with respect to any reference connection $A_0 \in \mathcal{A}(P)$  it holds
					$$\mathcal{A}(P) = \{ A_0 + a\, |\, a \in \Omega^1(\Sigma, \text{ad}(P) \}.$$
Similarly, the curvature $F_A$ of a connection $A$ is an equivariant and horizontal $2$-form on $P$ and can thus be identified with an element $F_A \in \Omega^2(\Sigma, \text{ad}(P))$.

Let $H$ be a Lie group and let $\tilde{\rho}: G \rightarrow H$ be a homomorphism of Lie groups. Then left-multiplication $\rho(g) := L_{\tilde{\rho}(g)} \in \text{Diff}(H)$ yields a representation of $G$ and the associated bundle $P_H := P\times_{\rho} H$ is a principal $H$ bundle. If $A \in \mathcal{A}(P)$, then $A$ induces a connection $\rho(A) \in  \mathcal{A}(P_H)$ by the formula
				$$\rho(A)( [p,h]; [\hat{p}, \hat{h}] ) := h^{-1} \hat{h} + h^{-1} \dot{\rho}( A(p; \hat{p}) ) h$$
where $\dot{\rho} := d\rho(\mathds{1}): \mathfrak{g} \rightarrow \mathfrak{h}$ denotes the induced homomorphism of Lie algebras. The curvature of the induced connection is given by
				$$F_{\rho(A)} = \dot{\rho}(F_A)$$
where $\dot{\rho}$ denotes the induced bundle map $\text{ad}(P) \rightarrow \text{ad}(P_H)$.

\paragraph{From principal bundles to vector bundles and back.}
Let $V$ be a vector space and let $\rho: G \hookrightarrow \text{GL}(V)$ be a faithful representation. The associated bundle $E := P\times_{\rho} V$ is then a vector bundle and the trivialization maps of $P$ yield a natural reduction of the structure group of $E$ to $G$. For a connection $A \in \mathcal{A}(P)$, the induced connection on $E$ is compatible with the linear structure and defines an affine $G$-connection in $\mathcal{A}_G(E)$. The bundles $\text{Aut}(E)$ and $\text{End}(E)$ can be described as associated bundles
				$$\text{Aut}(E) = P\times_{\text{Ad}(\rho)} \text{GL}(V) \qquad \text{and} \qquad \text{End}(E) = P\times_{\text{Ad}(\rho)} \text{End}(V)$$
where $\text{Ad}(\rho): G \rightarrow  \text{GL}(\text{End}(V))$ is defined as the composition of $\rho$ and the adjoint action of $\text{GL}(V)$ on $\text{End}(V)$. The induced map $\dot{\rho}: \mathfrak{g} \rightarrow \text{End}(V)$ provides an inclusion $\text{ad}(P) \rightarrow \text{End}(E)$ and with respect to this map holds
				$$F_{d_A} = \dot{\rho}(F_A)$$
for any connection $A \in \mathcal{A}(P)$.

Conversely, let $E \rightarrow \Sigma$ be a vector bundle with structure group $G \subset \text{GL}(n)$. The frame bundle of $E$ is defined by
		$$\text{Fr}(E) := \left\{ (z, e)\, | \, \text{$z \in \Sigma$, $e: V \rightarrow E_z$ such that $\text{pr}_2 \circ\psi_{\alpha}\circ e \in G$} \right\}$$
where $\psi_{\alpha} : E|_{U_{\alpha}} \rightarrow U_{\alpha}\times V$ is any trivialization of $E$ with $z \in U_{\alpha}$. It follows directly from the definition that $\text{Fr}(E)$ is a principal $G$ bundle. An affine $G$-connection $D \in \mathcal{A}_G(E)$ induces a connection $A_D \in \mathcal{A}(\text{Fr}(E))$ as follows: Let $\gamma: [0,1] \rightarrow \Sigma$ be a curve. We call $e \in \Omega^0( [0,1], \gamma^*\text{Fr}(E))$ a horizontal lift of $\gamma$ if for every $v \in V$ the section
				$$e_v \in \Omega^0([0,1], \gamma^* E), \qquad  e_v(t) := e(t) v \in E_{\gamma(t)} $$
satisfies $D_t (e_v) := D_{\dot{\gamma}(t)}( e_v(t)) = 0$. In a local trivialization this condition is equivalent to the ODE
							$$\dot{e}_{\alpha} + A_{\alpha}(\gamma)e_{\alpha} = 0.$$
This shows that horizontal lifts exist when the connection potentials $A_{\alpha}$ take values in $\mathfrak{g}$. The tangent vector along horizontal lifts trace out an equivariant horizontal distribution in $\text{Fr}(E)$ and hence determine a connection $A \in \mathcal{A}(\text{Fr}(E))$.

The frame bundle construction is inverse to the construction of associated bundles in the sense that 
			$$\text{Fr}(P\times_G V) \cong V \qquad \text{and} \qquad \text{Fr}(E)\times_G V \cong E$$
whenever $G \subset \text{GL}(V)$. This also provides a one-to-one correspondence between $\mathcal{A}(P)$ and $\mathcal{A}_G(E)$.

\paragraph{The Gauge group.}
The Gauge group of a principal $G$ bundle $P$ is defined as
				$$\mathcal{G}(P) := \Omega^0(\Sigma, \text{Ad}(P)).$$
This group is isomorphic to the group $\text{Aut}(P)$ of fiber preserving equivariant automorphism of $P$ under the map
				$$\psi: \Omega^0(\Sigma, \text{Ad}(P)) \cong \text{Aut}(P),\qquad  \psi_g(p) := pg(p).$$
It is useful think of $\mathcal{G}(P)$ as an infinite dimensional Lie group with Lie algebra
			$$ \text{Lie}(\mathcal{G}(P)) = \Omega^0(\Sigma,\text{ad}(P))$$
where all Lie theoretic operations are performed fiberwise. 
The Gauge group acts naturally on the space of connections via pull back
				$$ g(A) := \psi_{g^{-1}}^*A = -(dg)g^{-1} + g A g^{-1}.$$

The Gauge group of a vector bundle $E$ with structure group $G$ is the group
				$$\mathcal{G}(E) := \Omega^0(\Sigma, G(E)) \subset \Omega^0(\Sigma, \text{GL}(E))$$
which consists of all automorphisms of $E$ taking values in $G$ in any trivialization. We think of $\mathcal{G}(E)$ again as Lie group with Lie algebra $\Omega^0(\Sigma,\mathfrak{g}(E))$. The Gauge group acts naturally on the space of affine $G$-connection $\mathcal{A}_G(E)$ via pullback
				$$(g^{-1})^* D = g \circ D \circ g^{-1}.$$
This action is more explicitly described in terms of the connection potential by
				$$ (gA)_{\alpha} = - dg_{\alpha} g_{\alpha}^{-1} + g_{\alpha} A_{\alpha} g_{\alpha}^{-1} $$
where $g_{\alpha} := (\text{pr}_2\circ \psi_{\alpha})_*g: U_{\alpha} \rightarrow G$.

Suppose that $\rho: G \hookrightarrow \text{GL}(V)$ is a faithful representation and $E := P\times_{\rho} V$ is an associated vector bundle. Then $\rho$ induces an isomorphism $\text{Ad}(P) \cong G(E)$ and hence $\mathcal{G}(P) \cong \mathcal{G}(E)$. The derivative $\dot{\rho} := d\rho(\mathds{1}): \mathfrak{g} \hookrightarrow \text{End}(V)$ yields an isomorphism of $\text{ad}(P) \cong \mathfrak{g}(E)$ and hence an identification of the Lie algebras of $\mathcal{G}(P)$ and $\mathcal{G}(E)$. From the naturality of the gauge action it is clear that the identification $\mathcal{A}(P) \cong \mathcal{A}_G(E)$ is equivariant with respect to the action of $\mathcal{G}(P)$ and $\mathcal{G}(E)$.

\paragraph{The moment map picture.}

Fix an invariant inner product $\langle\cdot , \cdot \rangle$ on $\mathfrak{g}$. This induces an inner product on the fibers of $\text{ad}(P)$ and hence an invariant inner product on $\text{Lie}(\mathcal{G}(P)) = \Omega^0(\Sigma, \text{ad}(P))$ by the formula
						$$\langle \xi, \eta \rangle := \int_{\Sigma} \langle \xi, \eta \rangle \, dvol_{\Sigma}.$$
This provides a natural hermitian structure on the space $\mathcal{A}(P)$ as follows. Since $\mathcal{A}(P)$ is an affine space, it suffices to define the hermitian structure on the underlying linear space $\Omega^1(\Sigma, \text{ad}(P))$. For $a,b \in \Omega^1(\Sigma, \text{ad}(P))$ we define
			$$\omega_{\mathcal{A}}(a,b) := \int_{\Sigma} \langle a\wedge b \rangle, \qquad  \langle a , b \rangle := \int_{\Sigma} \langle a \wedge *b \rangle, \qquad J_{\mathcal{A}}a := *a = - a\circ j_{\Sigma}.$$
The following observation is due to Atiyah and Bott \cite{AtBott:YangMillsEq}.

\begin{Lemma} \label{PreLemma1}
Let $G$ be a compact connected Lie group. The action of the Gauge group is hamiltonian with moment map $\mu(A) := *F_A$. More explicitly, the infinitesimal action of $\xi \in \Omega^0(\Sigma, \text{ad}(P))$ on $A \in \mathcal{A}(P)$ is given by
			$$L_A \xi := \left.\frac{d}{dt}\right|_{t=0} \exp(t\xi) (A) = - d_A \xi$$
and satisfies
			$$d \langle \mu(A), \xi \rangle = \omega_{\mathcal{A}}(L_A \xi, \cdot).$$
\end{Lemma}
 
\begin{proof}
Let $\xi \in \Omega^0(\Sigma, \text{ad}(P))$ be given and think of it as an equivariant map $\xi: P \rightarrow \mathfrak{g}$. We then compute
				\begin{align*} \left.\frac{d}{dt}\right|_{t=0} \exp(t\xi) (A)  
										&= \left.\frac{d}{dt}\right|_{t=0} d\exp(t\xi)^{-1}\exp(t\xi) + \exp(t\xi) A \exp(- t\xi) \\
										&= \left.\frac{d}{dt}\right|_{t=0} \left(d\exp(t\xi)^{-1} \right) + \left[ \xi, A\right] \\
										&= - d \xi - \left[ A, \xi \right]
				\end{align*}
The last expression agrees with $ -d\xi$ along horizontal vectors in $P$ and vanishes along vertical vectors. Hence it coincides with $-d_A \xi$ for the induced affine connection $d_A$ on $\text{ad}(P)$ and this proves the formula for the infinitesimal action.

From the formula
			$$F_{A+a} = F_A + d_A a + \frac{1}{2} [a\wedge a]$$
we see that the variation of $F_A$ in the direction $a \in \Omega^1(\Sigma, \text{ad}(P))$ is given by $d_A a$. This yields
	$$\langle d\mu(A)[a]; \xi \rangle = \int_{\Sigma} \langle d_Aa , \xi \rangle =  \int_{\Sigma} \langle a \wedge d_A \xi \rangle  = \omega_{\mathcal{A}}(L_A \xi, a).$$
Here we used integration by part in the penultimate step and the formula 
				$$d \langle a, \xi \rangle = \langle d_A a ,\xi \rangle - \langle a \wedge d_A \xi \rangle$$
which follows from the $G$-invariance of the inner product.
\end{proof}

\subsubsection{The complexified gauge action}

Let $G$ be a compact connected Lie group and let $P \rightarrow \Sigma$ be a principal $G$ bundle. We denote by $G^c$ the complexification of $G$ and call $P^c := P\times_G G^c$ the complexification of $P$. The complexified gauge group of $P$ is defined as
				$$\mathcal{G}^c(P) := \mathcal{G}(P^c).$$
One can think of elements in $\mathcal{G}^c(P)$ as $G$-equivariant maps from $P$ to $G^c$. The Lie algebra bundle $\text{ad}(P^c)$ is the complexification of the bundle $\text{ad}(P)$ and since all Lie theoretic operations on the gauge group are defined fiberwise, it is reasonable to think of $\mathcal{G}^c(P)$ as the complexification of $\mathcal{G}(P)$. By the Peter-Weyl theorem, $G$ admits a faithful representation $G \hookrightarrow U(n)$. Identifying $G$ with its image in $U(n)$, we can describe its complexification $G^c \subset \text{GL}(n)$ explicitly as the image of $G\times \mathfrak{g}$ under the diffeomorphism
					$$U(n)\times \mathfrak{u}(n) \rightarrow \text{GL}(n), \qquad (u,\eta) \mapsto u \exp(\textbf{i}\eta).$$ 
In terms of the associated bundle $E := P\times_G \mathbb{C}^n$ the complexification of the gauge group is then given by
							$$\mathcal{G}^c(E) = \Omega^0(\Sigma, G^c(E)).$$
The goal of this section is to explain that the action of $\mathcal{G}(P)$ on the space of connections $\mathcal{A}(P)$ extends naturally to an holomorphic action of $\mathcal{G}^c(P)$. 

\begin{Proposition} \label{CpxGaugeAction}
There exists a natural action of $\mathcal{G}^c(P)$ on $\mathcal{A}(P)$ whose infinitesimal action satisfies
				\begin{align} \label{gaugeeq2} L_A (\xi + \textbf{i}\eta) = L_A \xi + * L_A\eta = - d_A \xi - *d_A \eta\end{align}
for all $\xi,\eta \in \Omega^0(\Sigma,\text{ad}(P))$ and $A \in \mathcal{A}(P)$.
\end{Proposition}

\begin{proof} See page \pageref{proofCpxGaugeAction}. \end{proof}

\paragraph{Holomorphic principal bundles.} An almost complex structure $J$ on a manifold $M$ is an endomorphism $J \in \text{End}(TM)$ satisfying $J^2 = -\mathds{1}$. It is called an integrable or holomorphic structure if it endows $M$ with the structure of a complex manifold. A holomorphic structure on the principal bundle $P^c = P\times_G G^c$ is an almost complex structure $J \in \text{End}(TP^c)$ of the total space, which is $G^c$ invariant and coincides with the canonical complex structure on the vertical subbundle, i.e. $J(p\zeta) = p(\textbf{i}\zeta)$ for any $p \in P^c$ and $\zeta \in \mathfrak{g}^c$. We denote by $\mathcal{J}(P^c)$ the space of all holomorphic structures on $P^c$. The next Lemma justifies this notation.

\begin{Lemma}
Every $J\in\mathcal{J}(P^c)$ is integrable.
\end{Lemma}

\begin{proof}
The Newlander-Nirenberg theorem states that an almost complex structure $J$ on a manifold $M$ is integrable if and only if the Nijenhuis-tensor $N_J : TM\otimes TM \rightarrow TM$ given by
				$$N_J(v,w) := [v,w]  + J[Jv,w] + J[v,Jw] - [Jv,Jw]$$
vanishes. We apply this to $M = P^c$. If $v,w \in T_p^{vert}(P^c)$ are both in the vertical bundle, we have $N_J(v,w) = 0$ as the fiber is a complex manifold. If $v \in T_p^{vert}(P^c)$ and $w \in T^{hor}_p (P^c)$ the Lie bracket $[v,w] = \mathcal{L}_v(w)$ vanishes, since the horizontal distribution is equivariant. In particular $N_J(v,w) = 0$ as all four terms vanish separately. Let finally $v,w \in T^{hor}_p(P^c)$ be horizontal vectors. We may assume that $p \in P$ and denote by $\bar{v} := d\pi(p)v$ and $\bar{w} := d\pi(p)w$ the projections onto $T\Sigma$. By definition of the curvature, we obtain the vertical component of the Nijenhuis tensor by
				\begin{align*}
					A_p(N_J(v,w)) &= F_A(\bar{v},\bar{w}) + \textbf{i} F_A(j_{\Sigma}\bar{v}, \bar{w}) + \textbf{i} F_A(\bar{v}, j_{\Sigma}\bar{w}) - F_A(j_{\Sigma}\bar{v}, j_{\Sigma}\bar{w}) \\
												&= 4 F_A^{0,2}(\bar{v},\bar{w}) = 0. 
				\end{align*}
In the last step we use that $\Sigma$ is a complex one-dimensional manifold and thus $\Omega^{0,2}(\Sigma) = 0$. The horizontal part of $N_J(v,w)$ gets identified under $d\pi(p)$ with $N_J(\bar{v},\bar{w})$ and vanishes as $\Sigma$ is a complex manifold. This completes the proof of $N_J = 0$.
\end{proof}

As a consequence, every holomorphic principal bundles admits holomorphic local trivializations with holomorphic transition maps. The next Lemma is due to Singer \cite{Singer:1959}.

\begin{Lemma} \label{PreLem2}
There exists a one to one correspondence between connections $A \in \mathcal{A}(P)$ and holomorphic structures $J \in \mathcal{J}(P^c)$.
\end{Lemma}

\begin{proof}
A connection $A \in \mathcal{A}(P)$ induces a connection on $P^c$ and thus determines for every $p \in P^c$ a splitting $T_p (P^c) = T_p^{hor}(P^c) \oplus T_p^{vert} (P^c)$. The vertical part is isomorphic to $\mathfrak{g}^{\mathbb{C}}$ and has a canonical complex structure. The differential of the projection $\pi: P^c \rightarrow \Sigma$ restricts to an isomorphism $d\phi(p): T_p^{hor} (P^c) \rightarrow T_{\pi(p)}\Sigma$ and induces a complex structure on $T_p^{hor}$.

Conversely, let $J \in \mathcal{J}(P^c)$ be given and think of $P \subset P^c$ as a subbundle. For $p \in P$ we define $H_p := T_p P \cap J_p(T_p P)$ and claim that $T_p P + J_p(T_p P) = T_p(P^c)$. Indeed, since $T^{vert}_p P \cong \mathfrak{g}$, the sum clearly contains the vertical fiber $T^{vert}_p (P^c)\cong \mathfrak{g}^c$ and $d\pi(p)$ maps $T_p P$ already onto $T_{\pi(p)}\Sigma$. It is immediate from the construction that $H_p$ is invariant under $J_p$ an defines a (real) two dimensional complement of $T_p^{vert} (P^c)$ in $T_p (P^c)$. As $p$ varies over $P$ we obtain an equivariant distribution along $P$ and hence a connection $A \in \mathcal{A}(P)$. 
\end{proof}

Let $A \in \mathcal{A}(P)$, $g \in \mathcal{G}(P)$ and let $J_A \in \mathcal{J}(P^c)$ be the holomorphic structure induced by $A$. Then $g(A)$ induces the holomorphic structure $(\psi_{g^{-1}})^* J_A$, since the construction above is clearly functorial. The action of $\mathcal{G}(P)$ on $\mathcal{J}(P^c)$ has a natural extension to the complexified gauge group via
					$$\mathcal{G}^c(P)\times \mathcal{J}(P^c) \rightarrow \mathcal{J}(P^c), \qquad g(J) := (\psi_{g^{-1}})^* J$$
where $\psi_{g^{-1}} \in \text{Aut}(P^c)$ is the automorphism corresponding to $g^{-1}$. Using the identification of $\mathcal{J}(P^c)$ with $\mathcal{A}(P)$ this yields the desired action of $\mathcal{G}^c(P)$ on $\mathcal{A}(P)$ and the quotient $\mathcal{A}(P)/\mathcal{G}^c(P)$ parametrizes the isomorphism classes of holomorphic structures on $P^c$.

\paragraph{Holomorphic vector bundles.} We consider the special case $G = U(n)$ and denote by $E:= P\times_{U(n)} \mathbb{C}^n$ the associated vector bundle. A holomorphic structure on $E$ is an almost complex structure $J\in \text{End}(TE)$ of the total space which restricts to the linear complex structure on the fibers. Similarly as in the case of principal bundles, one shows that every such structure is indeed integrable and that every holomorphic vector bundle admits holomorphic trivializations. It is then easy to see that every holomorphic vector bundle $E$ carries a natural operator
					$$\bar{\partial}_E: \Omega^0(\Sigma, E) \rightarrow \Omega^{0,1}(\Sigma, E)$$
which in any holomorphic trivialization agrees with the usual $\bar{\partial}$ operator on $\mathbb{C}^n$. This operator is a particular Cauchy-Riemann operator on $E$.

\begin{Definition}
Let $E \rightarrow \Sigma$ be a complex vector bundle. A Cauchy Riemann operator on $E$ is a linear operator
				$$D'': \Omega^0(\Sigma, E) \rightarrow \Omega^{0,1}(\Sigma, E)$$
which satisfies the Leibniz rule
				$$D''(fs) = \bar{\partial}f \otimes s + f \otimes D'' f$$
for all $f: \Sigma \rightarrow \mathbb{C}$ and $s \in \Omega^0(\Sigma,E)$.
\end{Definition}

The converse is also true: Every Cauchy-Riemann operator determines a holomorphic structure on the complex bundle $E$, whose local holomorphic sections are solutions of the Cauchy-Riemann equation $D'' s = 0$. This is another instance of the Newlander-Nirenberg theorem. In the case of Riemann surfaces a simpler proof of this result is given by Atiyah and Bott (\cite{AtBott:YangMillsEq}, Section 5).

Note that the associated vector bundle $E$ carries a canonical hermitian metric, which in any trivialization coincides with the standard hermitian metric on $\mathbb{C}^n$. We claim that there is a one to one correspondence between unitary connections on $E$ and Cauchy-Riemann operators. For an unitary connection $D$ we obtain a Cauchy-Riemann operator by the formula
			$$D'' s := \left( D s\right)^{0,1} := \frac{1}{2}\left( D s+ \textbf{i} (D s)\circ j_{\Sigma}\right) = \frac{1}{2}\left( D s - \textbf{i} *(D s)\right). $$
To show that this correspondence is bijective, it suffices to examine this correspondence locally. In an unitary trivialization $\psi: E|_{U} \rightarrow U \times \mathbb{C}^n$ the connection $D$ can be described in terms of a $1$-form $A \in \Omega^1(U,\mathfrak{u}(n))$ such that
			$$D s := ds + A s, \qquad D'' s := \bar{\partial} s + A^{0,1} s$$
holds for any section $s \in \Omega^0(U, \mathbb{C}^n)$ with $A^{0,1} := \frac{1}{2}(A + \textbf{i}A\circ j_{\Sigma})$. In particular, we recover $A$ as twice the skew-hermitian part of $A^{0,1}$ and therefore it is uniquely determined by $A^{0,1}$. Conversely, any Cauchy Riemann operator $D''$ is given in this local trivialization by
			$$D'' s := \bar{\partial} s + B s$$
for some $B \in \Omega^{0,1}(U,\mathfrak{gl}(n))$. Since $B$ satisfies $B(j_{\Sigma} v) = - \textbf{i} B(v)$ for any tangent vector $v\in T\Sigma|_{U}$, the skew-hermitian and hermitian part of $B$ interchange if we compose $B$ with $j_{\Sigma}$. This shows that $B$ has the form $B = \frac{1}{2}(A +\textbf{i}A\circ j_{\Sigma})$ for some $A \in \Omega^1(U,\mathfrak{u}(n))$ and this proves the claim.

On the level of Cauchy-Riemann operators the complexified Gauge group $\mathcal{G}^c(E) = \Omega^0(\Sigma, \text{GL}(E))$ acts naturally via
							$$g( \bar{\partial}_A) := g \circ \bar{\partial}_A \circ g^{-1} = \bar{\partial}_A - \bar{\partial}_A(g) g^{-1}.$$
The next Lemma summarizes the discussion above and provides an explicit formulas for this action on $\mathcal{A}(E)$.

\begin{Lemma} \label{PreLem1}
Let $E \rightarrow \Sigma$ be a complex vector bundle. 
	\begin{enumerate}
		\item	For every holomorphic structure $\bar{\partial}_E$ and hermitian metric $H$  exists a unique connection
			$$D := D(\bar{\partial}_E, H) =: D' + D'' \in \mathcal{A}^{1,0}(E)\oplus \mathcal{A}^{0,1}(E)$$
such that $D$ is unitary with respect to $H$ and $D'' = \bar{\partial}_E$.
		\item	Let $g \in \Omega^0(\Sigma, \text{GL}(E))$ and denote $h := g^*g$ (with respect  to $H$). Then
							$$D(g(\bar{\partial}_E), H) = g\left(D + h^{-1} D' (h) \right) g^{-1}$$
							$$F(g(\bar{\partial}_E), H) = g\left(F + D''(h^{-1} D'(h)) \right) g^{-1}.$$			
	\end{enumerate}
\end{Lemma}

\begin{proof}
For the first part, note that there is a one to one correspondence between hermitian metrics $H$ and reductions of the structure group of $E$ to $U(n)$: Using the Gram-Schmidt process we can always find local trivializations which identify $H$ with the standard hermitian product on $\mathbb{C}^n$ and the transition map between such trivializations are clearly unitary. The second part follows from the formula
			$$D(g(\bar{\partial}_E), H) = g(D) = g\circ D'' \circ g^{-1} + (g^{-1})^* \circ D' \circ g^*$$
and $F = D\circ D$.
\end{proof}

\begin{Remark}
Consider the general case and assume that $G \subset U(n)$ is a compact connected subgroup. The structure group of $E$ is then contained in $G$ and the explicit formula in the Lemma above shows that the subspace $\mathcal{A}_G(E)$ of $G$-connections is preserved by the action of $\mathcal{G}^c(E) = \Omega^0(\Sigma, G^c(E))$. Since holomorphic structures on $E$ and its frame determine one another, it is clear that this action corresponds to the action described on holomorphic principal bundles above.
\end{Remark}

We may now deduce the formula for the infinitesimal action (\ref{gaugeeq2}).

\begin{proof}[Proof of Proposition \ref{CpxGaugeAction}.] \label{proofCpxGaugeAction}  
 As in Lemma \ref{PreLemma1} one calculates
				$$\left. \frac{d}{dt} \right|_{t=0} \exp(t\zeta)(\bar{\partial}_A) = - \bar{\partial}_A \zeta$$
for $\zeta \in \Omega^0(\Sigma, \mathfrak{g}^c(E))$. Write $\zeta = \xi + \textbf{i}\eta$ with $\xi, \eta \in \Omega^0(\Sigma, \mathfrak{g}(E))$ and use the formula $\bar{\partial}_A (\textbf{i}\eta) = * \bar{\partial}_A \eta$ to deduce
		\begin{align*}
					L_A \zeta &= - \bar{\partial}_A \zeta + (\bar{\partial}_A \zeta)^* 
										= - (\bar{\partial}_A \xi - (\bar{\partial}_A \xi)^*) - *(\bar{\partial}_A \eta - (\bar{\partial}_A \eta)^*) \\
										&= - d_A\xi - * d_A \eta = L_A \xi + * L_A \eta.
		\end{align*}
\end{proof}

\subsubsection{Regularity assumptions}
Let $G$ be a compact connected Lie group and let $P \rightarrow \Sigma$ be a principal $G$ bundle. We shall always consider connections of Sobolev class $W^{1,2}$ and gauge transformations of Sobolev class $W^{2,2}$. More precisely, the space of $W^{1,2}$ connections on $P$ is defined with respect to some smooth reference connection $A_0$ as
				$$\mathcal{A}(P) := \{A_0 + a\,|\, a \in W^{1,2}(\Sigma, T^*\Sigma\otimes \text{ad}(P))$$
and the $W^{2,2}$ completion of the gauge group and its complexification are
				$$\mathcal{G}(P) := W^{2,2}(\Sigma, \text{Ad}(P)), \qquad \mathcal{G}^c(P) := W^{2,2}(\Sigma, \text{Ad}(P^c)).$$
We use the same notation as for the smooth groups, since all the results from the previous section carry over. In particular, the action of the gauge group and its complexification extend smoothly over these Sobolev completions, since $W^{2,2} \hookrightarrow C^0$ is in the good range of the Sobolev embedding. A connection still determines a holomorphic structure up to isomorphism due to the following regularity result.

\begin{Lemma}\label{LemmaSmoothGauge}
For every $W^{1,2}$ connection $A \in \mathcal{A}(P)$ exists a complex $W^{2,2}$ gauge transformation $g \in \mathcal{G}^c(P)$ such that $g(A)$ is smooth.
\end{Lemma}

\begin{proof}
This is Lemma 14.8 in \cite{AtBott:YangMillsEq}. By Proposition \ref{CpxGaugeAction}, the infinitesimal action of the complex gauge group is given by
			$$L_A: W^{2,2}(\Sigma, \text{ad}(P^c)) \rightarrow W^{1,2}(\Sigma, T^*\Sigma\otimes \text{ad}(P))$$
							$$L_A(\xi + \textbf{i}\eta) = -d_A \xi - *d_A \eta$$
For any smooth reference connection $A_0$, this is a compact perturbation of $L_{A_0}$ which is a Fredholm operator. Hence $L_A$ is also Fredholm and in particular its cokernel is finite dimensional.

It follows from the implicit function theorem in Banach spaces that we can choose a finite dimensional slice $N$ orthogonal to the $\mathcal{G}^c$-orbit through $A$. Say $\text{dim}(N) = r$ and fix $r+1$ connections $B_0, \ldots, B_{r} \in N$ which span an $r$-simplex containing $A$ in its interior. A small perturbation of the vertices yields smooth connections $\tilde{B}_0, \ldots, \tilde{B}_r$ and the simplex spanned by these connections will still intersect the orbit $\mathcal{G}^c(A)$. This intersection point yields a smooth connection in the $\mathcal{G}^c$ orbit of $A$.
\end{proof}

\subsection{Parabolic subgroups}
Let $G$ be a compact connected Lie group with Lie algebra $\mathfrak{g}$ and denote its complexification by $G^c$. Fix an invariant inner product on $\mathfrak{g}$. This induces a (real valued) inner product on $\mathfrak{g}^c = \mathfrak{g}\oplus \textbf{i}\mathfrak{g}$ where we define both factors to be orthogonal. We define parabolic subgroups of $G^c$ first by using toral generators of $\mathfrak{g}^c$. Then we recall briefly the root space decomposition of reductive Lie algebras and give an alternative intrinsic definition of parabolic subgroups. The first definition occures naturally in the geometric description of the weights in Chapter 5. The intrinsic version turns out to be useful in the proof of Proposition \ref{PropStabChar} which relates the algebraic notion of stability with the weights.

\subsubsection{Toral generators}
An element $\zeta \in \mathfrak{g}^c$ is called a \textbf{toral generator} if 
				$$T_{\zeta} := \overline{ \{\exp(t\zeta)\,| t \in \mathbb{R} \}} \subset G^c$$ 
is a compact torus. We denote by $\mathcal{T}^c$ the set of toral generators. Certainly $\mathfrak{g} \subset \mathcal{T}^c$. Since any maximal compact subgroup of $G^c$ is conjugated to $G$, for every $\zeta \in \mathcal{T}^c$ exists $g \in G^c$ such that $g^{-1} T_{\zeta} g^{-1} \subset G$. The relation $g T_{\zeta} g^{-1} = T_{g \zeta g^{-1}}$ then yields $g\zeta g^{-1} \in \mathfrak{g}$ and hence
			$$\mathcal{T}^c = \text{Ad}(G^c)(\mathfrak{g}) = \{ g\xi g^{-1}\,|\, g \in G^c,\, \xi \in \mathfrak{g}\}.$$

\begin{Definition} \label{DefnToralGen} A \textbf{parabolic subgroup} of $G^c$ is a subgroup of the form
			$$Q(\zeta) := \{ g \in G^c\,|\, \text{the limit $\lim_{t\rightarrow \infty} e^{\textbf{i}t\zeta} g e^{-\textbf{i}t\zeta}$ exists in $G^c$} \}$$
for some $\zeta \in \mathcal{T}^c$. The \textbf{Levi subgroup} of $Q(\zeta)$ is defined by
			$$L(\zeta) := \{ g \in G^c\,|\, e^{\textbf{i}\zeta} g e^{-\textbf{i}\zeta} = g \}.$$
\end{Definition}

\begin{Remark}
We consider $G^c = Q(0)$ as parabolic subgroup of itself. 
\end{Remark}

\begin{Lemma} \label{pblLem1}
Consider the setting described above and let $\zeta \in \mathcal{T}^c$.
		\begin{enumerate}
			\item $Q(\zeta)$ is a closed connected Lie subgroup of $G^c$ with Lie algebra
						$$\mathfrak{q}(\zeta) := \{ \rho \in \mathfrak{g}^c\,|\, \text{the limit $\lim_{t\rightarrow \infty} e^{\textbf{i}t\zeta} \rho e^{-\textbf{i}t\zeta}$ exists in $\mathfrak{g}^c$} \}.$$
			\item $L(\zeta)$ is a closed connected Lie subgroup of $G^c$ with Lie algebra
						$$\mathfrak{l}(\zeta) := \{ \rho \in \mathfrak{g}^c\,|\,  e^{\textbf{i}t\zeta} \rho e^{-\textbf{i}t\zeta} = \rho \}$$
			\item $L(\zeta)$ is a maximal reductive subgroup of $Q(\zeta)$.
			
			\item $Q(\zeta) = G^c$ if and only if $\zeta$ is contained in the center of $\mathfrak{g}^c$.
		\end{enumerate}
\end{Lemma}

\begin{proof}
Since $Q(g\zeta g^{-1}) = g Q(\zeta) g^{-1}$ and $L(g \zeta g^{-1}) = g L(\zeta) g^{-1}$, we may assume $\zeta = \xi \in \mathfrak{g}$. By the Peter-Weyl theorem, there exists a faithful representation $G \hookrightarrow U(n)$ and we may identify $G$ with a closed subgroup of $U(n)$. Then $\textbf{i}\xi$ yields a hermitian endomorphism of $\mathbb{C}^n$ which is diagonalizable with real eigenvalues $\lambda_1 < \cdots < \lambda_r$. Denote the eigenspace corresponding to $\lambda_j$ by $V_j$. They yield an orthogonal decomposition 
				$$\mathbb{C}^n = V_1\oplus \cdots \oplus V_r.$$
In this decomposition we can write $g \in G^c \subset \text{GL}(n,\mathbb{C})$ as
				$$ g = \begin{pmatrix} g_{11} & g_{12} & \cdots & g_{1r}\\
															 g_{21} & g_{22} & \cdots & g_{2r}\\
															\vdots  & \vdots & \ddots & \vdots \\
															 g_{r1} & g_{r2} & \cdots & g_{rr} \end{pmatrix} $$
with $g_{ij} \in \text{Hom}(V_j, V_i)$. Then
		$$e^{\textbf{i}t\xi} g e^{-\textbf{i}t\xi} = \begin{pmatrix} g_{11} & e^{(\lambda_1 - \lambda_2)t} g_{12} & \cdots & e^{(\lambda_1 - \lambda_r)t} g_{1r}\\
																																 e^{(\lambda_2 - \lambda_1)t} g_{21} & g_{22} & \cdots & e^{(\lambda_2 - \lambda_r)t} g_{2r}\\
																																		\vdots  & \vdots & \ddots & \vdots \\
																																 e^{(\lambda_r - \lambda_1)t} g_{r1} & e^{(\lambda_r - \lambda_2)t} g_{r2} & \cdots & g_{rr} 
																								\end{pmatrix} .$$
Thus $g \in Q(\xi)$ if and only if $g$ is upper triangular (i.e. $g_{ij} = 0$ for $i > j$) and $g \in L(\xi)$ if and only if $g$ is block diagonal (i.e. $g_{ij} = 0$ for $i \neq j$). This shows that $L(\xi)$ and $Q(\xi)$ are closed subgroups of $G^c$ and the formulas for $\mathfrak{l}(\xi)$ and $\mathfrak{q}(\xi)$ are immediate.

As the spaces $V_j$ are pairwise orthogonal, the intersection $G \cap Q(\xi)$ consists of block diagonal matrices and hence agrees with the centralizer of the torus $T_{\xi}$ in $G$. Since the centralizers of tori in compact groups are connected (see \cite{LiegpKnapp} Corollary 4.51) we conclude that $G\cap Q(\xi)$ is connected. Since $L(\xi)$ is the complexification of $G\cap Q(\xi)$ it is connected and reductive. Moreover $Q(\xi)/L(\xi)$ can be identified with the unipotent matrices in $Q(\xi)$ and hence $L(\xi)$ is a maximal reductive subgroup of $Q(\xi)$. We observe that
			$$Q(\xi) \rightarrow L(\xi), \qquad g \mapsto \lim_{t\rightarrow \infty} e^{\textbf{i}t\xi} g e^{-\textbf{i}t\xi}$$
defines a continuous retraction of $Q(\xi)$ onto $L(\xi)$ and hence $Q(\xi)$ is connected.

Finally, since $G^c$ is reductive, we have $G^c = Q(\xi)$ if and only if $G^c = L(\xi)$. The later is clearly equivalent to $\xi \in Z(\mathfrak{g})$.
\end{proof}

\subsubsection{The root-space decomposition}
We recall the necessary background on Lie theory briefly and refer to \cite{LiegpKnapp} for the proofs. Note that the discussion remains valid for any $G$-invariant inner product on $\mathfrak{g}$, which does not need to be the negative Killing form.

\paragraph{Reductive Lie groups.}
Using the invariant inner product on $\mathfrak{g}$, it is easy to show that the adjoint action of $\mathfrak{g}$ on itself is completely reducible. This yields an orthogonal decomposition
				$$\mathfrak{g} = \mathfrak{z}\oplus [\mathfrak{g},\mathfrak{g}]$$
where $\mathfrak{z}$ denotes the center of $\mathfrak{g}$ and the commutator $[\mathfrak{g},\mathfrak{g}]$ is a direct sum of simple ideals and hence a semisimple Lie algebra. The same decomposition is valid for the complexification. To see this extend the inner product on $\mathfrak{g}$ to a non-degenerated $\mathbb{C}$-bilinear form  $B: \mathfrak{g}^c\times \mathfrak{g}^c \rightarrow \mathbb{C}$ by
		$$B(\xi_1 + \textbf{i}\eta_1; \xi_2 + \textbf{i}\eta_2) = \langle \xi_1, \xi_2 \rangle - \langle \eta_1, \eta_2 \rangle + \textbf{i}(\langle \xi_1, \eta_2 \rangle + \langle \eta_1, \xi_2 \rangle).$$ 
This bilinear form is nondegenerate and $G^c$-invariant. Moreover, the $B$-orthogonal complement of a complex subspace $W \subset \mathfrak{g}^c$ is a $G^c$-invariant complement and the same argument as above yields the decomposition
				$$\mathfrak{g}^c = \mathfrak{z}^c \oplus [\mathfrak{g}^c,\mathfrak{g}^c].$$

\paragraph{Root space decomposition.} Fix a maximal torus $T \subset G$ with Lie algebra $\mathfrak{t}$ and decompose it orthogonally as $\mathfrak{t} = \mathfrak{z} \oplus \mathfrak{t}_0$. A nonzero imaginary valued real linear map
				$$\alpha = \textbf{i}a : \mathfrak{t}_0 \rightarrow \textbf{i}\mathbb{R}, \qquad a \in \text{Hom}(\mathfrak{t}_0, \mathbb{R})$$
is called a \textbf{root} of $G$ with respect to $T$ if there exists $e_{\alpha} \in [\mathfrak{g}^c, \mathfrak{g}^c]$ satisfying
				$$[t, e_{\alpha}] = \alpha(t) e_{\alpha} \qquad \text{for all $t\in\mathfrak{t}_0$}.$$
The element $e_{\alpha}$ is uniquely determined by $\alpha$ up to scaling. We denote by $\mathfrak{g}_{\alpha} := \mathbb{C}\cdot e_{\alpha}$ the one dimensional root space corresponding to $\alpha$ and denote by $R$ the set of all roots (relative to $T$).  The \textbf{root space decomposition} of $\mathfrak{g}^c$ is the vector space decomposition
				$$\mathfrak{g}^c = \mathfrak{z}^c \oplus \mathfrak{t}_0^c \oplus \bigoplus_{\alpha \in R} \mathfrak{g}_{\alpha}.$$
For a proof see \cite{LiegpKnapp} Chapters II.1-4 and IV.5.

\begin{Lemma}
Denote $\mathfrak{g}_0 := \mathfrak{t}_0^c$. 
	\begin{enumerate}
			\item For $\alpha, \beta \in R\cup\{0\}$ the Lie bracket satisfies the relation
				$$[\mathfrak{g}_{\alpha}, \mathfrak{g}_{\beta}] \subset \mathfrak{g}_{\alpha + \beta}$$
				where the right-hand side is defined to be zero when $\alpha + \beta \notin R \cup \{0\}$. 
			\item For $\alpha, \beta \in R\cup\{0\}$ with $\alpha \neq - \beta$ the subspaces $\mathfrak{g}_{\alpha}$ and $\mathfrak{g}_{\beta}$ are $B$-orthogonal.
			
			\item If $\alpha \in R$, then $- \alpha \in R$. Moreover, if $e_{\alpha} \in \mathfrak{g}_{\alpha}$ then $\bar{e}_{\alpha} \in \mathfrak{g}_{-\alpha}$ and 
				$$\left(\mathfrak{g}_{\alpha}\oplus \mathfrak{g}_{-\alpha}\right) \cap \mathfrak{g} = \mathbb{R}(e_{\alpha} + \bar{e}_{\alpha})\oplus \mathbb{R}(\textbf{i}e_{\alpha} - \textbf{i}\bar{e}_{\alpha})$$
	\end{enumerate}			
\end{Lemma}			
			
\begin{proof}
The first and the last statement follow directly from the definitions. For the second statement consider first the case $\beta = 0$ and $\alpha \in R$. Then follows for all $s,t \in \mathfrak{t}_0^c$
				$$B(\alpha(t)e_{\alpha}, s) = B([t, e_{\alpha}], s) = - B(e_{\alpha}, [t,s]) = 0$$
where we used in the second step that $B$ is $G^c$-invariant. This shows that $\mathfrak{t}_0^c$ is $B$-orthogonal to $\mathfrak{g}_{\alpha}$. Now consider $\alpha, \beta \in R$ with $\alpha + \beta \neq 0$. A similar calculation shows for all $s,t \in \mathfrak{t}_0^c$
				$$B(\alpha(t)e_{\alpha}, \beta(s)e_{\beta}) = B([t, e_{\alpha}], \beta(s)e_{\beta}) = - B(t, [e_{\alpha}, \beta(s)e_{\beta}]) = 0$$
where the last equality follows from the observation $[e_{\alpha}, \beta(s)e_{\beta}] \in \mathfrak{g}_{\alpha + \beta}$. 
\end{proof}

\paragraph{The Weyl group.} 
Using the inner product on $\mathfrak{g}$, we identify the roots $\alpha = \textbf{i}a \in R$ with vectors $t_{\alpha}\in \mathfrak{t}_0$ by the relation
				$$a(t) := \langle t_{\alpha}, t \rangle \qquad \text{for all $t \in \mathfrak{t}_0$}.$$
This yields a subset $\Phi_R = \{t_{\alpha}\,|\, \alpha \in R\}  \subset \mathfrak{t}_0$ which satisfies the properties of an abstract root system:
			\begin{enumerate}
						\item $\Phi_R$ is a spanning set for $\mathfrak{t}_0$.
						\item For every $t_{\alpha} \in \Phi_R$, the orthogonal reflection along $\ker{\alpha}$
											$$s_{\alpha}: \mathfrak{t}_0 \rightarrow \mathfrak{t}_0, \qquad s_{\alpha}(t) := t - \frac{2 \langle t, t_\alpha \rangle}{||t_{\alpha}||^2}$$ 
								carries $\Phi_R$ to itself.
						\item $\frac{2 \langle t_{\beta} , t_{\alpha} \rangle}{||t_{\alpha}||^2}$	is an integer for all $t_{\alpha}, t_{\beta} \in \Phi_R$.
			\end{enumerate}
This is discussed in \cite{LiegpKnapp} Chapters II.5. The subgroup $W$ generated by all the root reflection $s_{\alpha}$ inside the orthogonal group $O(\mathfrak{t}_0)$ is called the \textbf{Weyl group}. Since $\Phi_R$ is a spanning set of $\mathfrak{t}_0$, any orthogonal transformation which fixes $\Phi_R$ must be the identity and hence the Weyl group is always finite. After removing all hyperplanes $\text{ker}(\alpha)$ the Weyl group acts transitively and freely on $\mathfrak{t}_0 \backslash \cup\{\text{ker}(\alpha)\,|\, \alpha \in R\}$. The closure of a connected component of this space is called a \textbf{Weyl chamber} $\Omega_W \subset \mathfrak{t}_0$. In particular, $\Omega_W$ is the closure of a fundamental domain for the action of $W$. The Weyl group can alternatively be described as
							$$W \cong N_G(T)/Z_G(T).$$
Here the normalizer $N_G(T)$ acts on the maximal torus $T$ by conjugation. This action is trivial on the center $Z_0(G) \subset T$ and its derivative induces an action on $\mathfrak{t}_0$. Since the inner product on $\mathfrak{g}$ is $G$-invariant, this identifies $N_G(T)/Z_G(T)$ with a subgroup of the orthogonal group $O(\mathfrak{t}_0)$ and it is easy to check that this group permutes the roots $t_{\alpha}$. The equivalence of both descriptions of the Weyl-group is shown in \cite{LiegpKnapp} Chapters IV.6. Since any two maximal tori in $G$ are conjugated, this shows that the conjugation classes in $G$ are parametrized by $T/W$ and in particular any element $\xi \in \mathfrak{g}$ is conjugated to an element in the Weyl chamber $\Omega_W \subset \mathfrak{t}_0$.

\paragraph{Simple roots.}
Consider a notion of positivity on the set $R$ satisfying the properties
		\begin{enumerate} \label{rootseq0}
				\item For every root $\alpha \in R$ exactly one of $\alpha$ and $-\alpha$ is positive.
				\item If $\alpha$ and $\beta$ are positive, then $\alpha + \beta$ is positive.
		\end{enumerate}
An easy way to define such a notion goes as follows. Choose a real linear functional $\phi: \mathfrak{t}_0 \rightarrow \mathbb{R}$ such that $\ker{\phi}\cap \Phi_R = \emptyset$ and define a root $\alpha \in R$ to be positive whenever $\phi(t_{\alpha}) > 0$. We write $\alpha > 0$ for a positive root $\alpha$ and denote by $R^+$ the collection of positive roots. This induces a partial ordering on the roots according to the rule
					$$\alpha > \beta \quad \text{if and only if} \quad \alpha - \beta > 0.$$
A root $\alpha \in R^+$ is called simple if it cannot be decomposed as $\alpha = \beta + \gamma$ with $\beta, \gamma \in R^+$. In other words, a simple root is a minimal positive root. We denote by $R^+_0 = \{\alpha_1, \ldots, \alpha_r\}$ the set of simple roots. It is easy to deduce from the definitions that any root $\alpha$ can be written as
				\begin{align}\label{rootseq1} \alpha = \sum_{j=1}^r x_j \alpha_j \end{align}
with coefficients $x_1,\ldots,x_r \in \mathbb{Z}$ having all the same sign (or vanish). In particular $\Phi_{R_0^+}$ is a spanning set of $\mathfrak{t}_0$. A less obvious fact is that $\Phi_{R_0^+}$ is linear independent (see \cite{LiegpKnapp} II.5 Prop 2.49). Hence every root has a unique expression (\ref{rootseq1}) and a root is positive if and only if all the coefficients are nonnegative. This observation shows that the collection of simple roots and the partial ordering determine one another.

Any collection of simple roots $R^+_0 = \{\alpha_1, \ldots, \alpha_r\}$ determines a canonical Weyl chamber by the formula
			$$\Omega_W = \{ t \in \mathfrak{t}_0\,|\, \text{$a_j(t) \geq 0$ for all $j=1,\ldots,r$}\}$$
where we denote $\alpha_j = \textbf{i}a_j$ as above. Conversely, given a Weyl chamber $\Omega_W$ we can recover the collection of positive roots by the rule
			$$\text{$\alpha > 0$ if and only if $\langle t, t_\alpha \rangle \geq 0$ for all $t \in \Omega_W$}.$$
Hence the choice of a Weyl chamber and a partial ordering determine one another as well. Since any two Weyl chambers are conjugated by an element in $G$, this shows that all the choices in this section are canonical up to conjugation.

We denote the simple roots in $\Phi_{R_0^+}$ for convenience by $t_j := t_{\alpha_j}$. Since they define a basis of $\mathfrak{t}_0$, we can define a dual basis $\{\check{t}_1, \ldots, \check{t}_r\}$ by
							\begin{align} \label{rootseq2} \left\langle \check{t}_i , \frac{2 t_j}{||t_j||^2} \right\rangle = \delta_{ij} \end{align}
for $i,j = 1, \dots, r$. They are clearly contained in the Weyl chamber determined by the simple roots and yield the characterization
					$$t = \sum_{j=1}^r x_j \check{t}_j \in W_{\Omega} \qquad \Leftrightarrow \qquad \text{$x_j \geq 0$ for $j=1,\ldots,r$}.$$
The dual elements
				$$\lambda_j : \mathfrak{t}_0 \rightarrow \textbf{i}\mathbb{R}, \qquad \lambda_j(t) := \textbf{i} \langle \check{t}_j, t \rangle$$
are called the fundamental weights associated to the simple roots.

\subsubsection{An intrinsic definition of parabolic subgroups}

We provide an intrinsic definition of parabolic subgroups following the presentation \cite{Serre:1954} by Serre. Let $\xi \in [\mathfrak{g}, \mathfrak{g}]$ be given and choose a maximal torus $T \subset G$ such that $\xi \in \mathfrak{t}_0$. Moreover, let $R_0^+ := \{\alpha_1,\ldots, \alpha_r\}$ be a choice of simple roots such that $\xi$ is contained in the corresponding Weyl chamber. Denote
		\begin{align} \label{Req1}
		R(\xi) := \{ \alpha \in R\,|\, \langle \xi, t_{\alpha} \rangle \geq 0 \} \qquad \text{and} \qquad \tilde{R}(\xi) := \{ \alpha \in R\,|\, \langle \xi, t_{\alpha} \rangle = 0 \}.
		\end{align}
Define the Lie subalgebras
			\begin{align} \label{Peq1} \mathfrak{q}(\xi) := \mathfrak{z} \oplus \mathfrak{t}_0 \oplus \bigoplus_{\alpha \in R(\xi)} \mathfrak{g}_{\alpha} \end{align}
and
			\begin{align} \label{Peq2} \mathfrak{l}(\xi) := \mathfrak{z} \oplus \mathfrak{t}_0 \oplus \bigoplus_{\alpha \in \tilde{R}(\xi)} \mathfrak{g}_{\alpha}. \end{align}
The next Lemma shows that this notation is consistent with our definition in the section on toral generators.

\begin{Lemma}
Consider the setting from above and define $\mathfrak{q}(\xi)$ and $\mathfrak{l}(\xi)$ by (\ref{Peq1}) and (\ref{Peq2}) respectively. Then
			$$\mathfrak{q}(\xi) = \{ \rho \in \mathfrak{g}^c\,|\, \text{the limit $\lim_{t\rightarrow \infty} e^{\textbf{i}t\xi} \rho e^{-\textbf{i}t\xi}$ exists in $\mathfrak{g}^c$} \}$$
	and 
			$$\mathfrak{l}(\xi) = \{ \rho \in \mathfrak{g}^c\,|\,  e^{\textbf{i}t\xi} \rho e^{-\textbf{i}t\xi} = \rho \}.$$
\end{Lemma}

\begin{proof}
Decompose $\rho \in \mathfrak{g}^c$ with respect to the root space decomposition as
				$$\rho = \rho_0 + \sum_{\alpha \in R} \rho_{\alpha}$$
with $\rho_0 \in \mathfrak{t}$ and $\rho_{\alpha} \in \mathfrak{g}_{\alpha}$. By definition of the roots we have
				$$[ \textbf{i}\xi, \rho_a] = - a(\xi) \rho(\xi)\cdot \rho_{\alpha} =  - \langle t_{\alpha}, \xi \rangle \rho_{\alpha}$$
and hence
			$$e^{\textbf{i}t\xi} \rho e^{-\textbf{i}t\xi} = \rho_0 + \sum_{\alpha \in R} e^{- \langle t_{\alpha}, \xi \rangle t} \rho_{\alpha}.$$
This converges for $t \rightarrow \infty$ if and only if $\rho_{\alpha} = 0$ for all ${\alpha} \notin R(\xi)$. Similarly, we have $\rho = e^{\textbf{i}\xi} \rho e^{-\textbf{i}\xi}$ if and only if $\rho_{\alpha} = 0$ for all ${\alpha} \notin \tilde{R}(\xi)$.
\end{proof}

We could now define the parabolic subgroup $Q(\xi)$ and its Levi subgroup $L(\xi)$ as those connected subgroups of $G^c$ whose Lie algebras are given by $\mathfrak{q}(\xi)$ and $\mathfrak{l}(\xi)$ respectively. These are closed subgroups, since both agree with their normalizer in $G^c$.

\begin{Lemma} \label{PreLemma3}
Let $\check{t}_1, \ldots, \check{t}_r$ be defined by (\ref{rootseq2}) and let
					$$\xi = x_1 \check{t}_1 + \cdots + x_r \check{t}_r \in \Omega_W$$
with $x_j \geq 0$. Then $Q_j := Q(\check{t}_j)$ are maximal proper parabolic subgroups of $G^c$ and $Q(\xi) \subset Q(\check{t}_j)$ if and only if $x_j > 0$. Moreover, 
					$$Q(\xi) = \bigcap_{\{j\,|\, x_j > 0\}} Q(\check{t}_j).$$
\end{Lemma}

\begin{proof}
The proof is a simple matter of comparing $R(\check{t}_j)$ and $R(\xi)$.
\end{proof}


\section{Algebraic and symplectic stability}
Let $G$ be a compact connected Lie group and let $P \rightarrow \Sigma$ be a principal $G$ bundle over $\Sigma$. Denote by $G^c$ the complexification of $G$ and by $P^c := P \times_G G^c$ the complexified principal bundle.

The algebraic geometric construction of the moduli space of holomorphic structures on $P^c$, in the sense of Mumford's geometric invariant theory \cite{GIT:book}, depends on the notion of stable and semistable objects. For vector bundles this notion is due to Mumford \cite{Mumford:1963} and it was later extended by Ramanathan \cite{Ramanathan:1975} to principal bundles. We discuss these two definitions in the first subsection. We denote the corresponding GIT moduli space of holomorphic structures on $P^c$ by
								$$\mathcal{J}^{ss}(P^c)/\!/\mathcal{G}(P^c). $$
As mentioned in the introduction, this space is obtained by identifying two orbits in $\mathcal{J}^{ss}(P)/\mathcal{G}(P^c)$ when they cannot be separated.

The gauge group $\mathcal{G}(P)$ acts on $\mathcal{A}(P)$ in a hamiltonian way with moment map $\mu(A) = *F_A$ by Lemma \ref{PreLemma1}. For every central element $\tau \in Z(\mathfrak{g})$ one obtains the symplectic quotient
						$$\mathcal{A}(P)/\!/ \mathcal{G}(P) := \mu^{-1}(\tau)/\mathcal{G}(P).$$
Note that the moment map is not uniquely determined by the gauge action and  another moment map is given by $\mu_{\tau}(A) := *F_A - \tau$. This satisfies $\mu^{-1}_{\tau}(0) = \mu^{-1}(\tau)$ and we may therefore assume that the symplectic quotient is always obtained from the vanishing locus of some moment map. The symplectic version of GIT (see \cite{RobSaGeo}) yields a notion of stable and semistable objects in $\mathcal{A}(P)$ in terms of the moment map. We show in the second subsection that there exists a natural choice for $\tau \in Z(\mathfrak{g})$ which depends only on the topology of $P$. We call this the central type of $P$ and consider always the corresponding moment map $\mu_{\tau}$. It will follow from Theorem \ref{MLT} and Theorem \ref{ThmNUgen} in the next section that this yields
				$$\mu_{\tau}^{-1}(0)/\mathcal{G}(P) = \mathcal{A}^{ss}(P)/\!/\mathcal{G}^c(P).$$ 
The right hand side is again obtained by identifying orbits in $\mathcal{A}^{ss}(P)/\mathcal{G}^c(P)$ if they cannot be separated.

Using Lemma \ref{PreLem2} one can identify the space $\mathcal{J}(P^c)$ of holomorphic structure on $P^c$ with $\mathcal{A}(P)$. Theorem \ref{ThmNSRgen} shows that the notion of symplectic and algebraic stability are essentially equivalent. In particular, this yields an isomorphism
						$$\mathcal{J}^{ss}(P^c)/\!/\mathcal{G}(P^c) \cong \mathcal{A}^{ss}(P)/\!/\mathcal{G}^c(P) \cong \mu_{\tau}^{-1}(0)/\mathcal{G}(P)$$
if $\tau \in Z(\mathfrak{g})$ denotes the central type of $P$. The proof of this theorem will be based on the whole remainder of the exposition, namely on Proposition \ref{PropStabChar}, the moment-weight inequality (Theorem \ref{thmMWI}), the Harder-Narasimhan-Ramanathan theorem (Theorem \ref{ThmNSR}) and the dominant weight theorem (Theorem \ref{ThmDomWeight}).

\subsection{Algebraic stability}
Let $G$ be a compact connected Lie group and $P \rightarrow \Sigma$ be a principal $G$ bundle. We denote by $G^c$ the complexification of $G$ and by $P^c := P \times_G G^c$ the complexified principal bundle. We discuss the algebraic notion of stability on the space $\mathcal{J}(P^c)$  of holomorphic structures on the principal $G^c$ bundle $P^c$. This definition depends only on the complexified bundle $P^c$ itself and not on the reduction $P \subset P^c$. Consider as a warmup the case $G^c = \text{GL}(n)$. This allows us to identify $P^c$ with a complex vector bundle. The slope or normalized Chern class of a vector bundle $E \rightarrow \Sigma$ is defined as
					$$\mu(E) := \frac{c_1(E)}{\text{rk}(E)}.$$
The following Definition is due to Mumford \cite{Mumford:1963}.

\begin{Definition} \label{DefnStab1}
Let $E \rightarrow \Sigma$ be a holomorphic vector bundle. 
	\begin{enumerate}
			\item $E$ is called \textbf{stable} if for every proper holomorphic subbundle $0 \neq F \subset E$ we have $\mu(F) < \mu(E)$.
			\item $E$ is called \textbf{polystable} if $E$ is the direct sum of stable vector bundles all having the same slope.
			\item $E$ is called \textbf{semistable} if for every proper holomorphic subbundle $0 \neq F \subset E$ we have $\mu(F) \leq \mu(E)$.
			\item $E$ is called \textbf{unstable} if $E$ is not semistable.
	\end{enumerate}		
\end{Definition}

The analogue of this definition for general Lie groups was formulated by Ramanathan \cite{Ramanathan:1975}. Lemma \ref{RStabLem} below shows that Definition \ref{DefnStab1} corresponds to the special case $G^c= \text{GL}(n)$ in Definition \ref{DefnStab2}.

\begin{Definition}\label{DefnStab2}
Let $G^c$ be a connected reductive Lie group and $P^c \rightarrow \Sigma$ be a holomorphic principal $G^c$ bundle. 
	\begin{enumerate}
			\item $P^c$ is called \textbf{stable} if for every holomorphic reduction $P_{Q} \subset P^c$ to a maximal proper parabolic subgroup $Q \subset G^c$ the subbundle $\text{ad}(P_Q) \subset \text{ad}(P^c)$ satisfies $c_1(\text{ad}(P_Q)) < 0$.
			\item $P^c$ is called \textbf{polystable} if there exists a parabolic subgroup $Q \subset G^c$ and a holomorphic reduction $P_L \subset P^c$ to a Levi subgroup of $Q$ satisfying the following
						\begin{enumerate}
								\item $P_L$ is a stable principal $L$ bundle.
								\item For every character $\chi: L \rightarrow \mathbb{C}^*$, which is trivial on the center of $G^c$, the associated line bundle $\chi(P_L) := P_L\times_{\chi}\mathbb{C}$ satisfies $c_1(\chi(P_L)) = 0$.
						\end{enumerate}
			\item $P^c$ is called \textbf{semistable} if for every holomorphic reduction $P_{Q} \subset P^c$ to a maximal proper parabolic subgroup $Q \subset G^c$ the subbundle $\text{ad}(P_Q) \subset \text{ad}(P^c)$ satisfies $c_1(\text{ad}(Q)) \leq 0$.
			\item $P^c$ is called \textbf{unstable} if $\text{ad}(P^c)$ is not semistable.
	\end{enumerate}		
\end{Definition}

\begin{Remark}
Let $L_1,\cdots, L_r$ and $G^c$ be complex connected reductive Lie groups such that the product $L_1 \times \cdots \times L_r \subset G^c$ embeds as a subgroup. Let $P_j$ be stable principal $L_j$ bundles for $j=1,\ldots,r$. Then it is easy to see that $P_L := P_{L_1}\times \cdots \times P_{L_r}$ is a stable principal $L$ bundle. However, the extension $P^c := P_L \times_{L} G^c$ is in general not a semistable $G^c$-bundle. The second condition in the definition of polystability in needed to guarantee the semistability of $P^c$. To see this let $P_{Q'} \subset P^c$ be the reduction to a maximal parabolic subgroup and consider the determinant of the adjoint action of $Q' \subset G^c$ on its Lie algebra. This character is clearly trivial on the center of $G^c$ and either restricts to $L$ or to a maximal parabolic subgroup $Q'' = Q'\cap L \subset L$. In the first case, it follows from the definition of polystability that $c_1(\text{ad}(P_{Q'})) = 0$. In the other case observe that $P_{Q'}$ determines a maximal parabolic reduction $P_{Q''} \subset P_L$ and $c_1(\text{ad}(P_{Q'})) = c_1(\text{ad}(P_{Q''})) < 0$, since $P_L$ is stable.
\end{Remark}

\begin{Lemma} \label{RStabLem}
A holomorphic vector bundle $E$ is stable, polystable, semistable or unstable if and only its $\text{GL}(n)$-frame bundle $P^c := \text{Fr}(E)$ is stable, polystable, semistable or unstable respectively.
\end{Lemma}

\begin{proof}
We discuss the stable (resp. semistable) case first. A maximal parabolic subgroup of $GL(n)$ is the stabilizer a subspace $0 \neq V \subset \mathbb{C}^n$ and the holomorphic reduction $P_Q$ of the $\text{GL}(n)$-frame bundle to a maximal parabolic subgroup is thus the stabilizer of a holomorphic subbundle $F \subset E$. Consider the orthogonal splitting $E = F \oplus G$ with respect to some fixed hermitian metric on $E$. Then $\text{ad}(P_Q) \subset \text{End}(E)$ is given by the space of upper block diagonal matrices. We choose unitary connections $A_1$ on $E$ and $A_2$ of $G$ and denote by $A$ the induced connection of $E = F\oplus G$. This induces also a connection on $\text{ad}(P_Q)$ and the curvature of this connection is given by the endomorphism
			$$\xi \mapsto F_A \xi - \xi F_A$$
for $\xi \in \text{ad}(P_Q)$. Since $F_A = \text{diag}(F_{A_1}, F_{A_2})$ is block-diagonal, a short calculation shows that the trace of this map is given by $\text{rk}(G)\text{tr}(F_{A_1}) - \text{rk}(F)\text{tr}(F_{A_2})$ and Chern-Weyl theory yields
				\begin{align*}
						c_1(\text{ad}(P_Q)) &= \text{rk}(G)c_1(F) - \text{rk}(F)c_1(G)  \\
																&= \text{rk}(E/F) \text{rk}(F) \left(  \frac{c_1(F)}{\text{rk}(F)} - \frac{c_1(E/F)}{\text{rk}(E/F)}\right).
				\end{align*}
This expression is nonpositive if and only if $c_1(F)/\text{rk}(F) \leq c_1(E)/\text{rk}(E)$ and negative whenever strict inequality holds. This proves the equivalence of both definitions in the stable and semistable case. 

The unstable case is equivalent to semistable case and it remains to discuss the polystable case. A general parabolic subgroup of $\text{GL}(V)$ is the stabilizer of a filtration $V_1 \subset \cdots \subset V_r = V$ and a Levi subgroup in given as the stabilizer of a splitting $V = W_1 \oplus \cdots \oplus W_r$ with $V_j = W_1 \oplus\cdots \oplus W_j$. Hence, a holomorphic reduction $P_L \subset P^c$ to the Levi factor of a parabolic subgroup corresponds to the $L = \text{GL}(n_1)\times \cdots \times \text{GL}(n_r)$ frame bundle of a holomorphic splitting
				$$ E = E_1 \oplus \cdots \oplus E_r.$$
We claim that $P_L$ is a stable principal $L$ bundle if and only if all factors $E_j$ are stable holomorphic vector bundles. Indeed, a maximal parabolic subgroup of $L$ has the shape 
				$$Q = \text{GL}(n_1) \oplus \cdots \oplus \text{GL}(n_{j-1}) \oplus Q_j \oplus \text{GL}(n_{j+1}) \cdots \oplus \text{GL}(n_r)$$
where $Q_j \subset \text{GL}(n_j)$ is a maximal parabolic subgroup. Then 
			$$\text{ad}(P_Q) = \text{End}(E_1) \oplus \cdots \oplus \text{ad}(P_{Q_j}) \oplus \cdots \oplus \text{End}(E_r)$$
and hence $c_1(\text{ad}(P_Q)) = c_1(\text{ad}(P_{Q_j}))$. The claim follows now from our discussion of the stable case. 

It remains to verify that the slopes of all subbundles satisfy $\mu(E_j) = \mu(E)$ if and only if for every character $\chi: L \rightarrow \mathbb{C}^*$ which is trivial on the center of $\text{GL}(n)$ the associated line bundle $\chi(P_L)$ has degree zero. Every character $\chi : L \rightarrow \mathbb{C}^*$ factors as $\chi = \chi_1 \cdots \chi_r$ with $\chi_j: \text{GL}(n_i) \rightarrow \mathbb{C}^*$ and induces on the Lie algebra the representation
					$$\dot{\chi} = \dot{\chi}_1 + \cdots + \dot{\chi}_r$$
with $\dot{\chi_j} := d \chi_j(\mathds{1}): \mathfrak{gl}(n_j) \rightarrow \mathbb{C}$. Since every traceless matrix in $\mathfrak{gl}(n_j)$ is a commutator, there exist $\lambda_j \in \mathbb{C}$ such that 
			$$\dot{\chi}_j(\rho_j) = \lambda_j\text{tr}(\rho_j)$$
for all $\rho_j \in \mathfrak{gl}(n_j)$. We choose unitary connections $A_j \in \mathcal{A}(E_j)$ and denote by $A = A_1 \oplus \cdots \oplus A_r$ the induced unitary connection on $E$. Then follows from Chern-Weil theory
			\begin{align*}
					c_1(\chi(P_L)) &= \frac{\textbf{i}}{2\pi} \int_{\Sigma} F_{\chi(A)} \, dvol_{\Sigma} = \frac{\textbf{i}}{2\pi} \int_{\Sigma} \dot{\chi}(F_A)\, dvol_{\Sigma} \\
												 &= \lambda_1 c_1(E_1) + \cdots + \lambda_r c_1(E_r).
			\end{align*}
Note that $\chi$ vanishes on the center of $\text{GL}(n)$ if and only if $n_1\lambda_1 + \cdots + n_r\lambda_r = 0$ is satisfied. If in addition $\mu(E_j) = \mu(E)$ holds for all $j$, then
			$$c_1(\chi(P_L)) = \sum_{j=1}^r n_j \lambda_j \mu(E_j) = 0.$$ 
For the converse consider the character $\chi: \text{GL}(n_1)\times \cdots \times \text{GL}(n_r) \rightarrow \mathbb{C}^*$
				$$\chi(B_1,\ldots,B_r) := \det(B_j)^n \det(B)^{-n_j}.$$
This vanishes on the center of $\text{GL}(n)$ and satisfies $\dot{\chi}(\xi) = n \text{tr}(\xi_j) - n_i \text{tr}(\xi)$. Hence
			$$c_1(\chi(P_L)) = n c_1(E_j) - n_j c_1(E)$$
and this vanishes precisely when $\mu(E_j) = \mu(E)$ is satisfied.
\end{proof}

The next Lemma shows that we can always reduce to the case where $G^c$ has discrete center.

\begin{Lemma} \label{LemmaSSred1}
Let $G^c$ be a complex connected reductive Lie group and $P^c \rightarrow \Sigma$ be a principal $G^c$ bundle. Denote by $Z_0(G^c)$ the connected component of the center of $G^c$ containing the identity. Let $H^c := G^c/Z_0(G^c)$ and denote by 
							$$P_{H^c} := P^c/Z_0(G^c)$$ 
the associated $H^c$ bundle. This carries a natural induced holomorphic structure and $P^c$ is stable, polystable, semistable or unstable if and only if $P_{H^c}$ is stable, polystable, semistable or unstable respectively.
\end{Lemma}

\begin{proof}
The Lie algebra of $G^c$ splits as  $\mathfrak{g}^{c}=Z(\mathfrak{g}^c)\oplus[\mathfrak{g}^{c},\mathfrak{g}^{c}]$ and $[\mathfrak{g}^{c},\mathfrak{g}^{c}]$ can be identified with the semisimple Lie algebra of $H^c$. This splitting is preserved by the adjoint action of $G^c$ and produces a splitting $\text{ad}(P^c) = V \oplus \text{ad}(P_{H^c})$ where $V = \Sigma \times Z(\mathfrak{g}^c)$ is a trivial bundle. Parabolic subgroups $Q \subset G^c$ correspond bijectively to parabolic subgroups $\bar{Q} := Q/Z_0(G^c) \subset H$ and parabolic reductions $P_Q \subset P^c$ correspond bijectively to parabolic reductions $P_{\bar{Q}} := P_Q/Z_0(G^c) \subset P_{H^c}$. Since $\text{ad}(P_Q) = V \oplus \text{ad}(P_{\bar{Q}})$, we have $c_1(\text{ad}(P_Q)) = c_1(\text{ad}(P_{\bar{Q}}))$ and this shows that $P^c$ is stable (resp. semistable) if and only if $P_{H^c}$ is stable (resp. semistable).

If $L$ is a Levi subgroup of the parabolic subgroup $Q \subset G^c$, then $\bar{L} := L/Z_0$ is a Levi-subgroup of $\bar{Q} = Q/Z_0 \subset H^c$. Moreover, reductions $P_L \subset P^c$ to $L$ correspond bijectively to reductions $P_{\bar{L}} = P_L/Z_0(G^c) \subset P_{H^c}$. We have already shown that $P_{L}$ is stable if and only if $P_{\bar{L}}$ is stable. The characters $\chi: Q \rightarrow \mathbb{C}^*$ which are trivial on the center $Z(G^c)$ of $G^c$ correspond bijectively to the characters $\bar{\chi}: \bar{Q} \rightarrow \mathbb{C}^*$ which are trivial on $Z(G^c)/Z_0(G^c)$ and
				$$\chi(P_Q) \cong \bar{\chi}(P_{\bar{Q}}).$$
Thus $P^c$ is polystable if and only if $P_{H^c}$ is polystable.
\end{proof}

\subsection{Symplectic stability}
Let $G$ be a compact connected Lie group and $P \rightarrow \Sigma$ a principal $G$ bundle. Let $\chi: G \rightarrow S^1$ be a character and denote by $\dot{\chi} = d\chi(\mathds{1}) : \mathfrak{g} \rightarrow \textbf{i}\mathbb{R}$ the induced character on the Lie algebra. Since $\mathfrak{g} = Z(\mathfrak{g})\oplus[\mathfrak{g}, \mathfrak{g}]$, we may identify $-2\pi\textbf{i}\dot{\chi}$ with an element in $Z(\mathfrak{g})^* = \text{Hom}(Z(\mathfrak{g}), \mathbb{R})$. Denote by $\chi(P) := P \times_{\chi} \mathbb{C}$ the line bundle associated to $P$ via $\chi$. Then
		\begin{align} \label{eqstab1} c_1(\chi(P)) = \frac{\textbf{i}}{2\pi} \int_{\Sigma} \dot{\chi}(F_A) \end{align}
for any connection $A \in \mathcal{A}(P)$. The assignment $\dot{\chi} \mapsto -2\pi\textbf{i}c_1(\chi(P))$ extends to a unique element in $Z(\mathfrak{g})^{**}$, since the lattice of all infinitesimal characters spans $Z(\mathfrak{g})^*$ as a vector space. This corresponds under the canonical isomorphism $Z(\mathfrak{g})^{**} \cong Z(\mathfrak{g})$ to an element $\tau \in Z(\mathfrak{g})$ which satisfies
		\begin{align} \label{eqCT} \alpha(\tau) = \int_{\Sigma} \alpha(F_A) \qquad \text{for all $\alpha \in Z(\mathfrak{g})^*$ and $A \in \mathcal{A}(P)$.}\end{align}
Here we identify $Z(\mathfrak{g})^* \subset \mathfrak{g}^*$ with the subspace of linear functionals vanishing on $[\mathfrak{g},\mathfrak{g}]$. We call $\tau$ the \textbf{central type} of $P$.

\begin{Remark} \label{RmkCT}
Recall our standing assumption $\text{vol}(\Sigma) = 1$ and suppose that $A \in \mathcal{A}(P)$ satisfies $*F_A = \lambda$ for some $\lambda \in Z(\mathfrak{g})$. Then (\ref{eqCT}) yields
				$$\alpha(\tau) = \int_{\Sigma} \alpha(F_A) =  \int_{\Sigma} \alpha(\lambda) \, dvol_{\Sigma} = \alpha(\lambda)$$
for all $\alpha \in Z(\mathfrak{g})^*$ and hence $\lambda = \tau$.
\end{Remark}

Let $\tau \in Z(\mathfrak{g})$ be defined by (\ref{eqCT}). It follows from Lemma \ref{PreLem1} that
				\begin{align} \label{eqMM} \mu_{\tau}: \mathcal{A}(P) \rightarrow L^2(\Sigma, \text{ad}(P)), \qquad \mu(A) := *F_A - \tau \end{align}
is a moment map for the action of the gauge group $\mathcal{G}(A)$ on $\mathcal{A}(P)$. The following definition is the precise analogue of Definition 7.1 in \cite{RobSaGeo} with respect to this moment map.

\begin{Definition} \label{DefnStab3}
Let $G$ be a compact connected Lie group, let $P \rightarrow \Sigma$ be a principal $G$ bundle with central type $\tau \in Z(\mathfrak{g})$ defined by (\ref{eqCT}), and define $\mu_{\tau}$ by (\ref{eqMM}). For $A \in \mathcal{A}(P)$ denote by $\overline{\mathcal{G}^c(A)}$ the $W^{1,2}$-closure of the complex gauge orbit $\mathcal{G}^c(A)$. 
	\begin{enumerate}
		\item $A$ is called \textbf{$\mu_{\tau}$-stable} if and only if $\mathcal{G}^c(A) \cap \mu_{\tau}^{-1}(0)\cap \mathcal{A}^*(P) \neq \emptyset$ where $\mathcal{A}^*(P)$ denotes the irreducible connections on $P$.
		\item $A$ is called \textbf{$\mu_{\tau}$-polystable} if and only if $\mathcal{G}^c(A) \cap \mu_{\tau}^{-1}(0) \neq \emptyset$.
		\item $A$ is called \textbf{$\mu_{\tau}$-semistable} if and only if $\overline{\mathcal{G}^c(A)} \cap \mu_{\tau}^{-1}(0) \neq \emptyset$.
		\item $A$ is called \textbf{$\mu_{\tau}$-unstable} if and only if $\overline{\mathcal{G}^c(A)} \cap \mu_{\tau}^{-1}(0) = \emptyset$.
	\end{enumerate}
\end{Definition}

\begin{Remark}
We call $A \in \mathcal{A}(P)$ an irreducible connection if the map
		$$d_{A}: W^{2,2}(\Sigma, \text{ad}(P)) \rightarrow W^{1,2}(\Sigma, T^*\Sigma\otimes\text{ad}(P))$$
is injective. In particular, $\mu_{\tau}$-stable connections can only exist when the center of $G$ is discrete and $\tau = 0$.
Since the infinitesimal action of the gauge group is given by $L_A \xi := -d_A \xi$, a connection $A$ is irreducible if and only if the isotropy group of the orbit $\mathcal{G}(A)$ is discrete. Suppose that $A$ is an irreducible connection satisfying $\mu_0(A) = *F_{A} = 0$. The infinitesimal action of the complexified gauge group
		$$L_{A}(\xi + \textbf{i}\eta) = - d_{A}\xi - * d_{A} \eta$$
is readily seen to be injective in this case: Assume $L_{A}(\xi + \textbf{i}\eta) = 0$ and apply $d_{A}$ to the equation. Then follows $d_{A}^*d_{A} \eta = 0$ and hence $d_{A}\eta = 0$. Since $A$ is irreducible, we conclude $\eta = 0$ and then $\xi = 0$. This argument shows that the $\mu_{\tau}$-stable orbits are precisely the $\mu_{\tau}$-polystable orbits with discrete $\mathcal{G}^c(P)$ isotropy.	
\end{Remark}

The next Lemma relates the different notions of stability on $P$ and on the quotient bundle $P_H := P/Z_0(G)$ with fiber $H := G/Z_0(G)$. Note that $P_H$ has central type $0$ since its center is discrete.

\begin{Lemma} \label{LemmaSSred2}
Let $G$ be a compact connected Lie group and $P \rightarrow \Sigma$ a principal $G$ bundle of central type $\tau \in \mathfrak{g}$. Denote $P_H := P/Z_0(G)$ be the associated $H := G/Z_0(G)$ bundle. Let $A \in \mathcal{A}(P)$ and denote by $A_H \in \mathcal{A}(P_H)$ the induced connection. Then
		\begin{enumerate}
			\item $A_H$ is $\mu_0$-stable if and only if $A$ is $\mu_{\tau}$-polystable and the kernel of the infinitesimal action
							$$L_{A}: W^{2,2}(\Sigma, \text{ad}(P^c)) \rightarrow W^{1,2}(\Sigma, T^*\Sigma\otimes\text{ad}(P))$$
											$$L_A(\xi + \textbf{i}\eta) = - d_A \xi - * d_A \eta$$
						consists of constant central sections.
			\item $A_H$ is $\mu_0$-polystable if and only if $A$ is $\mu_{\tau}$-polystable.
			\item $A_H$ is $\mu_0$-semistable if and only if $A$ is $\mu_{\tau}$-semistable.
			\item $A_H$ is $\mu_0$-unstable if and only if $A$ is $\mu_{\tau}$-unstable.
		\end{enumerate}
\end{Lemma}

\begin{proof}
We begin with the polystable case. Every constant central curvature connections on $P$ clearly induces a flat connection on $P_{H}$. Conversely, assume that $A_1$ is a flat connection on $P_{H}$. As a general property of compact Lie groups, there exists an exact sequence
			\begin{align} \label{redSSeq1} 1 \rightarrow F \rightarrow Z_0(G)\times [G,G] \rightarrow G \rightarrow 1 \end{align}
where $F = Z_0(G)\cap[G,G]$ is a finite group. From this follows the exact sequence
			\begin{align} \label{redSSeq2} 1 \rightarrow F \rightarrow  G \rightarrow (G/Z_0(G))\times (Z_0(G)/F) \rightarrow 1. \end{align}
Consider the associated $(G/Z_0(G))\times (Z_0(G)/F)$ bundle 
			\begin{align} \label{redSSeq3} \tilde{P} = P\times_{G} ((G^c/Z_0)\times(Z_0/F)) = P_H \times_{\Sigma} P_2 \end{align}
where $P_2$ is a principal $Z_0(G)/F$-bundle over $\Sigma$. Since $Z_0(G)/F$ is connected and abelian, it is a torus and $P_2$ is isomorphic to the direct sum of $S^1$ bundles. It follows from Hodge theory that every line bundle admits a connection with constant central curvature and these yield a connection $A_2$ on $P_2$ with constant central curvature. Together with $A_1$ we obtain an induces a connection on $\tilde{P}$ which lifts to a connection on $P$ with constant central curvature. It follows from Remark \ref{RmkCT} that the curvature of this connection is given by $\tau$.

For the proof of the stable case observe that $\text{ad}(P^c) \cong V \oplus \text{ad}(P_H^c)$ where $ V = \Sigma \times Z(\mathfrak{g}^c)$ denotes the trivial $Z(\mathfrak{g}^c)$ bundle. The infinitesimal action
		$$L_{A}: W^{2,2}(\Sigma, \text{ad}(P^c)) \rightarrow W^{1,2}(\Sigma, T^*\Sigma\otimes\text{ad}(P))$$ 
agrees with $L_{A_H}$ on $\text{ad}(P_H^c)$. Since $d_A$ restricts to a flat connection on $V$ it follows $\text{ker}(L_A) \cong Z(\mathfrak{g}^c) \oplus \text{ker}(L_{A_H})$ and this shows the claim.

It remains to discuss the semistable case. Assume first that $A$ is $\mu_{\tau}$-semistable. Then exist connections $A^k \in \mathcal{G}^c(A)$ such that $A^k \rightarrow A^+$ for $k \rightarrow \infty$ and $\mu_{\tau}(A^+) = 0$. The induced connections $A^k_H \in \mathcal{A}(P_H)$ are clearly contained in $\mathcal{G}^c(A_H)$ and converge to the induced connection $A_H^+$. Since $\mu_{\tau}(A^+) = 0$, it follows that $\mu_0(A_H^+) = 0$ and hence $A_H$ is $\mu_0$-semistable.

For the converse, we consider the exact sequences (\ref{redSSeq1}) and (\ref{redSSeq2}) from above. Then (\ref{redSSeq3}) yields a finite covering 
				$$P \rightarrow \tilde{P} = P_H \times_{\Sigma} P_2$$
with covering group $F = Z_0(G) \cap [G,G]$. We have seen above that $P_2$ is a polystable $Z_0(G)/F$-bundle. Note that the natural identification $\mathcal{A}(\tilde{P}) = \mathcal{A}(P_H)\times\mathcal{A}(P_2)$ yields an inclusion
			\begin{align} \label{redSSeq4} \mathcal{G}^c(P_H) \times \mathcal{G}^c(P_2) \subset \mathcal{G}^c(\tilde{P}). \end{align}
Moreover, since $\text{Ad}(P^c) \rightarrow \text{Ad}(\tilde{P}^c)$ is a finite covering with covering group $F \subset Z_0(G^c)$, it is easy to see that every gauge transformation in $\mathcal{G}^c(\tilde{P})$ lifts to an element in $\mathcal{G}^c(P)$ and this lift commutes with the natural identification $\mathcal{A}(P) = \mathcal{A}(\tilde{P})$. 

Now assume that $A \in \mathcal{A}(P)$ induces a $\mu_0$-semistable connection $A_H \in \mathcal{A}(P_H)$. Since $P_2$ is polystable, it follows from (\ref{redSSeq4}) that there exists $g_0 \in \mathcal{G}^c(P)$ such that $g_0(A)$ induces $A_H \in \mathcal{A}(P_H)$ and a connection $A_2 \in \mathcal{A}(P_2)$ with constant central curvature. Since $A_H$ is $\mu_0$-semistable, using (\ref{redSSeq3}) again, there exists gauge transformations $g_k \in \mathcal{G}^c(P)$ such that $g_k(g_0(A))$ induce the same connection $A_2$ on $P_2$ and induce a sequence of connections $A_H^k$ on $P_H$ which converges to a flat connection $A_H^+$. Clearly, $g_k(g_0 A)$ converges to the connection $A^+$ which is induced by $A_2$ and $A_H^+$. Hence $A^+$ has constant central curvature and it follows from Remark \ref{RmkCT} that $*F_{A^+} = \tau$. This completes the proof of the semistable case.

\end{proof}

\subsection{Equivalence of algebraic and symplectic stability}

The following theorem shows that the algebraic notion of stability from Definition \ref{DefnStab2} and the symplectic notion of $\mu_{\tau}$-stability from Definition \ref{DefnStab2} are essentially equivalent. 

\begin{Theorem}[\textbf{Generalized Narasimhan-Seshadri-Ramanathan Theorem}] \label{ThmNSRgen}
Let $G$ be a compact connected Lie group and $P \rightarrow \Sigma$ a principal $G$ bundle with central type $\tau \in Z(\mathfrak{g})$ defined by (\ref{eqCT}). Let $A \in \mathcal{A}(P)$ and consider the complexified bundle $P^c := P \times_G G^c$ with the induced holomorphic structure $J_A$. Then
		\begin{enumerate}
				\item $(P^c, J_A)$ is stable if and only if $A$ is $\mu_{\tau}$-polystable and the kernel of 
							$$L_{A}: W^{2,2}(\Sigma, \text{ad}(P^c)) \rightarrow W^{1,2}(\Sigma, T^*\Sigma\otimes\text{ad}(P))$$
											$$L_A(\xi + \textbf{i}\eta) = - d_A \xi - * d_A \eta$$
						contains only constant central sections.
				\item $(P^c,J_A)$ is polystable if and only if $A$ is $\mu_{\tau}$-polystable.
				\item $(P^c,J_A)$ is semistable if and only if $A$ is $\mu_{\tau}$-semistable.
				\item $(P^c,J_A)$ is unstable if and only if $A$ is $\mu_{\tau}$-unstable.
		
		\end{enumerate}
\end{Theorem}

The stable case was first proven by Narasimhan and Seshadri \cite{NarasimhanSeshadri:1965} for $G=U(n)$ and later extended by Ramanathan \cite{Ramanathan:1975} to arbitrary compact Lie groups. They establish these results using algebraic geometric methods.The first analytic proof was given by Donaldson \cite{Donaldson:1983NS} for the case $G = U(n)$. We present a different approach given by Bradlow \cite{Bradlow:1991} and Mundet \cite{Mundet:2000} in Theorem \ref{ThmNSR}. The equivalence of both definitions for semistability is essentially contained in the work of Atiyah and Bott \cite{AtBott:YangMillsEq}.

\begin{proof}[Proof of Theorem \ref{ThmNSRgen}.] We assume the following results for the proof:
			\begin{itemize}
						\item the characterization of algebraic stability in Proposition \ref{PropStabChar}
						\item the moment-weight inequality (Theorem \ref{thmMWI})
						\item the Harder-Narasimhan-Ramanathan theorem (Theorem \ref{ThmNSR})
						\item the dominant weight theorem (Theorem \ref{ThmDomWeight})
			\end{itemize}
We establish these results independently in the remainder of the exposition.

We may assume by Lemma \ref{LemmaSSred1} and Lemma \ref{LemmaSSred2} that $Z(G)$ is discrete and $\tau = 0$. The stable case is equivalent to Theorem \ref{ThmNSR}.

We deduce the polystable case from the stable case: Assume first that $P^c$ is polystable. Then there exists a reductive subgroup $L \subset G^c$ and a holomorphic reduction $P_L \subset P^c$ which is stable. We may assume that $L = K^c$ is the complexification of a compact subgroup $K \subset G$. Then we have an induced reduction $P_K \subset P$ such that $P_K^c$ agrees with the complexification of $P_K$ and it follows from the construction in Lemma \ref{PreLem2} that $A$ restricts to a connection on $P_K$. Assuming the stable case (i.e. Theorem \ref{ThmNSR}) we conclude that there exists a gauge transformation $g \in \mathcal{G}^c(P_K) \subset \mathcal{G}^c(P)$ such that $*F_{gA} = \tau_K \in Z(\mathfrak{k})$. It remains to show that $\tau_K \in Z(\mathfrak{g}) = 0$ vanishes. If $\tau_K \neq 0$ then exists a character $\chi: L \rightarrow \mathbb{C}^*$ with $\dot{\chi}(\tau_K) \neq 0$. Since $Z(G)$ is finite, we may replacing $\chi$ by a suitable power and assume that it is trivial on $Z(G^c)$. Using the definition of polystability then yields the contradiction
			$$0 = c_1(\chi(P_L)) = \frac{\textbf{i}}{2\pi} \int_{\Sigma} \dot{\chi}(F_A)\, dvol_{\Sigma} = \frac{\textbf{i}}{2\pi} \dot{\chi}(\tau_K) \neq 0.$$	

For the converse, assume that $A \in \mathcal{A}(P)$ is a flat connection. Let $H \subset G$ be the holonomy subgroup and $P_H \subset P$ be a reduction to the holonomy. Let $K := C_G(Z(H))$ be the centralizer of the center of the holonomy and denote the induced connection on $P_K = P_H \times_H K$ again by $A$. It is well-known that the isotropy subgroup of $A$ consists of constant gauge transformations and is naturally isomorphic to the centralizer of its holonomy, i.e.
				$$\mathcal{G}_A := \{g \in \mathcal{G}(P_K)\,|\, g(A) = A\} \cong C_K(H).$$ 
Comparing the Lie algebras of both sides, one checks that $C_K(H) = Z_0(K)$ is satisfied and $A \in \mathcal{A}(P_K)$ has only trivial isotropy. It follows now from the stable case (i.e. Theorem \ref{ThmNSR}) that $P_K^c$ is a stable principal $L = K^c$ bundle. Note that $L$ is a Levi-subgroup of a parabolic subgroup of $G^c$, since $K$ is the centralizer of a torus in $G$. Since $F_A = 0$, we have for any character $\chi : L \rightarrow \mathbb{C}^*$
				$$c_1(\chi(P_L)) = \frac{\textbf{i}}{2\pi} \int_{\Sigma} \dot{\chi}(F_A)\, dvol_{\Sigma} = 0$$
and hence $P^c$ is polystable.

Assume that $P^c$ is unstable. By Proposition \ref{PropStabChar} there exists $\xi \in \Omega^0(\Sigma, \text{ad}(P))$ with $w_0(A,\xi) < 0$. The moment-weight inequality (Theorem \ref{thmMWI}) yields $\mu_{0}(gA) \geq - w_0(A,\xi)/||\xi|| > 0$ for all $g \in \mathcal{G}^c(P)$ and hence $A$ is $\mu_0$-unstable. Assume conversely that $A$ is $\mu_0$-unstable. The dominant-weight theorem (Theorem \ref{ThmDomWeight}) shows that there exists $\xi \in \Omega^0(\Sigma,\text{ad}(P))$ such that $w_0(A,\xi) < 0$ and hence $P^c$ is unstable by Proposition \ref{PropStabChar}. This completes the proof of the unstable case and the semistable case is equivalent to this case.
\end{proof}


\section{The Yang-Mills flow and symplectic stability}

Let $G$ be a compact connected Lie group and let $P \rightarrow \Sigma$ be a principal $G$ bundle of central type $\tau \in Z(\mathfrak{g})$ defined by (\ref{eqCT}). In the differential geometric approach towards GIT the moment map squared functional plays a crucial role. This is defined by
			\begin{align} \label{YMeq01}
		\mathcal{F}_{\tau}: \mathcal{A}(P) \rightarrow \mathbb{R}, \qquad \mathcal{F}_\tau(A) := \frac{1}{2} \int_{\Sigma} ||*F_A - \tau||^2 \, dvol_{\Sigma}. 
			\end{align}
Note that (\ref{eqCT}) implies $\int_{\Sigma} \langle F_A, \tau \rangle = ||\tau||^2$ for every connection $A \in \mathcal{A}(P)$ and hence
			\begin{align} \label{YMeq02}
		\mathcal{F}_{\tau}(A) = \frac{1}{2}\int_{\Sigma} ||*F_A - \tau||^2 dvol_{\Sigma} = \frac{1}{2} \left(\int_{\Sigma} ||F_A||^2 dvol_{\Sigma} - ||\tau||^2 \right).
			\end{align}
Thus $\mathcal{F}_{\tau}$ agrees up to a constant shift with the Yang-Mills functional
		\begin{align}
						\label{YMeq1} \mathcal{YM}: \mathcal{A}(P) \rightarrow \mathbb{R}, \qquad \mathcal{YM}(A) := \frac{1}{2}\int_{\Sigma} ||F_A||^2\, dvol_{\Sigma}.
		\end{align}
Rade showed in his thesis \cite{Rade1992} that the negative gradient flow of the Yang-Mills functional is well-defined and converges if the base manifold has dimension $2$ or $3$. We summarize his results in the first subsection. Recall that we always consider the $W^{1,2}$-topology on $\mathcal{A}(P)$ when nothing else is specified. 

A crucial observation is the following: Any solution of the Yang-Mills flow remains in a single complexified orbit and there exists a canonical lift of a solution $A(t)$ of the Yang-Mills flow under the projection $\mathcal{G}^c \rightarrow \mathcal{G}^c(A)$ to a curve in $\mathcal{G}^c(P)$. Since the Yang-Mills flow is $\mathcal{G}(P)$-invariant, the geometric importance lies within the projection of such curves in $\mathcal{G}^c(P)/\mathcal{G}(P)$. The fibers of this quotient coincide with the homogeneous space $G^c/G$ which is a complete, connected, simply connected Riemannian manifold of nonpositive sectional curvature (see \cite{RobSaGeo} Appendix A and B). This underlying geometry is crucial for the following application.

As a first application, we establish the moment limit theorem (Theorem \ref{MLT}) and the analogue of the Ness uniqueness theorem in Theorem \ref{ThmNUgen} following the line of arguments in \cite{RobSaGeo}. The first result says that the limit $A_{\infty} := \lim_{t\rightarrow \infty} A(t)$ of the Yang-Mills flow starting at $A_0 \in \mathcal{A}(P)$ minimizes the Yang-Mills functional over the complexified orbit $\mathcal{G}^c(A_0)$. The second result asserts that any connection in the $W^{1,2}$-closure of $\mathcal{G}(A_0)$ which minimizes the Yang-Mills functional over this orbit must be contained in $\mathcal{G}(A_{\infty})$. In particular, every $\mu_{\tau}$-semistable orbit contains a unique $\mu_{\tau}$-polystable orbit in its closure. This yields the identification
						$$\mathcal{A}^{ss}(P)/\!/\mathcal{G}^c(P) \cong \mathcal{A}^{ps}/\mathcal{G}^c(P) \cong \mu_{\tau}^{-1}(0)/\mathcal{G}(P)$$
where two semistable orbits on the left hand side are identified if and only if they contain the same polystable orbit in their closure.

In the last section we extend this observation and characterize in Theorem \ref{ThmYMStab} the $\mu_{\tau}$-stability of $A \in \mathcal{A}(P)$ in terms of the limit of the Yang-Mills flow starting at $A$. We observe in particular that $\mathcal{A}^{ss}(P)$ and $\mathcal{A}^s(P)$ are both open subsets of $\mathcal{A}(P)$ in the $W^{1,2}$-topology.

\subsection{Analytical foundations}

\subsubsection{The Yang Mills flow on low dimensional manifolds}

Recall for $A \in \mathcal{A}(P)$ and $a \in W^{1,2}(\Sigma, T^*\Sigma\otimes \text{ad}(P))$ the formula
					$$F_{A+a} = F_A + d_A a + \frac{1}{2}[a \wedge a ].$$
From this follows directly that the $L^2$-gradient of the Yang-Mills functional (\ref{YMeq1}) is given by
		$$\nabla \mathcal{YM}: \mathcal{A}(P) \rightarrow W^{-1,2}(\Sigma, \text{ad}(P)), \qquad \nabla \mathcal{YM}(A) := d_{A}^* F_A.$$
The critical points of the Yang-Mills functional (\ref{YMeq1}) are called Yang-Mills connections and satisfy the equation
				$$d_A^* F_A = 0.$$
It follows from the strong Uhlenbeck compactness result (see e.g. \cite{Wehrheim:Ucpct}	Theorem E) and elliptic regularity that every Yang-Mills connection is gauge equivalent to a smooth Yang-Mills connection and the set $\Lambda := \{ \mathcal{YM}(A)\,|\, d_A^* F_A = 0 \}$ of critical values is discrete. The negative gradient flow of the Yang-Mills functional is given by the degenerated parabolic equation
			\begin{align} \label{YMeq101}\partial_t A(t) + d_{A(t)}^* F_{A(t)} = 0. \end{align}

\begin{Definition}[\textbf{Weak solutions}]
Let $A_0 \in \mathcal{A}(P)$ be a connection of Sobolev class $W^{1,2}$. We call $A \in C^0([0,\infty), \mathcal{A}(P))$ a weak solution of the initial value problem
			\begin{align} \label{YMeq2} \partial_t A(t) + d_{A(t)}^* F_{A(t)} = 0, \qquad A(0) = A_0 \end{align}
if $A(0) = A_0$ and there exists a sequence $A_k: [0,\infty) \rightarrow \mathcal{A}(P)$ of smooth solutions of (\ref{YMeq101}) which converges in $C^0_{loc}([0,\infty), \mathcal{A}(P))$ to $A$, where $\mathcal{A}(P)$ is endowed with the $W^{1,2}$-topology.
\end{Definition}

The next two theorems state that the initial value problem (\ref{YMeq2}) has a unique (weak) solution for every initial data $A_0 \in \mathcal{A}(P)$ existing for all time and that this solution converges to a Yang-Mills connection.

\begin{Theorem}[\textbf{Long time existence}] \label{ThmRade1}
Let $G$ be a compact connected Lie group, $P \rightarrow \Sigma$ a principal $G$ bundle and $A_0 \in \mathcal{A}(P)$.
		\begin{enumerate}
			\item There exists a unique weak solution $A(t) \in C_{loc}^0([0,\infty),\mathcal{A}(P))$ for the initial value problem (\ref{YMeq2}). The curvature has the additional regularity properties $F_{A(t)} \in C_{loc}^0([0,\infty), L^2)$ and $F_{A(t)} \in L_{loc}^2([0,\infty), W^{1,2})$.
			\item The solution $A(t)$ and its curvature $F_{A(t)}$ depend smoothly on the initial data $A_0$ in these topologies.
			\item If $A_0$ is smooth, then the solution $A(t)$ is smooth and satisfies (\ref{YMeq101}).
		\end{enumerate}
\end{Theorem}

\begin{proof}
This is Theorem 1 in \cite{Rade1992}.
\end{proof}

\begin{Theorem}[\textbf{Convergence}]\label{ThmRade2}
Assume the setting of Theorem \ref{ThmRade1} and let $A(t) \in C_{loc}^0([0,\infty),\mathcal{A}(P))$ be a weak solution of (\ref{YMeq2}). Then exists a Yang-Mills connection $A_{\infty} \in \mathcal{A}(P)$ and constants $c,\beta >0$ such that
						$$||A(t) - A_{\infty}||_{W^{1,2}} \leq c t^{-\beta}$$
holds for all times $t > 0$.
\end{Theorem}

\begin{proof}
This is Theorem 2 in \cite{Rade1992}.
\end{proof}

The key ingredient in the proof of the convergence result is the appropriate analogue of the Lojasiewicz gradient inequality. This approach was introduced by Simon \cite{Simon:1983} for a general class of evolution equations.

\begin{Proposition}[\textbf{Lojasiewicz gradient inequality}] \label{PropLoj}
Let $A_{\infty} \in \mathcal{A}(P)$ be a Yang-Mills connection. There exist constants $\epsilon > 0$, $\gamma \in [\frac{1}{2}, 1)$ and $c > 0$ such that for every $A \in \mathcal{A}(P)$ with $||A - A_{\infty}||_{W^{1,2}} < \epsilon$ the estimate
				$$||d_A^*F_A||_{W^{-1,2}} \geq c |\mathcal{YM}(A) - \mathcal{YM}(A_{\infty})|^{\gamma}$$
is satisfied.				
\end{Proposition}

\begin{proof}
This is Proposition 7.2 and (9.1) in \cite{Rade1992}. 
\end{proof}

\begin{Remark}
This estimate is slightly stronger then one might expect, as the $L^2$-norm is the natural norm for the gradient. Thus one would expect the estimate
		$$||d_A^*F_A||_{L^2} \geq c |\mathcal{YM}(A) - \mathcal{YM}(A_\infty)|^{\gamma}$$
which is an immediate consequence of Proposition \ref{PropLoj}.
\end{Remark}

In finite dimensions the Lojasiewicz inequality always guarantees convergence by some standard arguments. We recall these arguments in the following and discuss additional technical difficulties arising in the infinite dimensional setting. Suppose that $A(t)$ satisfies (\ref{YMeq2}). It follows from the weak Uhlenbeck compactness (see \cite{Rade1992}, Proposition 7.1) that there exists a $\mathcal{G}(P)$-orbit $\mathcal{G}(A_{\infty})$ of Yang-Mills connections such that
					$$\inf_{t>0} \mathcal{YM}(A(t)) = \mathcal{YM}(A_{\infty})$$
and for every $\delta > 0$ exists $T >0$ and $g\in\mathcal{G}(P)$ such that 
					$$||A(T) - g(A_{\infty})||_{W^{1,2}} < \delta.$$
Since the Yang-Mills functional and the Lojasiewicz inequality are invariant under the action of $\mathcal{G}(P)$, the constant $\epsilon = \epsilon(g(A_{\infty})) > 0$ from the Lojasiewicz inequality does not depend on $g$. Now choose $\delta < \epsilon$ and define
			$$\overline{T} := \inf \{ t > T \,|\, ||A(t) - g(A_{\infty})||_{W^{1,2}} \geq \epsilon \}.$$
For any $s_1, s_2 \in (T, \overline{T})$ with $s_1 < s_2$ we obtain
		\begin{align*}
				||A(s_1) - A(s_2)||_{L^2} &\leq \int_{s_1}^{s_2} ||d_A^*F_A||_{L_2}\, dt \\
																	&\leq \int_{s_1}^{s_2} \frac{||d_A^*F_A||_{L_2}^2}{c |\mathcal{YM}(A) - \mathcal{YM}(A_{\infty})|^{\gamma}} \, dt \\
																	&\leq \frac{1}{c} \int_{s_1}^{s_2} \frac{\partial}{\partial s}\left( \mathcal{YM}(A) - \mathcal{YM}(A_{\infty}) \right)^{1-\gamma} \, ds \\
																	&\leq \frac{1}{c} \left( \mathcal{YM}(A(s_1)) - \mathcal{YM}(A(s_2)) \right)^{1-\gamma}.
		\end{align*}
To conclude the convergence result, one needs to show $\overline{T} = \infty$ and extend the estimate above to the $W^{1,2}$-norm. Both can be achieved by using the following Lemma.

\begin{Lemma} \label{LemmaLoj}
Let $A_{\infty} \in \mathcal{A}(P)$ be a Yang-Mills connection and $\epsilon = \epsilon(A_{\infty}) > 0$ as in Proposition \ref{PropLoj}. There exists a constant $c>0$ with the following significance: For every solution $A(t)$ of the Yang-Mills flow (\ref{YMeq2}) and real numbers $0 \leq s_1 \leq s_2 - 1$ such that $||A(t) - A_{\infty}||_{W^{1,2}} \leq \epsilon$ for all $t \in [s_1,s_2]$ we have
					$$\int_{s_1+1}^{s_2} ||d_A^*F_A||_{W^{1,2}} \, dt \leq c \int_{s_1}^{s_2} ||d_A^*F_A||_{L^2}\, dt.$$
\end{Lemma}

\begin{proof} This is Lemma 7.3 in \cite{Rade1992}. \end{proof}

An immediate consequence of the Lemma and our calculation above is the estimate
			\begin{align} \label{eqLoj1} ||A(s_1 + 1) - A(s_2)||_{W^{1,2}} \leq C \left( \mathcal{YM}(A(s_1)) - \mathcal{YM}(A(s_2)) \right)^{1-\gamma} \end{align}
for any $T \leq s_1 < s_1 + 1 \leq s_2 < \bar{T}$. Since $\mathcal{YM}(A(t))$ is decreasing in $t$ and bounded below, it only remains to show $\bar{T} = \infty$. Since the solutions of the Yang-Mills flow depend continuously on the initial condition in the $C^0_{loc}([0,\infty), W^{1,2})$ topology, there exists a constant $c_1 > 0$ such that
				$$||A(T+t) - g(A_{\infty})||_{W^{1,2}} \leq c_1 ||A(T) - g(A_{\infty})||_{W^{1,2}}$$
holds for all $t \in [0,1]$. This follows as we may view $g(A_{\infty})$ as constant flow line and the constant $c_1$ depends only on the orbit $\mathcal{G}(A_{\infty})$ by $\mathcal{G}(P)$-invariance of the Yang-Mills functional. In particular, for sufficiently small $\delta$, we have $\delta c_1 < \epsilon$ and hence $\bar{T} > 1$. Then (\ref{eqLoj1}) yields
				$$||A(T+1)  - A(t)||_{W^{1,2}} \leq C\left( \mathcal{YM}(A(T)) - \mathcal{YM}(A(t)) \right)^{1-\gamma} \leq C \delta^{1-\gamma}$$
for any $T+1 \leq t \leq \bar{T}$. For sufficiently small $\delta > 0$ the right hand side is smaller than $\epsilon$ and this yields $\bar{T} = \infty$. This argument proves the following Corollary.

\begin{Cor}\label{CorLoj}
Let $B \in \mathcal{A}(P)$ be a Yang-Mills connection and let $\epsilon > 0$. Then there exists $\delta > 0$ such that for every solution $A(t)$ of the Yang-Mills flow (\ref{YMeq2}) with $||A(0) - B||_{W^{1,2}} < \delta$ we have either
				$$\sup_{t \geq 0} ||A(0) - A(t)||_{W^{1,2}} < \epsilon $$
for all $t \geq 0$ or there exists $T > 0$ with $\mathcal{YM}(A(T)) < \mathcal{YM}(B)$.
\end{Cor}

\subsubsection{The Kempf-Ness flow}

By Proposition \ref{CpxGaugeAction} the infinitesimal action of the complexified Gauge action is given by
						$$L_{A} (\xi + \textbf{i}\eta) := \left. \frac{d}{dt}\right|_{t=0} \exp(t\xi + \textbf{i}t\eta) A = -d_A \xi - * d_A \eta$$ 
for $\xi, \eta \in W^{2,2}(\Sigma,\text{ad}(P))$ and $A \in \mathcal{A}(P)$. With this formula we can express the gradient of the Yang-Mills functional as
						$$\nabla\mathcal{YM}(A) = d_A^*F_A = -* d_A *F_A  = L_A (\textbf{i}*F_A).$$
This implies that any solution of the Yang-Mills flow (\ref{YMeq1}) remains in a single complexified orbit.

\begin{Proposition} \label{YMProp1}
	Let $A_0 \in \mathcal{A}(P)$ and let $A(t)$ be the (weak) solution of the Yang-Mills flow (\ref{YMeq2}) starting at $A_0$. Let $g(t)$ be the solution of the pointwise ODE
						\begin{align} \label{YMeq3} g(t)^{-1} \dot{g}(t) = \textbf{i} (* F_{A(t)}), \qquad g(0) = \mathds{1}. \end{align}
	Then holds $g \in C_{loc}^0([0,\infty), \mathcal{G}^c(P))$ and 
									$$A(t) = g(t)^{-1}A_0$$
	for all $t \in [0,\infty)$. Moreover, $g$ depends continuously on $g_0$ and $A_0$.
\end{Proposition}

\begin{proof}
Recall from Lemma \ref{PreLem1} the formula
		$$B(t) := g_t^{-1}(A_0) = A_0 + g_t^{-1} d_{A_0} g_t - g^{-1}_t(h_t^{-1} \partial_{A_0} h_t) g_t$$
with $h_t := (g_t^{-1})^* g_t^{-1}$. By Theorem \ref{ThmRade1} holds $F \in L^2_{loc}([0,\infty), W^{1,2})$ and hence $g \in W^{1,2}_{loc}([0,\infty), W^{1,2})$ and $B \in W^{1,2}_{loc}([0,\infty), L^1)$. The same calculation as in the smooth case shows
		$$\dot{B}(t) =  L_{B(t)} (g_t^{-1}\dot{g_t}) = - d_{B(t)}^* F_{A(t)}.$$
Approximation of $A_0$ with smooth connections shows $A \in W^{1,2}_{loc}([0,\infty), L^1)$ and
		$$\dot{A}(t) = - d_{A(t)}^* F_{A(t)}.$$
Define $C(t) := A(t) - B(t)$ and $\Psi(t) := * F_{A(t)} \in L^2_{loc}([0,\infty), L^1)$. The calculation above shows that $C$ solves the linear ODE
			$$\dot{C}(t) = *[C(t), \Psi(t)], \qquad C(0) = 0.$$
The unique solution of this ODE is $C = 0$ and from the Sobolev embedding $W^{1,2}([0,t_0], L^1) \hookrightarrow C^0([0,t_0], L^1)$ we conclude $A(t) = B(t) = g_t^{-1} A_0$ for all $t$.
		
Since $A$ maps continuously in $W^{1,2}$, it follows from the expression
				$$A(t) = g_t^{-1} A_0 = A_0 + g_t^{-1} d_{A_0} g_t - g^{-1}_t(h_t^{-1} \partial_{A_0} h_t) g_t$$
that $A(t)^{0,1} = A_0^{0,1} + g^{-1}_t \bar{\partial}_{A_0} g_t$ and $g^{-1}_t \bar{\partial}_{A_0} g_t$ maps continuously into $W^{1,2}$. Let $\tilde{A}$ be a smooth reference metric and write $A_0 = \tilde{A} + a_0$. Then $g^{-1}_t\bar{\partial}_{\tilde{A}} g^{-1}_t$ maps continuously into $L^p$ for any $p < \infty$ and by elliptic regularity, $g$ maps continuously in $W^{1,p}$. Since $W^{1,p} \hookrightarrow C^0$, we can rerun the argument where $g^{-1}_t\bar{\partial}_{\tilde{A}} g_t$ now maps continuously in $W^{1,2}$ and conclude $g \in C^0_{loc}([0,\infty), W^{2,2})$. Since $A$ and $F_A$ depend continuously on $A_0$, the solution $g$ depends continuously on $A_0$ in $W^{1,2}_{loc}([0,\infty), W^{1,2})$ and then by elliptic regularity also in $C^0_{loc}([0,\infty), W^{2,2})$.
\end{proof}

\begin{Remark} \label{RmkYM1}
Let $A_0 \in \mathcal{A}(P)$ and let $A(t)$ be as in Proposition \ref{YMProp1}. For $g_0 \in \mathcal{G}^c(P)$ consider the more general equation
	\begin{align} \label{YMeq4}  g(t)^{-1} \dot{g}(t) = \textbf{i}* F_{A(t)}, \qquad g(0) = g_0. \end{align}
Then $\tilde{g}(t) = g_0^{-1} g(t)$ solves equation (\ref{YMeq3}) with respect to $\tilde{A}_0 = g_0^{-1}(A_0)$. Hence (\ref{YMeq4}) has a unique solution in $C^0_{loc}([0,\infty), \mathcal{G}^c(P))$, which depends continuously on $g_0$ and $A_0$.
\end{Remark}

We shall consider the following variant of this equation.

\begin{Definition}[\textbf{Kempf-Ness flow}] \label{DefnWeak2}
Let $A_0 \in \mathcal{A}(P)$ and $g_0 \in \mathcal{G}^c(P)$. We say that $g(t) \in C_{loc}^0([0,\infty), \mathcal{G}^c(P))$ is a weak solution of the equation
					\begin{align} \label{YMeq6} g^{-1}(t)\dot{g}(t) = \textbf{i} * F_{g^{-1}(t) A_0}, \qquad g(0) = g_0 \end{align}
if there exist a sequence of smooth initial data $(A_k, g_0^k) \in \mathcal{A}(P)\times \mathcal{G}^c(P)$ converging to $(A_0, g_0)$ and smooth solutions $g_k(t)$ of the equation
					$$g_k^{-1}(t)\dot{g}_k(t) = \textbf{i} * F_{g_k^{-1}(t) A_k}, \qquad g_k(0) = g_0^k$$
such that $g_k(t)$ converges to $g(t)$ in $C_{loc}^0([0,\infty), W^{2,2})$. 
\end{Definition}

\begin{Remark}
We call a solution $g \in C_{loc}^0([0,\infty, \mathcal{G}^c(P))$ of (\ref{YMeq4}) a solution of the Kempf-Ness flow starting at $g_0$ (with respect to $A_0$). We show in Section 6 that there exists a $\mathcal{G}(P)$-invariant functional 
		$$\Phi_{A_0}: \mathcal{G}^c(P) \rightarrow \mathbb{R}$$
whose negative gradient flow lines correspond to solution of (\ref{YMeq6}).
\end{Remark}

\begin{Lemma}
For every initial data $(A_0, g_0) \in \mathcal{A}(P)\times\mathcal{G}^c(P)$ there exists a unique (weak) solution of (\ref{YMeq6}) which depends continuously on the initial data.
\end{Lemma}

\begin{proof}
We use the notation introduced in Definition \ref{DefnWeak2} and define $A_k(t) := g_k(t)^{-1}A_k$. Then
			$$\partial_t A_k(t) = L_{A_k(t)} \left(g_k^{-1}(t)\dot{g}_k(t) \right) = - L_{A_k(t)} \left(\textbf{i} * F_{A_k(t)}\right) = - d_{A_k(t)}^* F_{A_k(t)}$$
and thus $A_k(t)$ yields a smooth solution of the Yang-Mills flow. Conversely, the solution $A_k(t)$ is uniquely determined by the initial condition $(g_0^k)^{-1}A_k$ and we may recover $g_k(t)$ from this solution via Proposition \ref{YMProp1} and Remark \ref{RmkYM1}. Since solutions of the Yang-Mills flow and solutions of (\ref{YMeq4}) depend continuously on the initial data, it follows that the weak solution $g(t)$ of (\ref{YMeq6}) is uniquely determined by the weak solution $A(t)$ of the Yang-Mills flow starting at $g_0^{-1}A_0$.
\end{proof}

The next Proposition shows that solutions of the Kempf-Ness flow (\ref{YMeq6}) remain at bounded distance in the homogeneous space $\mathcal{G}^c/\mathcal{G}$.

\begin{Proposition}\label{YMProp2}
Let $A_0 \in \mathcal{A}(P)$ and let $g,\tilde{g} \in C_{loc}^0( [0,\infty), \mathcal{G}^c(P))$ be (weak) solutions of (\ref{YMeq6}) starting at $g_0, \tilde{g}_0 \in \mathcal{G}^c(P)$. Define $\eta(t) \in W^{2,2}(\Sigma,\text{ad}(P))$ and $u(t) \in \mathcal{G}(P)$ by the equation
		$$g(t) \exp(\textbf{i} \eta(t)) u(t) = \tilde{g}(t).$$
Then the following holds:
		\begin{enumerate}
			\item[\textbf{(i)}] $\rho(t) := ||\eta(t)||_{L^2}$ is non-increasing in $t$. More precisely, if $\eta(t) \neq 0$ then
					$$\dot{\rho}(t) = - \frac{1}{\rho(t)} \int_0^1 ||d_{A_{s,t}} \eta(t)||_{L^2}^2\, ds$$
				with $A_{s,t} := e^{-\textbf{i}s\eta(t)}g_t^{-1} A_0$.
			\item[\textbf{(ii)}] The differential inequality
								$$(\partial_t + \Delta)|\eta|^2 \leq 0$$
						is satisfied. In particular, $||\eta(t)||_{L^{\infty}}$ is non-increasing by the maximum principle for the heat equation 
						
			\item[\textbf{(iii)}] $\eta$ is uniformly bounded in $W^{2,2}$.
			
			\item[\textbf{(iv)}] $u$ is uniformly bounded in $W^{2,2}$.
		\end{enumerate}
\end{Proposition}

\begin{proof}
We prove (i) and (ii): By approximation, we can assume that $A_0$, $g$ and $\tilde{g}$ are all smooth. Let $\pi: G^c \rightarrow G^c/G$ denote the projection and define $\gamma(s,t) := \pi(g(t)e^{\textbf{i}s \eta(t)})$. Pointwise $\gamma(\cdot,t)$ is the unique geodesic of length $|\eta(t)|$ connecting $\pi(g)$ and $\pi(\tilde{g})$. The following calculation is pointwise valid:
		\begin{align*}
			\partial_t |\eta|^2 
										&= \partial_t \int_0^1 \langle \partial_s \gamma, \partial_s \gamma \rangle \, ds 
													= 2\int_0^1 \langle \nabla_t \partial_s \gamma, \partial_s \gamma \rangle  \,  ds  \\
										&= 2\int_0^1 \langle \nabla_s \partial_t \gamma, \partial_s \gamma \rangle  \,  ds
													= 2\int_0^1 \partial_s\langle \partial_t \gamma, \partial_s \gamma \rangle \,  ds		\\
										&= 2\left(\langle \partial_t \gamma(1,t), \partial_s \gamma(1,t) \rangle  - \langle \partial_t \gamma(0,t), \partial_s \gamma(0,t) \rangle \right) \\
										&= 2\langle \tilde{g}^{-1}(t) \dot{\tilde{g}}(t) - g^{-1}(t) \dot{g}(t), \textbf{i}\eta(t) \rangle \\
										&= 2\langle *F_{\tilde{g}(t)^{-1} A_0} - * F_{g(t)^{-1} A_0}, \eta(t) \rangle
			\end{align*}
With $A_{s,t} := e^{-\textbf{i}s\eta(t)}g^{-1} A_0$ this yields
		\begin{align*}
			\partial_t |\eta|^2	
										&= 2\int_0^1 \langle \eta(t), *d_{A_{s,t}} * d_{A_{s,t}} \eta(t) \rangle 
										= - \Delta |\eta|^2 - 2 \int_0^1 |d_{A_{s,t}} \eta(t)|^2\, ds.
	\end{align*}
This proves the second claim and the first one is obtained by integrating this inequality over $\Sigma$. \\

We prove (iii) and (iv): Recall that $\tilde{A}(t) := \tilde{g}^{-1}_t(A_0)$ and $A(t):= g^{-1}_t(A_0)$ are solutions of the Yang-Mills flow. Since they converge in $W^{1,2}$, they are both uniformly bounded in $W^{1,2}$. With $a(t) := e^{\textbf{i}\eta_t}u_t$ we have $\tilde{A} = a^{-1}(A)$ and hence
		\begin{align*}
				\tilde{A}^{0,1} = A^{0,1} + a^{-1} \bar{\partial}_A a
		\end{align*}
This shows that $\bar{\partial}_A a$ is uniformly bounded in $W^{1,2}$. Via elliptic bootstrapping follows that $a$ is uniformly bounded in $W^{2,2}$. From the formula $aa^* = e^{2\textbf{i}\eta}$ we conclude that $\eta$ is uniformly bounded in $W^{2,2}$ and then $u$ is also uniformly bounded in $W^{2,2}$. \\
\end{proof}

\subsection{Uniqueness of Yang-Mills connections}

We follow the arguments from (\cite{RobSaGeo}, Chapter 6) to proof the analog of the Ness uniqueness theorem and the moment limit theorem. These are originally due to Calabi-Chen \cite{CalabiChen:2002} and Chen-Sun \cite{ChenSun:2010} in the context of extremal Kähler metrics.

\begin{Proposition}\label{PropNUrest}
Let $A_0, A_1 \in \mathcal{A}(P)$ be Yang-Mills connections with $\mathcal{G}^c(A_0) = \mathcal{G}^c(A_1)$. Then holds $\mathcal{G}(A_0) = \mathcal{G}(A_1)$.
\end{Proposition}

This is a special case of Theorem \ref{ThmNUgen} below. The proof is essentially due to Calabi and Chen (\cite{CalabiChen:2002}, Corolary 4.1).

\begin{proof}
Choose $\tilde{g} \in \mathcal{G}^c(P)$ such that 
						$$A_1 = \tilde{g}^{-1}A_0.$$
holds. Since $A_0$ and $A_1$ are Yang-Mills connections, they generate constant flow lines $A_0(t) \equiv A_0$ and $A_1(t) \equiv A_1$. Let $g_0,g_1 \in C^0_{loc}([0,\infty), \mathcal{G}^c)$ be the solutions of the equation
			$$g_0^{-1}\dot{g}_0 = *F_{A_0}, \quad g_0(0) = \mathds{1} \quad \text{and} \quad g_1^{-1}\dot{g}_1 = *F_{A_1}, \quad g_1(0) = \tilde{g}$$
from Proposition \ref{YMProp1} and the following Remark. They satisfy 
			$$A_0 = g_0^{-1}(t)A_0 \quad \text{and} \quad A_1 = g_1^{-1}(t)A_0.$$
In particular, $g_0$ and $g_1$ are solutions of the Kempf-Ness flow (\ref{YMeq6}) with respect to $A_0$. Define $\eta(t) \in W^{2,2}(\Sigma, \text{ad}(P))$ and $u(t) \in \mathcal{G}$ by the equation
		$$g_0(t) \exp(\textbf{i} \eta(t)) u(t) = g_1(t)$$
as in Proposition \ref{YMProp2}. Using the notation from the Proposition, there exist $\eta_{\infty} \in W^{2,2}(\Sigma,\text{ad}(P))$, $u_{\infty} \in \mathcal{G}(P)$ and a sequence $t_i \rightarrow \infty$ such that
		$$\lim_{i\rightarrow \infty} \dot{\rho}(t_i) = 0, \qquad \eta(t_i) \stackrel{H^2}{\rightharpoonup} \eta_{\infty}, \qquad u(t_i) \stackrel{H^2}{\rightharpoonup} u_{\infty}.$$
By taking a further subsequence if necessary, we may assume that
		$$\lim_{i\rightarrow \infty} ||d_{A_{s,t_i}} \eta(t_i)||_{L^2} = 0$$
holds for almost every $s \in [0,1]$, where we defined 
				$$A_{s,t} = e^{-\textbf{i}s \eta(t)}(g_0^{-1}(t)A_0) = e^{-\textbf{i}s \eta(t)}A_0.$$
Moreover, by Rellich's theorem, $\eta(t_i)$ and $u(t_i)$ converge for every $p < \infty$ strongly in $W^{1,p}$ to $\eta_{\infty}$ and $u_{\infty}$. By continuity of the Gauge action $\mathcal{G}^{1,p} \times \mathcal{A}^p \rightarrow \mathcal{A}^p$ for $p > 2$, we conclude
			$$A_{s,t_i} \stackrel{L^p}{\rightarrow} A_{s,\infty} := e^{-\textbf{i}s \eta_{\infty}}A_0, \quad \text{and} \quad d_{A_{s, {t_i}}} \eta(t_i) \stackrel{L^p}{\rightarrow} d_{A_{s,\infty}} \eta_{\infty}.$$ 
This implies that for almost every $s \in [0,1]$, we must have $d_{A_{s,\infty}} \eta_{\infty} = 0$. For $s \rightarrow 0$ we conclude $d_{A_{0,\infty}} \eta_{\infty} = d_{A_0} \eta_{\infty} = 0$ and hence $e^{-\textbf{i}\eta_{\infty}}A_0 = A_0$. It follows now
		$$A_1 = g_1(t_i)^{-1}A_0 = u(t_i)^{-1} e^{-\textbf{i}\eta(t_i)} A_0 \stackrel{L^p}{\longrightarrow} u_{\infty}^{-1} e^{-\textbf{i}\eta_{\infty}}A_0 = u_{\infty}^{-1} A_0.$$
This shows $A_1 = u_{\infty}^{-1} A_0$ and thus $A_0$ and $A_1$ lie in the same $\mathcal{G}$-orbit. 
\end{proof}

\begin{Theorem}[\textbf{Moment Limit Theorem}] \label{MLT}
Let $A_0 \in \mathcal{A}(P)$ and $A: [0,\infty) \rightarrow \mathcal{A}(P)$ be the solution of the Yang-Mills flow starting at $A_0$. The limit $A_{\infty} := \lim_{t\rightarrow \infty} A(t)$ satisfies
			$$\mathcal{YM}(A_{\infty}) = \inf_{g \in \mathcal{G}^c(P)} \mathcal{YM}(g A_0).$$
Moreover, the $\mathcal{G}(P)$-orbit of $A_{\infty}$ depends only on the complexified orbit $\mathcal{G}^c(A_0)$.
\end{Theorem}

\begin{proof}
Let $g_0 \in \mathcal{G}^c(P)$ be given and define $g, \tilde{g} \in C^0_{loc}([0,\infty), \mathcal{G}^c)$ by
		$$g^{-1}\dot{g} = *F_{A}, \quad g(0) = \mathds{1} \quad \text{and} \quad \tilde{g}^{-1}\dot{\tilde{g}} = *F_{A} , \quad g(0) = g_0$$
as in Proposition \ref{YMProp1} and the following Remark. Let $A(t)$ and $\tilde{A}(t)$ be the solutions of the Yang-Mills flow starting at $A_0$ and $\tilde{A}_0 = g_0^{-1}A_0$. Then holds
			$$A_0(t) = g_t^{-1}(A_0), \qquad \tilde{A}(t) = \tilde{g}_t^{-1}(A_0)$$
and $g, \tilde{g}$ are solutions of the Kempf-Ness flow (\ref{YMeq6}) with respect to $A_0$. Define $\eta(t) \in W^{2,2}(\Sigma, \text{ad}(P))$ and $u(t) \in \mathcal{G}(P)$ by the equation
		$$g_0(t) \exp(\textbf{i} \eta(t)) u(t) = g_1(t)$$
as in Proposition \ref{YMProp2}. It follows that there exist $\eta_{\infty} \in W^{2,2}(\Sigma,\text{ad}(P))$ and $u_{\infty} \in \mathcal{G}(P)$ and a sequence $t_i \rightarrow \infty$ such that
		$$\eta(t_i) \stackrel{W^{2,2}}{\rightharpoonup} \eta_{\infty}, \qquad u(t_i) \stackrel{W^{2,2}}{\rightharpoonup} u_{\infty}.$$
By Rellich’s theorem we obtain strong convergence in $W^{1,p}$ and using the Sobolev embedding $W^{1,2} \hookrightarrow L^p$ for every $p < \infty$ we obtain:
	$$\tilde{A}_{\infty} \stackrel{W^{1,2}}{\longleftarrow} \tilde{A}(t_i) = u(t_i)^{-1} e^{\textbf{i}\eta(t_i)}A(t_i) \stackrel{L^p}{\longrightarrow} u_{\infty}^{-1} \eta_{\infty}^{-1} A_{\infty}.$$
Hence	$\tilde{A}_{\infty} = u_{\infty}^{-1} \eta_{\infty}^{-1} A_{\infty}$. Thus $\tilde{A}_{\infty}$ and $A_{\infty}$ are Yang-Mills connections lying in a common complexified orbit and Proposition \ref{PropNUrest} shows that in fact $\mathcal{G}(A_\infty) = \mathcal{G}(\tilde{A}_\infty)$. This shows $\mathcal{YM}(A_{\infty})= \mathcal{YM}(\tilde{A}_{\infty}) \leq \mathcal{YM}(g_0^{-1} A_0)$ and completes the proof.
\end{proof}

The following theorem is the analog of the Ness uniqueness theorem in finite dimensional GIT.

\begin{Theorem}[\textbf{Uniqueness of Yang-Mills connections}]\label{ThmNUgen}
Let $A_0 \in \mathcal{A}(P)$ and $A',A'' \in \overline{\mathcal{G}^c(A_0)}$ be in the $W^{1,2}$-closure of a single complexified orbit satisfying
			$$\mathcal{YM}(A') = \mathcal{YM}(A'') = \inf_{g \in \mathcal{G}^c} \mathcal{YM}(g A_0).$$
Then follows $\mathcal{G}(A') = \mathcal{G}(A'')$.
\end{Theorem}

\begin{Cor}
Let $P \rightarrow \Sigma$ be a principal $G$ bundle of constant central type $\tau \in Z(\mathfrak{g})$ defined by $(\ref{eqCT})$. Suppose $A \in \mathcal{A}(P)$ is $\mu_{\tau}$-semistable. Then the $W^{1,2}$-closure $\overline{\mathcal{G}^c(A)}$ contains a unique $\mu_{\tau}$-polystable orbit.
\end{Cor}

\begin{proof}
It follows from (\ref{YMeq02}) that solutions of the equation $*F_A = \tau$ correspond to global minima of the Yang-Mills functional on $\mathcal{A}(P)$. The Corollary follows thus from Theorem \ref{ThmNUgen}.
\end{proof}

The following proof of Theorem \ref{ThmNUgen} is essentially due to Chen and Sun (\cite{ChenSun:2010}, Theorem 4.1).

\begin{proof}[Proof of Theorem \ref{ThmNUgen}]
Let $A(t)$ be the solution of the Yang-Mills flow starting at $A_0$ and let $A_{\infty} := \lim_{t\rightarrow \infty} A(t)$. Then Theorem \ref{MLT} implies
				$$\mathcal{YM}(A_{\infty}) = \inf_{g \in \mathcal{G}^c} \mathcal{YM}(g A_0) =: m.$$
Since $A_{\infty} \in \overline{\mathcal{G}^c(A_0)}$, it suffices to show that any connection $B \in \overline{\mathcal{G}^c(A_0)}$ with $\mathcal{YM}(B) = m$ is contained in $\mathcal{G}(A_{\infty})$. For this let $A_i \in \mathcal{G}^c(A_0)$ be a sequence which converges to $B$ and denote by $A_i(t)$ the corresponding solutions of the Yang-Mills flow and set $B_i := \lim_{t\rightarrow \infty} A_i(t)$. Note that $B$ is necessarily a Yang-Mills connection, since
						$$\mathcal{YM}(B(t)) = \lim_{i\rightarrow \infty} \mathcal{YM}(A_i(t)) \geq m = \mathcal{YM}(B(0))$$
where we denote by $B(t)$ the solution of the Yang-Mills flow starting at $B$. Thus, we may apply Corollary \ref{CorLoj} with respect to $B$ and conclude that $||A_i - B_i||_{W^{1,2}}$ converges to zero and hence
					$$\lim_{i\rightarrow \infty} B_i = B.$$
By Theorem \ref{MLT} holds $\mathcal{G}(B_i) = \mathcal{G}(A_{\infty})$ and hence there exists $u_i \in \mathcal{G}(P)$ such that $u_i^{-1}(A_{\infty}) = B_{i}$. Since the connections $B_i$ are uniformly bounded in $W^{1,2}$, the gauge transformations $u_i$ are uniformly bounded in $W^{2,2}$. Thus there exists $u_{\infty}\in\mathcal{G}(P)$ such that after passing to a subsequence $u_i$ converges weakly in $W^{2,2}$ to $u_{\infty}$ and strongly in $W^{1,p}$ for any $p < \infty$. Using the continuity of the Gauge action $\mathcal{G}^{1,p}\times\mathcal{A}^p \rightarrow \mathcal{A}^p$ we conclude
					$$B_i = u_i^{-1}(A_{\infty}) \stackrel{L^p}{\longrightarrow} u_{\infty}^{-1}A_{\infty}$$
and in particular $B = u_{\infty}^{-1}A_{\infty} \in \mathcal{G}(A_{\infty})$

\end{proof}

\subsection{Yang-Mills characterization of $\mu_{\tau}$-stability}

We characterize the $\mu_{\tau}$-stability of a connection $A \in \mathcal{A}(P)$ in terms of the the limit $A_{\infty}$ of the Yang-Mills flow starting at $A$. This is Theorem \ref{ThmYMStab} below. The proof relies on the following proposition.

\begin{Proposition} \label{PropStabilityOpen}
Let $P \rightarrow \Sigma$ be a principal $G$ bundle of central type $\tau \in Z(\mathfrak{g})$ defined by (\ref{eqCT}). The subsets of $\mu_{\tau}$-semistable connections
			$$\mathcal{A}^{ss}(P) := \left\{ A \in \mathcal{A}(P)\,|\, \text{$A$ is $\mu_{\tau}$-semistable} \right\}$$
and $\mu_{\tau}$-stable connections
			$$\mathcal{A}^{s}(P) := \left\{ A \in \mathcal{A}(P)\,|\, \text{$A$ is $\mu_{\tau}$-stable} \right\}$$
			are open subsets of $\mathcal{A}(P)$ with respect to the $W^{1,2}$-topology.
\end{Proposition}

\begin{proof}
It follows from (\ref{YMeq02}) that
							$$\inf_{A \in \mathcal{A}(P)} \mathcal{YM}(A) \geq \frac{1}{2}||\tau||^2 =: m.$$
Moreover
			\begin{align} 
					\label{YMeq31}\mathcal{A}^{ss}(P) := \left\{ A \in \mathcal{A}(P)\,\left|\, \inf_{g\in \mathcal{G}^c(P)} \mathcal{YM}(gA) = m  \right.\right\}
			\end{align}
and $\mathcal{YM}(A) = m$ is equivalent to $*F_A = \tau$.\\

\textbf{Step 1:} \textit{$\mathcal{A}^{ss}(P)$ is open.} \\

Let $A_0 \in \mathcal{A}^{ss}(P)$ be given. Let $A(t)$ be the solution of the Yang-Mills flow starting at $A_0$ and $A_{\infty} := \lim_{t\rightarrow\infty} A(t)$. It follows from Theorem \ref{MLT} and (\ref{YMeq31}) that $A_{\infty}$ is a Yang-Mills connection satisfying $\mathcal{YM}(A_{\infty}) = m$. By the Lojasiewicz inequality (Proposition \ref{PropLoj}) there exists $\epsilon > 0$, $c > 0$ and $\gamma \in [\frac{1}{2}, 1)$ such that for all $B \in \mathcal{A}(P)$ with $||B - A_{\infty}||_{W^{1,2}} < \epsilon$ the inequality
				\begin{align} \label{YMapp2eq1} ||d_B^* F_B||_{L^2} \geq c ||\mathcal{YM}(B) - m||^{\gamma} \end{align}
is satisfied. By Corollary \ref{CorLoj} there exists $\delta > 0$ such that for every $B\in\mathcal{A}(P)$ with $||B - A_{\infty}||_{W^{1,2}} < \delta$ we have $||B_{\infty} - A_{\infty}||_{W^{1,2}} < \epsilon$. In particular,  (\ref{YMapp2eq1}) applies to $B_{\infty}$ and yields $\mathcal{YM}(B_{\infty}) = m$. This shows
			$$U := \{ B \in \mathcal{A}(P)\, |\, ||B - A_{\infty}||_{W^{1,2}} < \delta \} \subset \mathcal{A}^{ss}(P).$$
Now choose $T > 0$ such that $A(T) \in U$ and choose $g \in \mathcal{G}^c(P)$ with $A(T) = g^{-1}A_0$. By continuity of the gauge action there exists an open neighborhood $V$ of $A_0$ with $g^{-1}V \subset U$ and hence $V \subset \mathcal{A}^{ss}(P)$.\\

\textbf{Step 2:} \textit{Denote by $\mathcal{A}^*(P) \subset \mathcal{A}(P)$ the space of irreducible connections. This is an open subset and 
					$$\mathcal{Z} := \{A \in \mathcal{A}^*\,|\, \mathcal{YM}(A) = m \}/\mathcal{G}$$
	is a finite dimensional smooth submanifold of $\mathcal{A}^*/\mathcal{G}$.}\\

We may assume that $Z(G)$ is discrete, $\tau = 0$ and $m=0$, since otherwise $\mathcal{A}^*(P) = \emptyset$. Let $A_0 \in \mathcal{A}^*(P)$ be a smooth irreducible connection. The Laplacian $d_{A_0}^*d_{A_0}$ is then injective and by elliptic regularity there exists $c_0 > 0$ such that
			$$||d_{A_0}^*d_{A_0} \xi ||_{L^2} \geq c_0 ||\xi||_{W^{2,2}}$$
for all $\xi \in W^{2,2}(\Sigma, \text{ad}(P))$. For $a \in W^{1,2}(\Sigma, T^*\Sigma\otimes\text{ad}(P))$ expand
			$$d_{A_0+a}^*d_{A_0+a} \xi = d_{A_0}^* d_{A_0}\xi + d_{A_0}^* [a,\xi] - * [a,*d_{A_0}\xi] - * [a, *[a,\xi]].$$
Since $\text{dim}_{\mathbb{R}}(\Sigma) = 2$, we have the Sobolev estimate $||fg||_{L^2} \leq c ||f||_{W^{1,2}}||g||_{W^{1,2}}$ and $||fg||_{W^{1,2}} \leq ||f||_{W^{1,2}} ||g||_{W^{2,2}}$. This yields
			$$||d_{A_0+a}^*d_{A_0+a} \xi|| \geq c_0||\xi||_{W^{2,2}}  - c ||a||_{W^{1,2}} ||\xi||_{W^{2,2}}$$
and $A_0 + a$ is irreducible if $||a||_{W^{1,2}}$ is sufficiently small. Hence $\mathcal{A}^*(P)$ is open.

Now fix an irreducible connection $A_0$ with $*F_{A_0} = 0$. We may assume without loss of generality that $A_0$ is smooth and work in a Coulomb gauge relative to $A_0$. This allows us to identify a neighborhood of $[A_0]$ in $\mathcal{A}^*(P)/\mathcal{G}(P)$ with $a \in W^{1,2}(\Sigma, T^*\Sigma\otimes \text{ad}(P))$ satisfying $||a||_{W^{1,2}} < \epsilon$ and $d_{A_0}^* a = 0$ under the map $a \mapsto [A_0 + a]$. Consider
		$$\phi: \{a \in W^{1,2}(\Sigma, T^*\Sigma\otimes \text{ad}(P))\,|\, d_{A_0}^* a = 0, ||a||_{W^{1,2}} < \epsilon\} \rightarrow L^2(\Sigma,\text{ad}(P))$$
					$$\phi(a) := * F_{A_0 +a} $$
and define $Z_{A_0} := \phi^{-1}(0)$.  We claim that $0$ is a regular value for $\phi$ (after possibly shrinking $\epsilon$). Once this is established, the claim follows from the implicit function theorem. The derivative of $\phi$ at a point $a$ is given by
		$$d\phi(a): \{\hat{a} \in W^{1,2}(\Sigma, T^*\Sigma\otimes \text{ad}(P))\,|\, d_{A_0}^* \hat{a} = 0\} \rightarrow  L^2(\Sigma, \text{ad}(P))$$
					$$d\phi(a)\hat{a} = *d_{A_0} \hat{a} + *[a \wedge \hat{a}].$$
Since $d\phi(a)$ is the restriction of a compact perturbation of the Fredholm operator $*(d_{A_0}\oplus d_{A_0}^*)$, its kernel is finite dimensional. We denote by
			$$K := \{\hat{a} \in W^{1,2}(\Sigma, T^*\Sigma\otimes \text{ad}(P)) \,|\, d_{A_0} \hat{a} = 0, d_{A_0}^*\hat{a} = 0\}$$
the space of $A_0$-harmonic $1$-forms with values in $\text{ad}(P)$ and define $V$ by the $L^2$-orthogonal decomposition
			$$W^{1,2}(\Sigma, T^*\Sigma\otimes \text{ad}(P)) = V \oplus K.$$
Then the restriction of the Fredholm-operator $d_{A_0}\oplus d_{A_0}^*$ to $V$ defines an isomorphism
			$$d_{A_0}\oplus d_{A_0}^*: V \rightarrow L^2(\Sigma, \text{ad}(P)) \oplus L^2(\Sigma, \Lambda^2T^*\Sigma\otimes \text{ad}(P))$$
It is injective by definition of $V$ and to prove surjectivity let $f \in L^2(\Sigma, \text{ad}(P))$ and $\omega \in L^2(\Sigma, \Lambda^2T^*\Sigma\otimes \text{ad}(P))$ be given. Then by Hodge theory we can solve the equation 
			$$\Delta_{A_0} \hat{a} = d_{A_0}^* \omega + d_{A_0} f.$$
From this follows
			$$d_{A_0}^*(d_{A_0} \hat{a} - \omega) = d_{A_0} (f- d_{A_0}^* \hat{a}).$$
Since $*F_{A_0} = 0$, both sides of the equation are orthogonal and hence must vanish. Since $A_0$ is irreducible, it follows $d_{A_0} \hat{a} = \omega$ and $d_{A_0}^*\hat{a} = f$. In particular, for any $s \in L^2(\Sigma, \text{ad}(P))$ exists a solution $\hat{a} \in V$ of the equations
				\begin{align} \label{stabeq1} d_{A_0} \hat{a} + [a\wedge \hat{a}]= *s, \qquad d_{A_0}^* \hat{a} = 0 \end{align}
for $a = 0$. Since the equation is linear in $a$, another application of the inverse function theorem shows that after possibly shrinking $\epsilon$ the equation (\ref{stabeq1}) has a solution $\hat{a}(a) \in V$ for all $a$ with $||a||_{W^{1,2}} < \epsilon$.\\

\textbf{Step 3:} \textit{$\mathcal{A}^s$ is open.}\\

We may assume that $Z(G)$ is discrete, $\tau = 0$ and $m=0$, since otherwise $\mathcal{A}^s(P) = \emptyset$. Let $A \in \mathcal{A}^s(P)$ be given. By definition there exists $g \in \mathcal{G}^c(P)$ such that $A_0 = g^{-1}A$ is smooth and satisfies $\mathcal{YM}(A_0) = 0$. Let $Z_{A_0}$ be as in Step 2 and consider the map
		$$ \psi: Z_{A_0}\times W^{2,2}(\Sigma, \text{ad}(P))\times W^{2,2}(\Sigma,\text{ad}(P)) \rightarrow \mathcal{A}$$
					$$\psi(A,\xi,\eta) := e^{\textbf{i} \eta} e^{\xi} A.$$
We have seen that $Z_{A_0}$ is a smooth manifold with tangent space
				$$T_{A_0} Z_{A_0} = \{ \hat{a} \in W^{1,2}(\Sigma, \text{ad}(P))\, |\, d_{A_0}^* \hat{a} = 0, \, d_{A_0} \hat{a} = 0\}.$$
The differential of $\psi$ at the point $(A_0, 0,0)$ is given by
			$$d\psi(A_0,0,0) [\hat{a}, \hat{\xi}, \hat{\eta}] := \hat{a} - d_{A_0} \hat{\xi} - * d_{A_0} \hat{\eta}.$$
Since $F_{A_0} = 0$, it follows as in Step 2 from Hodge theory that $d\psi(A_0,0,0)$ is an isomorphism. The implicit function theorem yields thus an open neighborhood $U$ of $A_0$ with
			$$A_0 \in U \subset \text{Im}(\psi) \subset \mathcal{A}^{s}.$$
Finally, by continuity of the gauge action, there exists an open neighborhood $V$ of $A$ with $g^{-1}V \subset U$ and hence $\mathcal{A}^s(P)$ is open.

\end{proof}

\begin{Theorem} \label{ThmYMStab}
Let $P \rightarrow \Sigma$ be a principal $G$ bundle of central type $\tau \in Z(\mathfrak{g})$ defined by (\ref{eqCT}) and denote $m := \frac{1}{2}||\tau||^2$. Let $A_0 \in \mathcal{A}(P)$ and denote by $A_{\infty}$ the limit of the the Yang-Mills flow $A(t)$ starting at $A_0$.
	\begin{enumerate}
		\item $A_0$ is $\mu_{\tau}$-stable if and only if $A_{\infty}$ is irreducible.
		\item $A_0$ is $\mu_{\tau}$-polystable if and only if $\mathcal{YM}(A_{\infty}) = m$ and $A_{\infty} \in \mathcal{G}^c(A_0)$.
		\item $A_0$ is $\mu_{\tau}$-semistable if and only if $\mathcal{YM}(A_{\infty}) = m$.
		\item $A_0$ is $\mu_{\tau}$-unstable if and only if $\mathcal{YM}(A_{\infty}) > m$.
	\end{enumerate}
\end{Theorem}

\begin{proof}
It follows from (\ref{YMeq02}) that $m$ is a lower bound for the Yang-Mills functional on $\mathcal{A}(P)$ and $A \in \mathcal{A}(P)$ satisfies $\mathcal{YM}(A) = m$ if and only if $*F_A = \tau$. Thus the characterization for $\mu_{\tau}$-unstable and $\mu$-semistable orbits follows from Theorem \ref{MLT}.

Suppose next that $A_0$ is $\mu_{\tau}$-polystable. Then exists $g_0 \in \mathcal{G}^c(P)$ such that $\tilde{A}_0 := g_0^{-1}(A_0)$ satisfies $\mathcal{YM}(\tilde{A}_0) = m$. The Yang-Mills flow line $\tilde{A}(t)$ starting at $\tilde{A}_0$ is constant and it follows from Theorem \ref{MLT} and Theorem \ref{ThmNUgen} that $A_{\infty} \in \mathcal{G}(\tilde{A}_0) \subset \mathcal{G}^c(A_0)$. The converse is immediate and this proves the criterion for $\mu_{\tau}$-polystable orbits.

Suppose now that $A_0$ is $\mu_{\tau}$-stable. Then the orbit $\mathcal{G}^c(A_0)$ has only discrete $\mathcal{G}^c(P)$ isotropy. Since $A_0$ is in particular $\mu_{\tau}$-polystable, we have $A_{\infty} \in \mathcal{G}^c(A_0)$. Hence the infinitesimal action $L_{A_{\infty}}: \xi \mapsto - d_{A_{\infty}} \xi$ is injective and $A_{\infty}$ is irreducible. Suppose conversely that $A_{\infty}$ is irreducible. Since $A_{\infty}$ is a Yang-Mills connection, it satisfies $d_{A_{\infty}} *F_{A_{\infty}} = 0$ and hence $F_{A_{\infty}} = 0$. This shows that $\mathcal{G}^c(A_{\infty})$ is stable. By the Proposition \ref{PropStabilityOpen}, the subset $\mathcal{A}^s(P)$ of $\mu_{\tau}$-stable connections is open and hence $A(t) \in \mathcal{A}^s(P)$ for all sufficiently large $t$. Since the notion of $\mu_{\tau}$-stability is $\mathcal{G}^c(P)$-invariant, and since $A(t) \in \mathcal{G}^c(A_0)$, we conclude that $A_0$ is $\mu_{\tau}$-stable. 
\end{proof}


\section{Maximal weights}

Let $G$ be a compact connected Lie group, let $P \rightarrow \Sigma$ be a principal $G$ bundle and let $\tau \in Z(\mathfrak{g})$ denote the central type of $P$ defined by (\ref{eqCT}). It follows from Lemma \ref{PreLemma1} that $\mu_{\tau}(A) = *F_A - \tau$ defines a moment map for the action of $\mathcal{G}(P)$ on $\mathcal{A}(P)$. The \textbf{weights} associated to the gauge action with respect to this moment map are defined by
			\begin{align} \label{eqW} w_{\tau}(A,\xi) := \lim_{t\rightarrow \infty} \langle *F_{e^{\textbf{i}t\xi}A} - \tau, \xi \rangle \end{align}
for every $\xi \in W^{2,2}(\Sigma, \text{ad}(P))$ and $A \in \mathcal{A}(P)$. Differentiating the right hand side in time yields
			\begin{align} \label{MWeq0}
					\frac{d}{dt} \langle *F_{e^{\textbf{i}t\xi}A} - \tau, \xi \rangle = \langle - * d_{e^{\textbf{i}t\xi}A} * d_{e^{\textbf{i}t\xi}A} \xi , \xi \rangle = ||d_{e^{\textbf{i}t\xi}A} \xi||_{L^2}^2 \geq 0
			\end{align}		
and therefore $w_{\tau}(A,\xi) \in \mathbb{R}\cup\{+\infty\}$ is well-defined.

\begin{Remark}
The weights can be defined when $\xi$ is only of Sobolev class $W^{1,2}$. The calculation above shows
		\begin{align}
				\label{Weq1}
						w_{\tau}(A,\xi) = \langle *F_A - \tau, \xi \rangle + \int_0^{\infty} ||d_{e^{\textbf{i}t\xi}A} \xi||_{L^2}^2\,  dt 
		\end{align}
and the right hand side is well-defined for $\xi \in W^{1,2}(\Sigma,\text{ad}(P))$.
\end{Remark}

We show in Proposition \ref{PropMaxWeight} and Lemma \ref{WLemma1} that there exists a one to one correspondence between finite weights $w_{\tau}(A,\xi) < \infty$ and
				$$\left\{ (P_Q, \xi_0)\, \left| \, \begin{array}{c} \text{$\xi_0 \in \mathfrak{g}$,  $Q = Q(\xi_0)$} \\
																													\text{$P_Q$ is a principal $Q$ bundle} \\
																													\text{$P_Q \subset (P^c, J_A)$ is a holomorphic reduction}		
																				 \end{array} 
																				 \right.\right\}.$$
In particular, using a deep regularity result of Uhlenbeck and Yau \cite{UYau:1986}, we note that for every finite weight the section $\xi \in \Omega^0(\Sigma,\text{ad}(P))$ is smooth provided $A$ is a smooth connection.

Using this geometric description, we show in Proposition \ref{PropStabChar} that the algebraic stability of $(P^c, J_A)$ is equivalent to the conditions on the weights $w_{\tau}(A,\xi)$ required in the Hilbert-Mumford criterion. In the last subsection we prove the moment weight inequality
					$$- \frac{w_{\tau}(A,\xi)}{||\xi||_{L^2}} \leq \inf_{g \in \mathcal{G}^c(P)} ||\mu_{\tau}(gA)||.$$
As an immediate Corollary, $A$ is $\mu_{\tau}$-unstable whenever there exists a negative weight. Proposition \ref{PropStabChar} shows that the later is true if and only if $(P^c,J_A)$ is unstable.

\subsection{Finite weights}

It is more convenient to describe the weights in the language of vector bundles. For this, fix a faithfull representation $G \hookrightarrow U(n)$ and identify $G$ with a subgroup of $U(n)$. Let $E := P\times_G \mathbb{C}^n$ denote the associated vector bundle with structure group $G$. Then there are well-defined subbundles
				$$G(E),\, \mathfrak{g}(E),\, G^c(E),\, \mathfrak{g}^c(E) \,\subset\, \text{End}(E)$$
which consist of endomorphisms that in any trivialization are contained in $G$, $\mathfrak{g}$, $G^c$ and $\mathfrak{g}^{c}$ respectively. There are canonical identifications
				$$\mathcal{G}(P) \cong \mathcal{G}(E) = \Omega^0(\Sigma, G(E)), \qquad \text{ad}(P) \cong \mathfrak{g}(E) \subset \text{End}(E)$$
and 
				$$\mathcal{G}(P^c) \cong \mathcal{G}^c(E) = \Omega^0(\Sigma, G^c(E)), \qquad \text{ad}(P^c) \cong \mathfrak{g}^c(E) \subset \text{End}(E).$$
We denote by $\mathcal{A}_G(E)$ the space of $G$-connections on $E$ which is canonically isomorphic to $\mathcal{A}(P)$. Assume for convenience that the invariant inner product on $\mathfrak{g}$ is obtained by restriction of the standard inner product 
						$$\langle \xi ,\eta \rangle := \text{tr}(\xi \eta^*)$$ 
on $\mathfrak{u}(n)$.

\begin{Proposition} \label{PropMaxWeight}
Consider the setting described above. Let $A \in \mathcal{A}_G(E)$ be a smooth connection and let $\xi \in W^{1,2}(\Sigma, \mathfrak{g}(E))\backslash\{0\}$. If $w_{\tau}(A,\xi) < \infty$, then the following holds:		
			\begin{enumerate}
				\item The endomorphism $\textbf{i}\xi$ has constant eigenvalues $\lambda_1 < \cdots < \lambda_r$. The corresponding eigenspaces are unitary subbundles $D_j$ and decompose $E$ as orthogonal direct sum $E = D_1\oplus\cdots \oplus D_r$.
				\item Each partial sum $E_j := D_1\oplus\cdots\oplus D_j$ is a holomorphic subbundle of $E$. This yields a holomorphic filtration
															$$0 < E_1 < E_2 < \cdots < E_r = E.$$
				\item The weight of $\xi$ is given by the formula
							$$w_{\tau}(A,\xi) = 2\pi \sum_{j=1}^r \lambda_j c_1(D_j) - \langle \tau, \xi \rangle$$ 				
			\end{enumerate}			
\end{Proposition}

This is Lemma 4.2 in \cite{Mundet:2000}. Before giving the proof, we need to discuss the regularity of weakly holomorphic subbundles. 

\begin{Definition}
Let $E$ be a holomorphic hermitian vector bundle. A weakly holomorphic subbundle of $E$ is a section $\pi \in W^{1,2}(\Sigma, \text{End}(E))$ satisfying $\pi = \pi^2 = \pi^*$ and $(\mathds{1} - \pi)\bar{\partial} (\pi) = 0$.
\end{Definition}

The following theorem is a special case of a more general result of Uhlenbeck and Yau \cite{UYau:1986}. They prove that weakly holomorphic subbundles of holomorphic hermitian vector bundles over arbitrary Kähler manifolds correspond to torsion-free coherent subsheaves. Since any torsion-free coherent sheaf over a Riemann surface is locally free, this reduces to the following:

\begin{Theorem}[Uhlenbeck and Yau \cite{UYau:1986}] \label{ThmUY}
If $\pi \in W^{1,2}(\Sigma, \text{End}(E))$ is a weakly holomorphic subbundle, then $\pi$ is the projection on a smooth holomorphic subbundle $E' \subset E$.
\end{Theorem}

\begin{proof}[Proof of Proposition \ref{PropMaxWeight}]
Let $0\neq \xi \in W^{1,2}(\Sigma, \mathfrak{g}(E))$ be given and assume $w_{\tau}(A,\xi) < \infty$. Since $\mathfrak{g}^c = \mathfrak{g} \oplus \textbf{i}\mathfrak{g}$ is per definitionem an orthogonal decomposition we have
		\begin{align*}
				||d_{e^{\textbf{i}t\xi}A} \xi||^2 &= \frac{1}{2}||\bar{\partial}_{e^{\textbf{i}t\xi}A} \xi||^2 
									= \frac{1}{2}||\text{Ad}\left(e^{\textbf{i}t\xi}\right) \circ \bar{\partial}_{A} \circ \text{Ad}\left(e^{-\textbf{i}t\xi}\right) (\xi)||^2 \\
									&= \frac{1}{2}||e^{\textbf{i}t\xi} \bar{\partial}_{A}( \xi) e^{-\textbf{i}t\xi}||^2.
		\end{align*}
and from (\ref{Weq1}) follows
		\begin{align}\label{Weq2} w_{\tau}(A,\xi) = \int_{\Sigma} \langle *F_A - \tau, \xi \rangle\, dvol_{\Sigma} + \frac{1}{2} \int_0^{\infty} ||e^{\textbf{i}t\xi} \bar{\partial}_A(\xi) e^{-\textbf{i}t\xi}||^2\, dt. \end{align}

Denote $A_t := e^{\textbf{i}t\xi}(A)$ and let $k \geq 1$ be an integer. Then follows
			$$\bar{\partial} \text{tr}(\xi^k) = \text{tr}(\bar{\partial}_{A_t} (\xi^k)) = k \text{tr}(\xi^{k-1} \bar{\partial}_{A_t} (\xi))$$
and the Cauchy-Schwarz inequality $|\text{tr}(AB)| \leq ||A||\cdot ||B||$ yields
			\begin{align*}
						\int_{\Sigma} |\bar{\partial}\, \text{tr}(\xi^k)| \, dvol_{\Sigma}
						&\leq  k\int_{\Sigma} ||\xi^{k-1}|| \cdot || \bar{\partial}_{A_t} \xi ||  \, dvol_{\Sigma}\\
						&=  k ||\xi^{k-1}||_{L^2} \cdot || e^{\textbf{i}\xi t}  \bar{\partial}_A (\xi) e^{-\textbf{i}\xi t} ||_{L^2}.
			\end{align*}
Since $w_{\tau}(A,\xi)$ is finite, it follows from (\ref{Weq1}) that there exists a sequence $t_j \rightarrow \infty$ such that
		\begin{align} \label{Weq3} \lim_{j\rightarrow \infty} || e^{\textbf{i}\xi t_j}  \bar{\partial}_A (\xi)  e^{-\textbf{i}\xi t_j} ||_{L^2} = 0. \end{align}
Hence $\bar{\partial} \text{tr}(\xi^k) = 0$ and it follows from the maximum principle that $\text{tr}(\xi^k)$ is constant. Denote the eigenvalues of $\textbf{i}\xi$ with repetition according to their multiplicity by $\lambda_1' \leq \cdots \leq \lambda_n'$. Then
			$$\text{tr}(\xi^k) = (\lambda'_1)^k + \ldots + (\lambda'_n)^k$$
is constant for every $k \geq 1$. This is only possible if all the functions $\lambda_j'$ are constant and hence $\textbf{i}\xi$ has constant eigenvalues.

Let $\lambda_1 < \cdots < \lambda_r$ be the distinct eigenvalues of $\textbf{i}\xi$. Since $\textbf{i}\xi$ is a normal (hermitian) operator, it follows from basic linear algebra that the eigenspaces are orthogonal. Moreover, if $\Gamma_j$ is a small loop around the eigenvalue $\lambda_j$ in the complex plane, we can explicitly write down the orthogonal projection $\pi_j': E \rightarrow D_j$ onto the eigenspace of $\lambda_j$ as
			$$\pi_j' := \frac{1}{2\pi\textbf{i}} \int_{\Gamma_j} (z\mathds{1} - \textbf{i}\xi)^{-1} \,dz.$$
These projections have regularity $\pi_j' \in W^{1,2}(\text{End}(E))$ and clearly satisfy $\pi_j' = (\pi_j')^2 = (\pi_j')^*$.

We show next that the projections $\pi_j := \pi_1' + \cdots + \pi_j': E \rightarrow E_j$ define weakly holomorphic subbundles. By construction
					\begin{align} \label{Weq4} \textbf{i}\xi = m_1 \pi_1 + \cdots + m_r \pi_r \end{align}
for some $m_1, \ldots, m_r \in \mathbb{R}$. Write $\bar{\partial}_A(\xi) = [\hat{\xi}_{ij}]$ with respect to the splitting $E = D_1\oplus \cdots \oplus D_r$. Then holds
				$$\left[ e^{\textbf{i}t\xi} \hat{\xi} e^{-\textbf{i}t\xi} \right]_{ij} = e^{(\lambda_j - \lambda_i)t} \hat{\xi}_{ij}$$
and (\ref{Weq3}) implies $\hat{\xi}_{ij} = 0$ for $i > j$. Thus $\bar{\partial}_A (\textbf{i}\xi)$ is upper triangular and (\ref{Weq4}) yields
				\begin{align} \label{Weq5} 0 = (\mathds{1} - \pi_j) (\bar{\partial}_A \xi)  \pi_j = \sum_{k=1}^r m_k (\mathds{1} - \pi_j) \bar{\partial}_A (\pi_k) \pi_j. \end{align}
The Leibniz rule provides the formula
	\begin{align*}	
		(\mathds{1} - \pi_j) \bar{\partial}_A (\pi_k) \pi_j = 
					\begin{cases} (\mathds{1}-\pi_k)(\mathds{1} - \pi_j) \bar{\partial}_A (\pi_j) & \text{for $k > j$} \\
																					(\mathds{1} - \pi_j) \bar{\partial}_A (\pi_j) & \text{for $k = j$} \\
																					(\mathds{1} - \pi_j) (\bar{\partial}_A(\pi_k) - \pi_k \bar{\partial}_A(\pi_j)) & \text{for $k<j$}
					\end{cases}.
	\end{align*}				
This implies together with (\ref{Weq5}) the formula $(\mathds{1} - \pi_j)\bar{\partial}_A(\pi_j) = 0$ by induction on $j$. Hence $\pi_j$ defines a weakly holomorphic subbundle and $E_j$ is smooth by Theorem \ref{ThmUY}. This proves the first two parts of the theorem.

Write $\bar{\partial}_A$ with respect to the splitting $E = D_1 \oplus \cdots \oplus D_r$ as
	\begin{align*}
			\bar{\partial}_A = 
				\left(\begin{array}{cccc}
									\bar{\partial}_{A_1}  & 			A_{1 2} 					& \ldots & A_{1 r} \\
																0       & \bar{\partial}_{A_2}		& \ldots & A_{2 r} \\
							 					 \vdots	  		  &      \vdots							&	\ddots & \vdots  \\
																0			  &      		0								& \ldots & \bar{\partial}_{A_r}
							\end{array}\right)
		\end{align*}
where $A_{ij} \in \Omega^{0,1}(D_i \otimes D_j^*)$ and $\bar{\partial}_{A_j}$ is the Cauchy-Riemann operator corresponding to the induced unitary connection $A_j \in \mathcal{A}(D_j) \cong \mathcal{A}(E_j/E_{j-1})$. Decompose $\bar{\partial}_A = \bar{\partial}_{A_+} + A_0$ with
		\begin{align*}
			\bar{\partial}_{A_+} = 
										\left(\begin{array}{cccc}
												\bar{\partial}_{A_1}  & 			0 								& \ldots & 0 \\
																			0       & \bar{\partial}_{A_2}		& \ldots & 0 \\
															 \vdots	  		  &      \vdots							&	\ddots & \vdots  \\
																			0			  &      0									& \ldots & \bar{\partial}_{A_r}
										\end{array}\right),
				\qquad					
				A_0= 
										\left(\begin{array}{cccc}
																				0					  & 			A_{1 2} 					& \ldots & A_{1 r} \\
																			  0           & 				0								& \ldots & A_{2 r} \\
																		 \vdots	  		  &      \vdots							&	\ddots & \vdots  \\
																			  0  	   		  &         0    						& \ldots & 		0
									\end{array}\right)
		\end{align*}
We claim that $e^{\textbf{i}t\xi}(A)$ converges uniformly to $A_+ := A_1 \oplus \cdots \oplus A_r$ as $t \rightarrow \infty$. In fact
				$$\bar{\partial}_{A_t} - \bar{\partial}_{A_+} = e^{\textbf{i}t\xi} A_0 e^{-\textbf{i}t\xi}$$
and
				$$[e^{\textbf{i}t\xi} A_0 e^{-\textbf{i}t\xi}]_{ij} = -\textbf{i} e^{t(\lambda_i - \lambda_j)} (\lambda_j - \lambda_i) A_{ij} $$				
decays exponentially to zero, since $A_0$ is strictly upper triangular. This in turn implies that $e^{\textbf{i}t\xi}A$ converges to $A_+$ and hence
			\begin{align*}
						w_{\tau}(A,\xi) &= \lim_{t\rightarrow \infty} \langle *F_{e^{\textbf{i}t\xi} A}, \xi \rangle = \langle *F_{A_+} - \tau, \xi \rangle 
																= \sum_{j=1}^r \langle*F_{A_j}, \xi\rangle - \langle \tau, \xi \rangle \\
														&= \sum_{j=1}^r \textbf{i}\lambda_j \int_{\Sigma} \text{tr}(F_{A_j})\, dvol_{\Sigma} - \langle \tau, \xi \rangle \\
														&= 2\pi \sum_{j=1}^r \lambda_j c_1(D_j) - \langle \tau, \xi \rangle.
			\end{align*}

\end{proof}

\begin{Cor}\label{CorSplit}
Suppose $\xi \in \Omega^0(\Sigma, \mathfrak{g}(E))$ yields a finite weight $w_{\tau}(A,\xi)$. Then the limit
		$$A_+ := \lim_{t\rightarrow \infty} e^{\textbf{i}t\xi} A$$
exists in $\mathcal{A}_G(E)$. Moreover, the splitting 
					$$E =  D_1\oplus \cdots \oplus D_r$$ 
is holomorphic with respect to $A_+$ and on each factor the holomorphic structure agrees with the one induced by the isomorphism $D_j \cong E_j/E_{j-1}$
\end{Cor}

\begin{proof} This follows directly from the proof of Proposition \ref{PropMaxWeight}.\end{proof}

\begin{Remark} The Corollary shows that $A_+$ remains in the complexified orbit $\mathcal{G}^c(A)$ if and only if the holomorphic filtration determined by $\xi$ splits holomorphically.
\end{Remark}

We reformulate the characterization of the finite weights in intrinsic terms. Let $A \in \mathcal{A}(P) \cong \mathcal{A}_G(E)$ and suppose that $\xi$ is a smooth section of $\text{ad}(P) \cong \mathfrak{g}(E)$ which yields a finite weight $w_{\tau}(A,\xi)$. By Proposition \ref{PropMaxWeight} this defines a filtration
				$$0 < E_1 < E_2 < \cdots < E_r = E$$
and there exist unitary trivializations of this filtration such that $\xi =  \xi_0$ where $\xi_0 = - \textbf{i}\text{diag}(\lambda_1, \ldots, \lambda_r)$ is a block diagonal matrix with $\lambda_1 < \lambda_2 < \cdots < \lambda_r$. This trivialization yields a reduction $P_{K(\xi)} \subset P$ to $K(\xi_0) := C_G(\xi_0)$. Note that $\xi_0$ gives rise to a constant central section of $\text{ad}(P_{K(\xi)}) \subset \text{ad}(P)$ and agrees with $\xi$ in $\text{ad}(P)$. We can rewrite the formula for the weight as
					$$w_{\tau}(A,\xi) := \int_{\Sigma} \langle *F_{A_+}, \xi \rangle \, dvol_{\Sigma} - \langle\tau, \xi \rangle$$ 
where $A_+ \in \mathcal{A}(P_{K(\xi)})$ is a $K(\xi)$-connection. It follows from Chern-Weyl theory that the right hand side does not change when we replace $A_+$ by another $K(\xi)$-connection. The weight depends therefore only on the reduction $P_{K(\xi)} \subset P$ and $\xi$. The complexification yields a reduction $P_{K(\xi)}^c = P_{L(\xi)} \subset P^c$ to the Levi subgroup $L(\xi_0) \subset G^c$ (see Definition \ref{DefnToralGen}). The reduction $P_{L(\xi)} \subset P^c$ is holomorphic if and only if $\bar{\partial}_A$ takes values in $\mathfrak{l}(\xi_0)$ and this is the case if and only if the filtration determined by $\xi$ splits holomorphically. In contrast, the extension $L(\xi_0) \subset Q(\xi_0)$ yields a reduction $P_{Q(\xi_0)} \subset P^c$ to the stabilizer of the filtration determined by $\xi_0$ within $G^c$. This reduction is always holomorphic, since $\bar{\partial}_A$ is upper block triangular.

Conversely, let $P_Q \subset P^c$ be a holomorphic reduction to a parabolic subgroup $Q = Q(\xi_0) \subset G^c$. This yields a canonical reduction $P_K \subset P$ to $K = C_G(\xi_0)$, since $G^c/Q(\xi_0) \cong G/C_G(\xi_0)$. Since $\xi_0$ is contained in the center of $K$, it gives rise to a constant section in $\text{ad}(P_K)$ and its image under the embedding $\text{ad}(P_K) \subset \text{ad}(P)$ yields a section $\xi \in \Omega^0(\Sigma, \text{ad}(P)$ which gives rise to a finite weight $w_{\tau}(A,\xi)$. We summarize our discussion in the following Lemma.

\begin{Lemma} \label{WLemma1}
Let $P \rightarrow \Sigma$ be a principal $G$ bundle, let $A \in \mathcal{A}(P)$ be a smooth connection and let $P^c := P\times_G G^c$ denote the complexification of $P$ endowed with the holomorphic structure determined by $A$. There exists a one-to-one correspondence between
		$$\{ \xi \in \Omega^0(\Sigma, \text{ad}(P))\, |\, w_{\tau}(A,\xi) < \infty \}$$
and
			$$\left\{ (P_Q, \xi_0)\, \left| \, \begin{array}{c} \text{$\xi_0 \in \mathfrak{g}$,  $Q = Q(\xi_0)$} \\
																													\text{$P_Q$ is a principal $Q$ bundle} \\
																													\text{$P_Q \subset P^c$ is a holomorphic reduction}		
																				 \end{array} 
																				 \right.\right\}$$
More precisely, every reduction $P_Q \subset P^c$ yields a canonical reduction $P_K \subset P$ to $K = C_G(\xi_0)$. The toral generator $\xi_0$ yields a constant section of $\text{ad}(P_K)$ and its image in $\text{ad}(P)$ yields $\xi$. Moreover, the weight is given by the formula
						$$w_{\tau}(A,\xi) = \int_{\Sigma} \langle *F_B - \tau, \xi \rangle \, dvol_{\Sigma}$$
for any connection $B \in \mathcal{A}(P_K)$
\end{Lemma}

\begin{proof} This follows directly from the preceding discussion \end{proof}

The next Lemma describes how the weights behave under an extension $G \hookrightarrow H$ of the structure group.

\begin{Lemma} \label{LemWeightExt}
Let $H$ be a compact connected Lie group and fix an invariant inner product on its Lie algebra $\mathfrak{h}$. Suppose that there exists a monomorphism $G \hookrightarrow H$ which identifies $G$ with a subgroup of $H$ and assume that the invariant inner product on $\mathfrak{g}$ is obtained by restriction of the one on $\mathfrak{h}$. Let $P \rightarrow \Sigma$ be a pricipal $G$ bundle of central type $\tau \in Z(\mathfrak{g})$ defined by (\ref{eqCT}) and denote by $P_H := P \times_G H$ the associated $H$ bundle.

	\begin{enumerate}
		\item The central type $\tau_H \in Z(\mathfrak{h})$ of $P_H$ is the image of $\tau$ under the orthogonal projection
							$$Z(\mathfrak{g}) \hookrightarrow \mathfrak{h} \cong Z(\mathfrak{h}) \oplus [\mathfrak{h}, \mathfrak{h}] \rightarrow Z(\mathfrak{h}).$$
	
		\item Let $\xi \in \Omega^0(\Sigma, \text{ad}(P))$ and denote by $\xi_H \in \Omega^0(\Sigma, \text{ad}(P_H))$ the image of $\xi$ under the embedding $\text{ad}(P) \subset \text{ad}(P_H)$. Then
							$$w_{\tau}(A,\xi) = w_{\tau_H}(A,\xi_H) + \int_{\Sigma}\langle \tau_H - \tau, \xi_0 \rangle \, dvol_{\Sigma}.$$
							
		\item Let $\xi_H \in \Omega^0(\Sigma, \text{ad}(P_H))$ be a section with $w_{\tau_H}(A,\xi) < \infty$ and denote by $\xi \in \Omega^0(\Sigma, \text{ad}(P))$ the image of $\xi_H$ under the orthogonal projection $\text{ad}(P_H) \rightarrow \text{ad}(P)$. Then
						$$w_{\tau}(A,\xi) = w_{\tau_H}(A,\xi_H) + \int_{\Sigma}\langle \tau_H - \tau, \xi \rangle \, dvol_{\Sigma}.$$
	\end{enumerate}

\end{Lemma}

\begin{proof}
For the first part, note that $\mathfrak{h} = Z(\mathfrak{h})\oplus [\mathfrak{h}, \mathfrak{h}]$ yields an orthogonal decomposition with respect to any invariant inner product of $\mathfrak{h}$. The orthogonal projection of $\tau$ onto $Z(\mathfrak{h})$ does therefore depend only on the embedding of $G$ into $H$ and it is easy to verify that it satisfies (\ref{eqCT}) for $P_H$.

By Lemma \ref{WLemma1} there exists $\xi_0 \in \mathfrak{g}$ and a reduction $P_K \subset P$ to a principal $K = C_G(\xi_0)$ bundle such that $\xi$ is the image of the constant section $\xi_0$ under the embedding $\text{ad}(P_K) \subset \text{ad}(P)$. Moreover, 
					$$w_{\tau}(A,\xi) = \int_{\Sigma} \langle *F_B - \tau, \xi \rangle \, dvol_{\Sigma}$$
for any connection $B \in \mathcal{A}(P_K)$. Define $\tilde{K} = C_H(\xi_0)$ and $P_{\tilde{K}} := P_K \times_K \tilde{K} \subset P_H$. Then $\xi_H$ agrees with the image of $\xi_0$ under the embedding $\text{ad}(P_{\tilde{K}}) \subset \text{ad}(P_H)$ and Lemma \ref{WLemma1} yields
					$$w_{\tau_H}(A,\xi) = \int_{\Sigma} \langle *F_B - \tau_H, \xi_0 \rangle $$
for any connection $B \in \mathcal{A}(P_{\tilde{K}})$. In particular, for $ B \in \mathcal{A}(P_K) \subset \mathcal{A}(P_{\tilde{K}})$, we get
					$$w_{\tau}(A,\xi) - w_{\tau_H}(A,\xi) = \int_{\Sigma} \langle \tau_H - \tau, \xi_0 \rangle \, dvol_{\Sigma}$$
and this proves the second part.
			
The third part follows by a similar argument. Note that the proof of Proposition \ref{PropMaxWeight} implies that there exists a connection $B = A_+ \in \mathcal{A}(P)\cap \mathcal{A}(P_{\tilde{K}})$ for the reduction $P_{\tilde{K}} \subset P_H$ associated to $\xi_H$. For such a connection holds $\langle \xi_H, F_B \rangle = \langle \xi, F_B \rangle$ and the claim follows as in the second part.
\end{proof}

\subsection{Weights and algebraic stability}

The following Proposition characterizes the (algebraic) stability of the holomorphic principal bundle $(P^c, J_A)$ in terms of the associated weights $w_{\tau}(A,\xi)$. 

\begin{Proposition}[\textbf{Characterization of Stability}] \label{PropStabChar}
Let $P$ be a principal $G$ bundle of central type $\tau \in Z(\mathfrak{g})$ defined by (\ref{eqCT}). Let $A \in \mathcal{A}(P)$ be a smooth connection and let $P^c := P \times_G G^c$ be the complexified principal bundle endowed with the induced holomorphic structure $J_A$.
		\begin{enumerate}
			\item $(P^c, J_A)$ is stable if and only if $w_{\tau}(A,\xi) > 0$ for all $\xi \in W^{1,2}(\Sigma, \text{ad}(P))$ which are not constant central sections.
			
			\item $(P^c, J_A)$ is polystable if and only if $w_{\tau}(A,\xi) \geq 0$ for all $\xi \in W^{1,2}(\Sigma, \text{ad}(P))$ and whenever $w_{\tau}(A,\xi) = 0$ the associated (smooth) reduction $P_{L(\xi)} \subset P_{Q(\xi)} \subset P^c$ is holomorphic.
			
			\item $(P^c, J_A)$ is semistable if and only if $w_{\tau}(A,\xi) \geq 0$ for all $\xi \in W^{1,2}(\Sigma, \text{ad}(P))$.
		
			\item $(P^c, J_A)$ is unstable if and only if there exists $\xi \in W^{1,2}(\Sigma, \text{ad}(P))$ with $w_{\tau}(A,\xi) < 0$.
		\end{enumerate}
\end{Proposition}

\begin{proof}
Using the geometric interpretation of the finite weights in Lemma \ref{WLemma1} we can reduces the proof to a lemma of Ramanathan \cite{Ramanathan:1975}. The proof will be given on page \pageref{proofHM} below.
\end{proof}

\subsubsection*{Reduction argument}
We reduce the theorem to the case where $Z(G)$ is discrete and $\tau = 0$. Recall that the invariant inner product on $\mathfrak{g}$ yields the decomposition $\mathfrak{g} = Z(\mathfrak{g})\oplus [\mathfrak{g}, \mathfrak{g}]$ of the Lie algebra into its center and a semisimple subalgebra. The center yields a trivial $Z(\mathfrak{g})$ subbundle $V \subset \text{ad}(P)$ and its orthogonal complement can be identified with $\text{ad}(P/Z_0(G))$.

\begin{Lemma} \label{WeightsPblLemma}
Assume the setting of Proposition \ref{PropStabChar}. Let $\xi \in \Omega^0(\Sigma, \text{ad}(P))$ with $w_{\tau}(A,\xi) < \infty$ and decompose $\xi = \xi^z + \xi^{ss}$ with respect to the splitting $\text{ad}(P) = V \oplus \text{ad}(P/Z_0(G))$. Then
			$$w_{\tau}(A,\xi) = w_0(\bar{A},\xi^{ss})$$
where $\bar{A}\in \mathcal{A}(P/Z_0(G))$ denotes the induced connection on $P/Z_0(G)$.
\end{Lemma}

\begin{proof}
By Lemma \ref{WLemma1} exists a reduction $P_K \subset P$ and an element $\xi_0 \in \mathfrak{g}$ which gives rise to a constant central section in $\text{ad}(P_K)$ and such that $\xi$ is the image of $\xi_0$ under the embedding $\text{ad}(P_K) \subset \text{ad}(P)$. Decompose $\xi_0 = \xi_0^z + \xi_0^{ss}$ with respect to $\mathfrak{g} = Z(\mathfrak{g})\oplus [\mathfrak{g},\mathfrak{g}]$. Then $\xi_0^z$ yields $\xi^z$ and $\xi_0^{ss}$ yields $\xi^{ss}$ under the embedding $\text{ad}(P_K) \subset \text{ad}(P)$.
By Lemma \ref{WLemma1} the weight is given by
		$$w_{\tau}(A,\xi) = \int_{\Sigma} \langle *F_B - \tau, \xi^{ss}\rangle  \,dvol_{\Sigma} + \int_{\Sigma} \langle *F_B - \tau, \xi^z \rangle\, dvol_{\Sigma}.$$
for any connection $B \in \mathcal{A}(P_K)$. Since $\tau \in Z(\mathfrak{g})$ is orthogonal to $[\mathfrak{g},\mathfrak{g}]$, we have
		$$\int_{\Sigma} \langle *F_B - \tau, \xi^{ss} \rangle \, dvol_{\Sigma} = \langle *F_B, \xi^{ss} \rangle \, dvol_{\Sigma}= w_0(\bar{A}, \xi^{ss}).$$
The definition of the central type (\ref{eqCT}) shows $\int_{\Sigma} \langle *F_B - \tau, \xi^z \rangle \, dvol_{\Sigma} = 0$ and this completes the proof of the Lemma.
\end{proof}

\subsubsection*{The main argument}

The following result is a reformulation of Lemma 2.1 in \cite{Ramanathan:1975} by Ramanathan. 

\begin{Lemma}\label{WeightPblLemma2}
Assume the setting of Proposition \ref{PropStabChar} and suppose in addition that $Z_0(G)$ is discrete and $\tau = 0$. $P^c$ is stable (resp. semistable) with respect to Definition \ref{DefnStab2} if and only if $w_0(A,\xi) > 0$ (resp. $w_0(A,\xi) \geq 0$) for all $\xi \in W^{1,2}(\Sigma, \text{ad}(P))$.
\end{Lemma}

\begin{proof}
Let $\xi \in \Omega^0(\Sigma, \text{ad}(P))$ with $w_0(A,\xi) < \infty$ be given. By Lemma \ref{WLemma1} exists a reduction $P_K \subset P$ and an element $\xi_0 \in \mathfrak{g}$ such that $K = C_G(\xi_0)$ and $\xi$ is the image of $\xi_0$ under the embedding $\text{ad}(P_K) \subset \text{ad}(P)$.

Let $T \subset G$ be a maximal torus whose Lie algebra contains $\xi_0$ and let $R_0^+ = \{\alpha_1, \ldots, \alpha_r\}$ be a system of simple roots with respect to $T$ whose Weyl-chamber contains $\xi_0$. Recall that $\alpha_j = \textbf{i}a_j$ with $a_j \in \text{Hom}(\mathfrak{t},\mathbb{R})$ and define $t_j \in \mathfrak{t}$ by $a_j = \langle t_j, \cdot \rangle$. The elements $\check{t}_1,\ldots,\check{t}_r \in \mathfrak{t}$ defined by (\ref{rootseq2}) yield a basis of $\mathfrak{t}$ and $\xi_0$ has the shape
					$$\xi_0 = \sum_{j=1}^r x_j \check{t}_j$$
with $x_j \geq 0$. Note that $\check{t}_j$ lies in the center of the Lie algebra of $K = C_G(\xi_0)$ when $x_j > 0$. Then $\check{t}_j$ gives rise to a constant central section of $\text{ad}(P_K)$ and
				\begin{align} \label{HMeq1} w_0(A,\xi) = \sum_{j=1}^r x_j w_0(A,\check{t}_j). \end{align}

Fix $1 \leq j \leq r$ with $x_j > 0$ and denote $Q_j := Q(\check{t}_j)$. This is a maximal parabolic subgroup of $G^c$ which contains $Q(\xi_0)$ and the extension $P_{Q_j} := P_{Q(\xi)}\times_{Q(\xi)} Q_j  \subset P^c$ yields a maximal parabolic reduction. Let $\chi: Q_j \rightarrow \mathbb{C}^*$ be the determinant of the action of $Q_j$ on its Lie algebra and denote by $\dot{\chi}: \mathfrak{q}_j \rightarrow \mathbb{C}$ the induced map on the Lie algebra. Chern-Weyl theory yields the relation
					$$c_1(\text{ad}(P_{Q_j})) = \frac{\textbf{i}}{2\pi} \int_{\Sigma} \dot{\chi}(F_{B}) $$
for a connection $B \in \mathcal{A}(P_K)$. For $\eta \in \mathfrak{q}_j$ the value of $\dot{\chi}(\eta)$ is given as the trace of $\text{ad}(\eta) := [\eta,\cdot]$ acting on
					\begin{align} \label{HMeq0} \mathfrak{q}_j = \mathfrak{t} \oplus \bigoplus_{\alpha \in R(\check{t}_j)} \mathfrak{g}_{\alpha}. \end{align}
where $R(\check{t}_j)$ is defined by (\ref{Req1}). This decomposition is unitary and by definition of the roots we have $\text{ad}(t) e_{\alpha} = \alpha(t) e_{\alpha}$ for $t \in \mathfrak{t}$. This shows
				\begin{align} \dot{\chi}(\eta) = \sum_{\alpha \in R(\check{t}_j)} \alpha(\eta) \label{pblweightseq1} \end{align}
for all $\eta \in \mathfrak{t}$. Since $\dot{\chi}$ vanishes on $[\mathfrak{q}_j, \mathfrak{q}_j]$ it vanishes on all root space $\mathfrak{g}_{\alpha}$ with $\{\alpha, -\alpha\} \subset R(\check{t}_j)$. These are the roots in $\tilde{R}(\check{t}_j)$ which produce the Levi subgroup $L(\check{t}_j)$. The remaining root spaces $\mathfrak{g}_{\alpha}$ with $\alpha \in R(\check{t}_j) \backslash \tilde{R}(\check{t}_j)$ form a nilpotent subalgebra. This shows that (\ref{pblweightseq1}) remains valid for all $\eta \in \mathfrak{q}_j$ if one extends the roots by complex linearity over $\mathfrak{t}^c$ and by zero over the root spaces.

Denote by $R^+$ the positive roots and by $R^-(\check{t}_{j}) = R(\check{t}_{j}) \backslash R^+$ the negative roots whose root spaces are contained in $\mathfrak{q}_j$. Then $\dot{\chi} = \gamma_1 + \gamma_2$ with
			$$\gamma_1 = \sum_{\alpha \in R^+} \alpha, \qquad \gamma_2 := \sum_{\alpha \in R^-(\check{t}_{j}) } \alpha$$
and
			$$\langle \alpha_i, \gamma_1 \rangle = \sum_{\alpha \in R^+} \langle t_i, t_{\alpha} \rangle = |t_i|^2 + \sum_{\alpha \in R^+ \backslash\{\alpha_j\}} \langle t_i, t_{\alpha} \rangle$$
holds for every simple root $\alpha_i$.	The root reflection 
			$$s_j : \mathfrak{t} \rightarrow \mathfrak{t}, \qquad s_j(t) := t - \frac{2 \langle t, t_j\rangle}{|t_j|^2} t_j$$
restricts to a permutation of $R^+ \backslash\{\alpha_j\}$. Indeed, any root has a unique representation $t_{\alpha} = \sum_{k=1}^r c_k t_k$ and all coefficients happen to have the same sign. Applying the reflection $s_j$ changes only the coefficient $c_j$ and thus $s_j(\alpha)$ remains positive if $c_k > 0$ for some coefficient $k\neq j$. Using this symmetry we conclude
					\begin{align} \label{HMeq2}\langle \alpha_i, \gamma_1 \rangle= |t_i|^2. \end{align}
A similar argument shows for $i \neq j$
					\begin{align} \label{HMeq3} \langle \alpha_i, \gamma_2 \rangle = \sum_{\alpha \in R^-(\check{t}_{j})} \langle t_i, t_{\alpha} \rangle = -|t_i|^2 + \sum_{\alpha \in R^-(\check{t}_{j}) \backslash\{\alpha_j\}} \langle t_i, t_{\alpha} \rangle = -|t_i|^2. \end{align}
This shows $\dot{\chi}(t_i) = 0$ for $i\neq j$. As a general property of root systems (see \cite{LiegpKnapp} Lemma 2.51) it holds $\langle t_j, t_i \rangle \leq 0$ for distinct simple roots $\alpha_i, \alpha_j$ and thus
					\begin{align} \label{HMeq4} \dot{\chi}(t_j) = |t_j|^2 + \sum_{\alpha \in R^-(\check{t}_{j})} \langle t_j, t_{\alpha} \rangle  > 0 \end{align}
Combining (\ref{HMeq2}), (\ref{HMeq3}) and (\ref{HMeq4}) we conclude
			$$\dot{\chi}(t) = \textbf{i}m\langle\check{t}_j, t \rangle$$
for some $m > 0$. Hence 
			\begin{align} \label{HMeq5}
					c_1(\text{ad}(P_{Q_j})) = \frac{-m}{2\pi} \int_{\Sigma} \langle F_B, \check{t}_j \rangle = \frac{-m}{2\pi} w_0(A, \check{t}_j).
			\end{align}

Suppose now that $P^c$ is stable (resp. semistable). Then the left hand side in (\ref{HMeq5}) is negative (resp. nonpositive) and (\ref{HMeq1}) implies $w_0(A,\xi) > 0$
(resp. $w_0(A,\xi) \geq 0$). Conversely, Lemma \ref{WLemma1} show that every holomorphic reduction $P_Q \subset P^c$ to a proper maximal parabolic subgroup $Q(\xi_0) \subset Q$ is induced by some $\xi \in \Omega^0(\Sigma, \text{ad}(P))$ with $w_0(A,\xi) < \infty$. Lemma \ref{PreLemma3} shows that in (\ref{HMeq0}) exactly one coefficient $x_j$ does not vanishes. Hence (\ref{HMeq5}) implies that $c_1(\text{ad}(P_Q))$ is negative or vanishes if and only if $w_0(A,\xi)$ is positive or vanishes respectively. This establishes the converse direction and completes the proof of the Lemma.
\end{proof}

\subsubsection*{Completion of the proof}

\begin{proof}[Proof of Proposition \ref{PropStabChar}.] \label{proofHM}
We may assume by Lemma \ref{LemmaSSred1} and Lemma \ref{WeightsPblLemma} that $Z_0(G)$ is discrete and $\tau = 0$. The stable and semistable case follow then from Lemma \ref{WeightPblLemma2} and the unstable case is equivalent to the semistable case.

Assume that $P^c$ is polystable. Then there exists a holomorphic reduction $P_L \subset P^c$ to a Levi subgroup of $L \subset G^c$ and $P_L$ is a stable $L$ bundle. Let $\xi \in \Omega^0(\Sigma, \text{ad}(P))$ with $w_0(A,\xi) = 0$ be given. By Lemma \ref{WLemma1} exists $\xi_0 \in \mathfrak{g}$ and a reduction $P_K \subset P$ to an principal $K = C_G(\xi_0)$ bundle such that $\xi$ agrees with the image of $\xi_0$ under the embedding $\text{ad}(P_K) \subset \text{ad}(P)$. Using the notation from the proof of Lemma \ref{WeightPblLemma2} above, write $\xi_0$ with respect to a system of simple roots as
					$$\xi_0 = \sum_{j=1}^r x_j \check{t}_j$$
with $x_j \geq 0$. Since $P^c$ is in particular semistable, the proof of Lemma \ref{WeightPblLemma2} shows that $w_0(A,\xi) = 0$ if and only if
			$$x_j > 0 \quad \Rightarrow \quad c_1(\text{ad}(P_{Q_j})) = 0$$
where $Q_j := Q(\check{t}_j)$. We may assume (after conjugation) that $L = L(\eta_0)$ for some $\eta_0 \in \mathfrak{g}$ and $\eta_0$ is contained in the Weyl-chamber determined by our choice of simple roots. If $L$ is not contained in $Q_j$, then $Q_j' := L\cap Q_j$ is a maximal parabolic subgroup of $L$ and we have a natural holomorphic reduction $P_{Q_j'} \subset P_L$. Since $L$ and $G^c$ are reductive, the Lie algebra bundles $\text{ad}(P_L)$ and $\text{ad}(P^c)$ carry a non degenerated symmetric $\mathbb{C}$-bilinear form. Hence they are both self-dual and have vanishing first Chern-class. This shows
			\begin{align*} c_1(\text{ad}(P_{Q_j} )) &= - c_1(\text{ad}(P^c)/\text{ad}(Q_j)) = - c_1(\text{ad}(P_L)/\text{ad}(Q_j')) \\
																							&= c_1(\text{ad}(P_{Q_j'})) < 0 \end{align*}
where the last step follows from the stability of $P_L$. We have thus proven that $L \subset Q_j$ whenever $x_j > 0$ and this yields $L \subset L(\xi_0)$. Since the reduction to $L$ is holomorphic, so is the reduction to $L(\xi_0)$.

Assume conversely, that all weights are nonnegative and if $\xi \in \Omega^0(\Sigma,P^c)$ is a section with $w_0(A,\xi) = 0$ then $P_{L(\xi)} \subset P^c$ is a holomorphic reduction (where $P_{L(\xi)} = P_{K(\xi)}^c$ and $P_{K(\xi)}$ is determined by Lemma \ref{WLemma1}). It follows from Lemma \ref{WeightPblLemma2} that $P^c$ is semistable. If $P^c$ is in fact stable, then we are done. Otherwise there exists a vanishing weight $w_0(A,\xi) = 0$ and by assumption this yields a holomorphic reduction $P_{L(\xi)} \subset P^c$. In particular $A$ restricts to a connection on $P_{K(\xi)} \subset P$ and $P_{K(\xi)}$ is again of central type $0$. For the later claim let $\eta \in \mathfrak{g}$ be contained in the center of the Lie algebra of $K$ and consider its image $\eta'$ under the embedding $\text{ad}(P_K) \subset \text{ad}(P)$. Then follows
					$$\int_{\Sigma} \langle *F_B, \eta \rangle \, dvol_{\Sigma} =  w_0(A,\eta') \geq 0.$$
for any connection $B \in \mathcal{A}(P_K)$. Replacing $\eta$ by $- \eta$ shows that this expression must vanish and hence $P_{K(\xi)}$ is of central type $0$. Now Lemma \ref{WeightPblLemma2} shows that $P_{L(\xi)}$ is again semistable. If $P_{L(\xi)}$ is not stable, then there exists $\tilde{\xi} \in \Omega^0(\Sigma, \text{ad}(P_{K(\xi)})$ with $w_0(A,\tilde{\xi}) = 0$. We can consider $\tilde{\xi}$ as section $\xi'$ of $\text{ad}(P)$ which then satisfies $w_0(A,\xi') = 0$ and thus yields a strictly smaller holomorphic reduction $P_{L(\xi')} \subset P_{L(\xi)}$. If we replace $\xi$ by $\xi'$ and rerun the argument from above we obtain after finitely many iterations a section $\xi$ which satisfies $w_0(A,\xi) = 0$ and yields a stable holomorphic reduction $P_{L(\xi)} \subset P^c$.

Let $\chi: L \rightarrow \mathbb{C}^*$ be a character. We need to show $c_1(\chi(P_{L(\xi)})) = 0$. Decompose $\xi_0 = \sum_{j=1}^r x_j \check{t}_j$ as above and denote
				$$S := \{j\,|\, x_j > 0\}.$$
Since $\dot{\chi}: \mathfrak{l}(\xi) \rightarrow \mathbb{C}$ vanishes on $[\mathfrak{l}(\xi), \mathfrak{l}(\xi)]$, it vanishes on all the root spaces $\mathfrak{g}_{\alpha}$ belonging to $\mathfrak{l}(\xi)$ and the dual vectors $t_{\alpha} \in \mathfrak{t}$. In particular, $\dot{\chi}$ vanishes on the simple roots $t_j$ with $j \notin S$ and hence has the shape
				$$\dot{\chi}(\eta) = \sum_{j\in S} \textbf{i} r_j \langle \eta, \check{t}_j \rangle$$ 
for some $r_j \in \mathbb{R}$. Chern-Weyl theory yields
		$$c_1(\chi(P_{L(\xi)})) = \frac{\textbf{i}}{2\pi} \int_{\Sigma} \dot{\chi}(F_B) = \frac{\textbf{i}}{2\pi} \sum_{j\in S} \textbf{i}r_j \int_{\Sigma} \langle *F_B, \check{t}_j \rangle$$
for some connection $B \in \mathcal{A}(P_{K(\xi)})$. We claim that each summand vanishes separately in the last expression. This follows from the assumption
			$$0 = w_0(A,\xi) = \sum_{j=1}^r x_j \int_{\Sigma} \langle *F_B, \check{t}_j \rangle\, dvol_{\Sigma} = \sum_{j \in S} x_j \int_{\Sigma} \langle *F_B, \check{t}_j \rangle\, dvol_{\Sigma}$$
and 
			$$w_0(A,\check{t}_j) = \int_{\Sigma} \langle *F_B, \check{t}_j \rangle \, dvol_{\Sigma} \geq 0$$
since $P^c$ is semistable.
\end{proof}

\subsection{The moment weight inequality}

The moment-weight inequality provides a lower bound for the norm of the moment-map $\mu_{\tau}(A) = *F_A - \tau$ on the complexified orbit $\mathcal{G}^c(A)$. 

\begin{Theorem}[\textbf{The moment-weight inequality}]\label{thmMWI}
Let $P \rightarrow \Sigma$ be a principal $G$ bundle of central type $\tau \in Z(\mathfrak{g})$ defined by (\ref{eqCT}). Let $A \in \mathcal{A}(P)$ be a smooth connection and $\xi \in W^{1,2}(\Sigma,\text{ad}(P))$. Then
			\begin{align} \label{eqMW} - \frac{w_{\tau}(A,\xi)}{||\xi||_{L^2}} \leq \inf_{g \in \mathcal{G}^c} ||*F_{g(A)} - \tau||_{L^2} .\end{align}
\end{Theorem}

The moment weight-inequality is essentially proven by Atiyah and Bott (\cite{AtBott:YangMillsEq}, Prop. 8.13 and Prop. 10.13). They explicitly determine the infimum of the Yang-Mills functional over $\mathcal{G}^c(A)$ in terms of the Harder-Narasimhan filtration of the holomorphic vector bundle $\text{ad}(P^c)$. It follows from the proof of the dominant weight theorem (Theorem \ref{ThmDomWeight}) in the next section that the same description yields the supremum over the left-hand side whenever it is positiv. We provide a different approach following the arguments in \cite{RobSaGeo} for the finite dimensional case which are essentially due to Chen \cite{Chen:2009, Chen:2008} and Donaldson \cite{Donaldson:2005}.

\begin{proof}
We reduce the proof to the case where $Z(G)$ is discrete and $\tau = 0$. Denote by $\bar{A} \in \mathcal{A}(P/Z_0(G))$ the induced connection on the quotient bundle and decompose $\xi = \xi^{ss} + \xi^z$ as in Lemma \ref{WLemma1}. Let $g \in \mathcal{G}^c(P)$ be given and decompose $F_{gA} = F^{ss} + F^z$ in the same way. Note that $F_{g\bar{A}} = F^{ss}$. Suppose that the moment-weight inequality is satisfied on $P/Z_0(G)$, i.e.
						$$-\frac{w_0(\bar{A},\xi^{ss})}{||\xi^{ss}||} \leq ||*F_{g\bar{A}}||.$$
If $w_{\tau}(A,\xi) \leq 0$ then Lemma \ref{WLemma1} implies
		\begin{align*}
				- \frac{w_{\tau}(A,\xi)}{||\xi||} &\leq - \frac{w_0(\bar{A},\xi^{ss})}{||\xi^{ss}||} \leq ||*F_{g\bar{A}}|| \leq \sqrt{ ||F^{ss}||^2 + ||F^{z} - \tau||^2} \\
																					&= ||F_{gA} - \tau||
		\end{align*}
If $w_{\tau}(A,\xi) > 0$ then (\ref{eqMW}) is trivial and this completes the reduction argument.

Now assume that $Z(G)$ is discrete and $\tau= 0$. Let $\xi \in W^{1,2}(\Sigma, \text{ad}(P))$ with $w_0(A,\xi) < \infty$. Then $\xi$ is smooth by Proposition \ref{PropMaxWeight} and  $e^{\textbf{i}t\xi}A$ converges uniformly to some connection $A_{+}$ by Corollary \ref{CorSplit}. In particular, there exists a constant $C>0$ such that
		\begin{align} \label{MWeq1} \sup_{t\geq 0} ||d_{e^{\textbf{i}t\xi}A} \xi||_{L^2}^2 \leq C. \end{align}
Let $g_0 = u_0 e^{\textbf{i}\eta_0} \in \mathcal{G}^c(P)$ be given and define $\eta(t) \in W^{2,2}(\Sigma,\text{ad}(P))$ and $u(t) \in \mathcal{G}$ by the equation
			$$e^{\textbf{i}\xi t} = e^{\textbf{i}\eta(t)} u(t) g_0.$$ 
From this follows pointwise the estimate
		\begin{align} \label{MWeq2}|\eta(t) - t\xi| \leq |\eta_0|. \end{align}
To see this, denote by $\pi: G^c \rightarrow G^c/G$	the canonical projection and recall that $G^c/G$ is a complete simply-connected Riemannian manifold with nonpositive sectional curvature. For a fixed time $t$ and $z \in \Sigma$ define $p := \pi(e^{\textbf{i}t\xi(z)})$ and $q := \pi(e^{\textbf{i} \eta(t,z)})$. Then 
			$$\gamma:[0,1] \rightarrow G^c/G,\qquad \gamma(s) := \pi(e^{\textbf{i}t\xi(z)} e^{-\textbf{i}s \eta_0(z)})$$
is the unique geodesic from $p$ to $q$ in $G^c/G$ of length $|\eta_0(z)|$. It is well-known that the exponential map on a complete simply-connected Riemannian manifold is distance increasing. This yields
			$$|\eta(t,z) - t\xi(z)| \leq \text{dist}_{G^c/G}(p,q) = |\eta_0(z)|$$
and hence (\ref{MWeq2}). With this estimate we get
		\begin{align*}
				\left|\left| \frac{\xi}{||\xi||} - \frac{\eta(t)}{||\eta(t)||} \right|\right|
						&\leq \left|\left| \frac{t\xi - \eta(t)}{t||\xi||} + \frac{\eta(t)}{t||\xi||} - \frac{\eta(t)}{||\eta(t)||} \right|\right| \\
						&\leq \frac{||t\xi - \eta(t)||}{t||\xi||} + \frac{|\,||\eta(t)|| - t||\xi||\,|}{t ||\xi||} \\
						&\leq 2 \frac{||\eta_0||}{t||\xi||}
		\end{align*}
and hence
		\begin{align} \label{MWeq3}
			\lim_{t\rightarrow \infty} \left|\left|\frac{\eta(t)}{||\eta(t)||} - \frac{\xi}{||\xi||} \right|\right| = 0.
		\end{align}
By (\ref{MWeq0}) the map
				$$s \mapsto \langle *F_{e^{\textbf{i}su^{-1} \eta u}g_0A} , u^{-1} \eta u \rangle$$
is nondecreasing in $s$. With the relation $e^{\textbf{i}u^{-1} \eta u} g_0 = u^{-1} e^{\textbf{i}t\xi}$ follows
		\begin{align*}
			- ||*F_{g_0 A}|| &\leq \frac{1}{||\eta||} \left\langle * F_{g_0 A}, u^{-1} \eta u \right\rangle 
													\leq \frac{1}{||\eta||} \left\langle *F_{e^{\textbf{i} u^{-1}\eta u}g_0A} , u^{-1} \eta u \right\rangle \\
											 &\leq \frac{1}{||\eta||} \left\langle *F_{u^{-1}e^{\textbf{i}t\xi}A} , u^{-1} \eta u \right\rangle
												= \frac{1}{||\eta||} \left\langle *F_{e^{\textbf{i}t\xi}A}, \eta \right\rangle \\
											 &\leq \frac{\left\langle *F_{e^{\textbf{i}t\xi}A}, \xi \right\rangle}{||\xi||}  + \left\langle *F_{e^{\textbf{i}t\xi}A}, \frac{\eta}{||\eta||} - \frac{\xi}{||\xi||} \right\rangle
		\end{align*}
It follows from (\ref{MWeq1}) and (\ref{MWeq2}) that the right and side converges to $\frac{w_0(A,\xi)}{||\xi||}$ for $t \rightarrow \infty$ and this proves the theorem.

\end{proof}

\section{The Kempf-Ness functional}

Let $G$ be a compact connected Lie group and let $P \rightarrow \Sigma$ be a principal $G$ bundle of central type $\tau \in Z(\mathfrak{g})$ defined by (\ref{eqCT}). Let $A \in \mathcal{A}(P)$ be a smooth connection. The \textbf{Kempf-Ness functional} associated to $A$ is the $\mathcal{G}(P)$-invariant functional 
			\begin{align} \label{KNTeq1}
					\Phi_A : \mathcal{G}^c(P) \rightarrow \mathbb{R}, \qquad \Phi_A(e^{\textbf{i} \xi} u) = \int_0^1 \langle *F_{e^{-\textbf{i}t\xi}A} - \tau, - \xi \rangle\, dt.
			\end{align}
We show in Lemma \ref{KNLemma0} below that the derivative of $\Phi_A$ is given by
			\begin{align} \label{KNTeq2}
					\alpha_A(g ; \hat{g}) = - \langle *F_{g^{-1}A} - \tau, \text{Im}(g^{-1} \hat{g}) \rangle.
			\end{align}
The asymptotic slope of $\Phi_A$ along the geodesic ray $t \mapsto e^{\textbf{i}t\xi}$ yields the weight $w_{\tau}(A, -\xi)$. This is related to the stability of the associated holomorphic principal bundle $(P^c, J_A)$ by Proposition \ref{PropStabChar}. On the other hand, it follows directly from (\ref{KNTeq2}) that $g \in \mathcal{G}^c(P)$ is a critical point of $\Phi_A$ if and only if $*F_{g^{-1}A} = \tau$. The analog of the Kempf-Ness theorem in classical GIT is Theorem \ref{ThmKNgen} below. It characterizes the different notions of $\mu_{\tau}$-stability in terms of the global behaviour of $\Phi_A$ and thus provides a link between the algebraic and the symplectic notion of stability. We can deduce from this the Narasimhan-Seshadri-Ramanathan theorem in the second subsection.

\subsection{The generalized Kempf-Ness theorem}

\begin{Lemma} \label{KNLemma0}
Let $P \rightarrow \Sigma$ be a principal $G$ bundle and define $\Phi_A : \mathcal{G}^c(P) \rightarrow \mathbb{R}$ by (\ref{KNTeq1}).
		\begin{enumerate}
				\item The derivative of $\Phi_A$ is given by
							$$\alpha_A(g ; \hat{g}) = - \langle *F_{g^{-1}A} - \tau, \text{Im}(g^{-1} \hat{g}) \rangle.$$
				\item Let $g, h \in \mathcal{G}^c(P)$, then
							$$\Phi_{h^{-1}A}(h^{-1}g) = \Phi_A(g)- \Phi_A(h).$$
		\end{enumerate}
\end{Lemma}

\begin{proof}
Let $g \in \mathcal{G}^c(P)$, $\hat{g} \in T_g \mathcal{G}^c(P)$ and let $u \in \mathcal{G}(P)$ be given. Then
		\begin{align*}
				\alpha_A(gu^{-1}, \hat{g}u^{-1}) &= \langle * F_{u g^{-1}A}, \text{Im} (u g^{-1}\hat{g}u^{-1}) \rangle \\
																				 &= \langle u * F_{g^{-1}A} u^{-1}, u \text{Im}(g^{-1} \hat{g}) u^{-1} \rangle \\
																				 &= \alpha_A(g, \hat{g})
		\end{align*}
shows that $\alpha_A$	is invariant under the right-action of $\mathcal{G}(P)$ and hence descends to a $1$-form on $\mathcal{G}^c(P)/\mathcal{G}(P)$. 

We claim that $\alpha_A$ is closed. Denote by $\pi: G^c \rightarrow G^c/G$ the canoncial projection and let $\hat{g}_1 = d\pi(g)g\textbf{i}\xi$ and $\hat{g}_2 := d\pi(g) g\textbf{i} \eta$ be two tangent vectors in $T_{\pi(g)} \mathcal{G}^c(P)/\mathcal{G}(P)$. Then
			\begin{align*}
					d\alpha_A(g; \hat{g}_1, \hat{g}_2) 
																				&= d \alpha_A(g; \hat{g}_2) [\hat{g}_1] - d \alpha_A(g; \hat{g}_1) [\hat{g}_2] - \alpha_A(g; [\hat{g}_1, \hat{g}_2]) \\
																				&= d \langle F_{g^{-1}A} - \tau, \eta \rangle [ g \textbf{i}\xi] - d \langle F_{g^{-1}A} - \tau, \xi \rangle [ g \textbf{i}\eta] \\
																				&= \langle d_{g^{-1}A}^* d_{g^{-1}A} \xi, \eta  \rangle - \langle d_{g^{-1}A}^* d_{g^{-1}A} \eta, \xi  \rangle = 0.
			\end{align*}
We used in the second step that $[\hat{g}_1, \hat{g}_2] \in T_{g}\mathcal{G}(P)$ is tangent to the real gauge orbit and thus lies in the kernel of $\alpha_A(g;\cdot)$.

Denote for $p,q \in \mathcal{G}^c(P)/\mathcal{G}(P)$ by $[p,q]$ the geodesic segment connecting $p$  to $q$. Then (\ref{KNTeq1}) can be reformulated as
						\begin{align} \label{KNalt} \Phi_A(g) = \int_{[\pi(\mathds{1}), \pi(g)]} \alpha_A. \end{align}
For $h \in \mathcal{G}^c(P)$ we have $\alpha_{h^{-1} A}(h^{-1}g, h^{-1} \hat{g}) = \alpha_A(g, \hat{g})$ and hence
						$$\Phi_{h^{-1}A}(h^{-1}g) = \int_{[\pi(\mathds{1}), \pi(h^{-1}g)]} \alpha_{h^{-1}A} = \int_{[\pi(h), \pi(g)]} \alpha_A.$$
Since $\alpha_A$ is closed we have
						$$\int_{[\pi(h), \pi(g)]} \alpha_A = \int_{[\pi(\mathds{1}), \pi(g)]} \alpha_A - \int_{[\pi(\mathds{1}), \pi(h)]} \alpha_A = \Phi_A(g) - \Phi_A(h)$$
and this establishes the second part of the lemma.

Using the second part, we can can reduce the proof of the first part to the case $g = \mathds{1}$ and in this case the claim follows directly from (\ref{KNTeq1}).

\end{proof}

The difficult part of the following theorem is the stable case. The proof of this case is due to Bradlow \cite{Bradlow:1991} and Mundet \cite{Mundet:2000} in the context of more general moduli problems.

\begin{Theorem}[\textbf{Generalized Kempf-Ness Theorem}] \label{ThmKNgen}
Let $G$ be a compact connected Lie group, let $P \rightarrow \Sigma$ be a principal $G$ bundle with central type $\tau \in Z(\mathfrak{g})$ defined by (\ref{eqCT}) and let $A \in \mathcal{A}(P)$.
		\begin{enumerate}
				\item $A$ is $\mu_{\tau}$-stable if and only if $\mathcal{G}^c(A)$ has discrete $\mathcal{G}^c(P)$ isotropy and for every $R > 0$ such that
									$$M_R := \{ \xi \in W^{2,2}(\Sigma, \text{ad}(P)\,|\, ||*F_{e^{-\textbf{i}\xi}A} - \tau ||_{L^2} \leq R \}$$
				is nonempty, there exist constants $c_1, c_2 > 0$ such that
								\begin{align} \label{KNeq3} \Phi_A(e^{\textbf{i}\xi}) \leq c_1||\xi||_{L^{\infty}} + c_2 \qquad \text{for all $\xi \in M_R$.} \end{align}
									
				\item	$A$ is $\mu_{\tau}$-polystable if and only if	$\Phi_A$ has a critical point.
				
				\item $A$ is $\mu_{\tau}$-semistable if and only if $\Phi_A$ is bounded below.
				
				\item $A$ is $\mu_{\tau}$-unstable if and only if $\Phi_A$ is unbounded below.
		\end{enumerate}
\end{Theorem}

\begin{proof}
We consider both implications of the stable case in the following lemmas first. The proof will then be given on page \pageref{ProofThmKN} below.
\end{proof}

\begin{Lemma} \label{KNLemma1}
Assume the setting of Theorem \ref{ThmKNgen}. Suppose that the orbit $\mathcal{G}^c(A) \subset \mathcal{A}^*(P)$ contains only irreducible connections and that there exist $c_1,c_2, R > 0$ such that $M_R$ is nonempty and (\ref{KNeq3}) holds. Then exists $\xi_0 \in M_R$ such that 
						\begin{align} \label{KNeq4} \Phi_A(e^{\textbf{i}\xi_0}) \leq \Phi_A(e^{\textbf{i}\xi}) \qquad \text{for all $\xi \in M_R$} \end{align}
and $B := e^{-\textbf{i}\xi_0} A$ satisfies $F_B = 0$.
\end{Lemma}

\begin{proof}
Suppose first that $\xi_0 \in M_R$ satisfies (\ref{KNeq4}). Let $B := e^{-\textbf{i}\xi_0}A$ and let $\eta \in W^{2,2}(\Sigma, \text{ad}(P))$ be a solution of the equation
				$$\Delta_B \eta = d_B^* d_B \eta = * F_B$$
which exists since $B$ is irreducible. Then follows
		\begin{align*}
				\left.\frac{d}{dt}\right|_{t=0} \Phi_A(e^{\textbf{i}\xi_0} e^{\textbf{i}\eta t}) 
										= \alpha_A( e^{\textbf{i}\xi_0}, e^{\textbf{i}\xi_0}\textbf{i}\eta ) 
										= - \langle *F_B, \eta \rangle = - ||d_B \eta ||^2
		\end{align*}
and
		\begin{align*}
				\left.\frac{d}{dt}\right|_{t=0} \left|\left| *F_{e^{-\textbf{i}\eta t} e^{-\textbf{i}\xi_0}A}\right|\right|^2
												&= 2 \left\langle *F_B, * \left.\frac{d}{dt}\right|_{t=0} F_{e^{-\textbf{i}\eta t} B} \right\rangle \\
																= 2 \langle *F_B , * d_B * d_B \eta \rangle 
												&= -2 \langle *F_B, \Delta_B \eta \rangle = - 2 || * F_B||^2
		\end{align*}
Now decompose $e^{\textbf{i}\xi_0} e^{\textbf{i}\eta t} = e^{\textbf{i}\xi_1}u$. Then the calculation shows that for sufficiently small $t$ we have $\xi_1 \in M_B$ and $\Phi_A(e^{\textbf{i}\xi_1}) \leq \Phi_A(e^{\textbf{i}\xi_0})$ with equality if and only if $F_B = 0$. Since (\ref{KNeq4}) yields the converse inequality, we have indeed equality and hence $F_B = 0$.

It remains to prove the existence of a minimizer $\xi_0 \in M_R$. Let $\{\xi_k\} \subset M_R$ be minimizing sequence satisfying
				\begin{align} \label{KNeq5}
						\lim_{k\rightarrow \infty} \Phi_A(e^{\textbf{i}\xi_k}) = \inf_{\xi \in M_R} \Phi_A(e^{\textbf{i}\xi}).
				\end{align}
By definition of $M_R$, the curvature $F_{e^{\textbf{i}\xi_k} A}$ is uniformly bounded in $L^2$. Hence the Uhlenbeck compactness theorem asserts that there exists $u_k \in \mathcal{G}(P)$ such that $A_k := u_ke^{\textbf{i}\xi_k} A$ converges weakly in $W^{1,2}$. For $g_k := u_k e^{\textbf{i}\xi_k}$ the expression
					$$\bar{\partial}_{A_k} - \bar{\partial}_A = g_k^{-1} \bar{\partial}_A {g_k}$$
is thus uniformly bounded in $W^{1,2}$. Since $\xi_k$ is uniformly bounded in $L^{\infty}$ by (\ref{KNeq3}) and (\ref{KNeq5}), we conclude that $g_k$ is uniformly bounded in $W^{2,2}$. Then $e^{2\textbf{i}\xi_k} = g_k^*g_k$ is also uniformly bounded in $W^{2,2}$ and since $\xi_k$ is uniformly bounded in $L^{\infty}$ this implies that $\xi_k$ is uniformly bounded in $W^{2,2}$. Hence, after taking a subsequence, there exists $\xi_0 \in M_R$ such that $\xi_k$ converges to $\xi_0$ weakly in $W^{2,2}$ and strongly in $W^{1,p}$ for every $p < \infty$. From this follows
			$$\lim_{k \rightarrow \infty} \langle * F_{e^{-\textbf{i}t\xi_k}A}, - \xi_k \rangle = \langle *F_{e^{-\textbf{i}t\xi_0}} A, - \xi_0 \rangle.$$
Hence $\lim_{k \rightarrow \infty} \Phi_A(e^{\textbf{i}\xi_k})	= \Phi_A(e^{\textbf{i}\xi_0})$ and $\xi_0$ satisfies (\ref{KNeq4}).
\end{proof}

\begin{Lemma} \label{KNLemma2}
Assume the setting of Theorem \ref{ThmKNgen}. Suppose that $Z_0(G)$ is discrete, $\tau = 0$ and $w_0(A,\xi) > 0$ for all nonzero $\xi \in W^{1,2}(\Sigma, \text{ad}(P))$. Let $R > 0$ be given such that $M_R$ is nonempty. Then exist constants $c_1, c_2 > 0$ such that (\ref{KNeq3}) is satisfied.
\end{Lemma}

\begin{proof}
The proof consists of several steps.\\

\textbf{Step 1:} \textit{There exists $C > 0$ such that 
											$$||\xi||_{C^0} \leq C \left(||\xi||_{L^1} + 1 \right) \qquad \text{for all $\xi \in M_R$.}$$}

We observe that
	\begin{align*}
			\langle *F_{e^{\textbf{i}\xi}A} - * F_A, \xi \rangle 
						&= \int_0^1 \langle \Delta_{e^{\textbf{i}t\xi}A} \xi, \xi \rangle \, dt = \Delta ||\xi||^2 + 2 \int_0^1 ||d_{e^{\textbf{i}t\xi}A} \xi||^2\, dt \\
						&\geq \Delta ||\xi||^2 \geq 2 ||\xi|| \Delta ||\xi||
	\end{align*}
and hence
	$$\Delta ||\xi|| \leq \frac{1}{2} || *F_{e^{\textbf{i}\xi}A} - * F_A ||  \leq \frac{1}{2}\left(||*F_A|| + R\right)$$
by definition of $M_R$. Hence there exists a constant $c_0 > 0$ for which the following pointwise estimate holds
					\begin{align} \label{KNeq7} \Delta ||\xi|| \leq c_0. \end{align} 
An argument due to Simpson (\cite{Simpson:1987}, Prop 2.1) shows that this implies the claim. For this denote
			$$f: \Sigma \rightarrow \mathbb{R}, \qquad f(z) := ||\xi(z)||.$$ 
For $z_0 \in \Sigma$ choose a local coordinate which identifies $z_0$ with the origin in $\mathbb{C}$. Let $B_{r_0}(0)$ be a ball contained in the image of this local coordinate and let $r \in (0,r_0)$. Let $w, h$ be solutions of 
			$$\Delta w = c_0,\, w|_{\partial B_r(0)} = 0 \qquad \Delta h = 0,\,  h|_{\partial B_r(0)} = f|_{B_r(0)}.$$
Here we consider the Laplacian of $\Sigma$ which agrees with the Laplacian on $\mathbb{C}$ up to a positive factor. Hence (\ref{KNeq7}) and the maximum principle show that $f - w - h \leq 0$ and the mean value theorem yields
		$$f(0) - w(0) \leq h(0) = \frac{1}{2\pi r} \int_{\partial B_r(0)} f.$$
Moreover, elliptic regularity and the Sobolev embedding $W^{2,2} \hookrightarrow C^0$ yield
					$$|w(0)| \leq C ||w||_{W^{2,2}} \leq C ||\Delta w||_{L^2} = C r$$ 
and hence
		$$f(z_0) \leq C \left( r + \frac{1}{r} \int_{\partial B_r(0)} \tilde{f} \right).$$
Now choose $r \in (r_0/2 , r_0)$ such that $\pi r_0 \int_{\partial B_r(0)} \tilde{f} \leq ||f||_{L^1}$ holds. Then follows
				$$f(z_0) \leq C\left(r_0 + \frac{1}{r_0^2} ||f||_{L^1} \right).$$	
Since $\Sigma$ is compact, we can perform this argument within finitely many charts and choose the final constant $C$ to be independent of $z_0$.\\

\textbf{Step 2:} \textit{There exist $c_1, c_2 > 0$ such that 
				$$||\xi||_{L^1} \leq c_1 \Phi_A(e^{\textbf{i}\xi}) + c_2 \qquad \text{for all $\xi \in M_R$.}$$}

Suppose the claim is false. Then exists $C_k > 0$ and $\xi_k \in M_R$ such that
		$$\lim_{k \rightarrow \infty} C_k = \infty, \quad \lim_{k \rightarrow \infty} ||\xi_k||_{L^1} = \infty \quad \text{and} \quad \text{ $||\xi_k||_{L^1} \geq C_k \Phi_A(e^{\textbf{i}\xi_k})$.}$$
It follows from Step 2 that $\eta_k := - \xi_k/||\xi_k||_{L^1}$ is uniformly bounded in $L^{\infty}$. Denote $\ell_k := ||\xi_k||_{L^2}$. Then
		\begin{align*}
			\frac{1}{C_k} \geq \frac{\Phi_A(e^{\textbf{i}\xi_k})}{||\xi_k||_{L^1}} = \int_0^1 \langle * F_{e^{\textbf{i}t\xi_k}A} , \eta_k \rangle \, dt = \frac{1}{\ell_k} \int_0^{\ell_k} \langle * F_{e^{\textbf{i}t\eta_k} A}, \eta_k \rangle \, dt
		\end{align*}
The integrand is increasing by (\ref{MWeq0}). Hence, for any fixed $t > 0$ follows
		\begin{align} \label{KNeq8}
				\frac{1}{C_k} \geq \frac{l_k - t}{l_k} \langle *F_{e^{\textbf{i}t\eta_k}}, \eta_k \rangle + \frac{t}{l_k} \langle F_A, \eta_k \rangle.
		\end{align}		
It follows from (\ref{Weq2}) that
		$$\langle *F_{e^{\textbf{i}t\eta_k}A}, \eta_k \rangle = \langle *F_A, \eta_k \rangle + \frac{1}{2} \int_0^t || e^{i\eta_k s} (\bar{\partial}_A \eta_k) e^{-\textbf{i}\eta_k s} ||^2\, ds$$
and, since $\eta_k$ is uniformly bounded in $C^0$, we conclude that $\bar{\partial}_A \eta_k$ is uniformly bounded in $L^2$. Since $A$ is irreducible and $||\bar{\partial}_A\eta_k||^2 = \frac{1}{2} ||d_A \eta_k||^2$ this shows that $\eta_k$ is uniformly bounded in $W^{1,2}$. Hence, after taking a subsequence, there exists $\eta \in W^{1,2} \cap L^{\infty}$ such that $\eta_k \rightarrow \eta$ converges weakly in $W^{1,2}$ and strongly in $L^p$ for every $1 \leq p < \infty$. In particular $||\eta||_{L^1} = 1$ shows that $\eta \neq 0$ and it is straightforward to check that
			$$\lim_{k \rightarrow \infty} \langle *F_{e^{\textbf{i}t\eta_k}A}, \eta_k \rangle = \langle *F_{e^{\textbf{i}t\eta}A}, \eta \rangle.$$
Now (\ref{KNeq8}) implies $\langle *F_{e^{\textbf{i}t\eta}A}, \eta \rangle \leq 0$ and as $t \rightarrow \infty$ we obtain $w_0(A,\eta) \leq 0$. This contradicts our assumptions and proves Step 2.\\

\textbf{Step 3:} \textit{There exist $c_1, c_2 > 0$ such that 
				$$||\xi||_{L^{\infty}} \leq c_1 \Phi_A(e^{\textbf{i}\xi}) + c_2 \qquad \text{for all $\xi \in M_R$.}$$}

This follows directly from Step 1 and Step 2.
		
\end{proof}

\begin{proof}[Proof of Theorem \ref{ThmKNgen}.] \label{ProofThmKN}
Suppose $A$ is $\mu_{\tau}$-stable. Then $Z_0(G)$ is discrete, $\tau = 0$ and $\mathcal{G}^c(A)$ has discrete $\mathcal{G}^c(P)$ isotropy. We claim that $w_0(A,\xi) > 0$ for all $\xi \in W^{1,2}(\Sigma, \text{ad}(P))$. By Proposition \ref{PropStabChar} this condition is equivalent to the stability of the induced holomorphic structure $J_A$ on $P^c := P \times_G G^c$. In particular, this condition is invariant under the action of $\mathcal{G}^c(P)$ and we may assume that $F_A = 0$. Then (\ref{Weq1}) shows 
		$$w_0(A,\xi) = \int_0^{\infty} ||d_{e^{\textbf{i}t\xi}A} \xi||_{L^2}^2\,  dt > 0$$						
since $A$ is irreducible. Thus Lemma \ref{KNLemma2} applies and shows that the estimate (\ref{KNeq3}) is satisfied. The converse direction follows from Lemma \ref{KNLemma1}. 

The characterization of the $\mu_{\tau}$-polystable case follows directly from (\ref{KNTeq2}).

In the following let $A(t)$ denote the solution of the Yang-Mills flow (\ref{YMeq2}) starting at $A$ and let $A_{\infty} := \lim_{t\rightarrow \infty} A(t)$. Suppose $A$ is $\mu_{\tau}$-unstable. Theorem \ref{MLT} and (\ref{YMeq02}) show that 
					$$||F_{gA} - \tau|| \geq ||F_{A_{\infty}} - \tau|| = c > 0$$
for all $g \in \mathcal{G}^c(P)$.	Now define $g(t)$ by (\ref{YMeq3}). Then $A(t) = g(t)^{-1}(A)$ and
			\begin{align*}
					\frac{d}{dt} \Phi_A (g(t)) &= \alpha_A(g(t), \dot{g}(t)) = - \langle *F_{A(t)} - \tau, *F_{A(t)}  \rangle \\
																		 &= - ||*F_{A(t)} - \tau||^2 \leq -c
			\end{align*}
where the penultimate step follows from (\ref{eqCT}). This shows that $\Phi_A$ is unbounded below.

Suppose conversely that $A$ is $\mu_{\tau}$-semistable. It follows from Theorem \ref{MLT} and (\ref{YMeq02}) that $A_{\infty}$ is a global minima for the Yang-Mills functional on $\mathcal{A}(P)$ and $*F_{A_\infty} = \tau$. It follows from the Lojasiewicz inequality (Lemma \ref{LemmaLoj}) that there exists $\gamma \in [\frac{1}{2}, 1)$ and $C, t_0 > 0$ such that 
		\begin{align*}
					||*F_{A(t)} - \tau||^2 &\leq 2 |\mathcal{YM}(A(t)) - \mathcal{YM}(A_{\infty})|^{\gamma}
																					\leq C || d_{A(t)}^* F_{A(t)}||_{L^2} \\
																 &\leq C \frac{||d_{A(t)}^* F_{A(t)}||_{L^2}^2}{(\mathcal{YM}(A(t)) - \mathcal{YM}(A_{\infty}))^{\gamma}} \\
																 &=  \frac{d}{dt} C\left( \mathcal{YM}(A(t)) - \mathcal{YM}(A_{\infty}) \right)^{1-\gamma}					
		\end{align*}																	
for all $t > t_0$. Since the right hand side is integrable, the solution $g(t)$ of (\ref{YMeq3}) satisfies
			$$\lim_{t\rightarrow \infty} \Phi_A(g(t)) = - \int_0^{\infty} ||*F_{A(t)} - \tau||^2\, dt =: a > - \infty.$$
We claim that $a$ is a global minimum for $\Phi_A$. For this let $\tilde{g}_0 \in \mathcal{G}^c(P)$ and let $\tilde{g}(t)$ be the solution of (\ref{YMeq6}) starting at $\tilde{g}_0$. This is a negative gradient flow line of $\Phi_A$ and satisfies
		$$ \frac{d}{dt} \Phi_A(\tilde{g}(t)) = -\alpha_A( \tilde{g}(t), \dot{\tilde{g}}(t) ) = - ||*F_{\tilde{g}(t)^{-1}A} - \tau||^2 \leq 0.$$
Define $\eta(t) \in W^{2,2}(\Sigma, \text{ad}(P))$ and $u(t) \in \mathcal{G}(P)$ by the equation		
		$$g(t) \exp(\textbf{i} \eta(t)) u(t) = \tilde{g}(t)$$
and
			$$\beta_t: [0,1] \rightarrow \mathcal{G}^c(P),\qquad \beta_t(s) = g(t)e^{\textbf{i}s \eta(t)}.$$
Then $(\Phi_A\circ \beta_t)$ satisfies
		\begin{align*}
				\left. \frac{d}{ds}\right|_{s=0} (\Phi_A\circ \beta_t)(s) &= \alpha_A(g(t), \partial_s \beta_t(s) \rangle = - \langle * F_{g(t)^{-1} A} - \tau, \eta \rangle \\
																										&\geq - ||* F_{g(t)^{-1} A} - \tau||\cdot ||\eta(t)||_{L^2}
		\end{align*}
and
		\begin{align*}
				 \frac{d^2}{ds^2} (\Phi_A\circ \beta_t)(s) 
													&=  - \frac{d}{ds} \langle * F_{e^{-\textbf{i}\eta(t)s}g(t)^{-1} A} - \tau, \eta(t) \rangle \\
													&= \langle d_{A_{s,t}}^* d_{A_{s,t}} \eta(t), \eta(t) \rangle = || d_{A_{s,t}} \eta(t)||^2 \geq 0
		\end{align*}
where we abbreviated $A_{s,t} := e^{-\textbf{i}\eta(t)s}g(t)^{-1} A$. In particular, $\Phi_A\circ \beta_t$ is convex and since $\eta(t)$ is uniformly bounded in $L^{\infty}$ by Proposition \ref{YMProp2} there exists a constant $C > 0$ such that
				$$\Phi_A(\tilde{g}(t)) \geq \Phi_A(g(t)) - C ||F_{g(t)^{-1} A} - \tau||.$$
Since $\Phi_A(\tilde{g}_0) \geq \Phi_A(\tilde{g}(t))$ for all $t$ and the right hand side converges to $a$ as $t \rightarrow \infty$ we conclude the $\Phi_A(\tilde{g}_0) \geq a$. This establishes the claim and completes the proof of the theorem.
\end{proof}

\subsection{The Narasimhan-Seshadri-Ramanathan theorem}

The Narasimhan-Seshadri-Ramanathan theorem relates the notion of stable objects in Definition \ref{DefnStab2} and Definition \ref{DefnStab3}. This was first proven by Narasimhan-Seshadri \cite{NarasimhanSeshadri:1965} in the case $G=U(n)$ and later extended by Ramanathan \cite{Ramanathan:1975} to general compact connected Lie groups. Both of these proofs are entirely of algebraic geometric nature. 

In the case $G  = U(n)$ Donaldson \cite{Donaldson:1983NS} gave an analytic proof of this result. His argument uses the moment weight inequality and an induction argument which is based on the Harder-Narasimhan filtration. We present a different proof which is due to Bradlow \cite{Bradlow:1991} and Mundet \cite{Mundet:2000}. The main step in their proof consists of establishing the stable case in Theorem \ref{ThmKNgen}.

\begin{Theorem} [\textbf{Narasimhan-Seshadri-Ramanathan}] \label{ThmNSR}
Let $G$ be a compact connected Lie group and $P \rightarrow \Sigma$ a principal $G$ bundle with central type $\tau \in Z(\mathfrak{g})$ defined by (\ref{eqCT}). Let $A \in \mathcal{A}(P)$ and consider the complexified bundle $P^c := P \times_G G^c$ with the holomorphic structure induced by $A$. Then $P^c$ is stable if and only if there exists a complex gauge transformation $g \in \mathcal{G}^c(P)$ such that $*F_{gA} = \tau$ and the kernel of 
							$$L_{A}: W^{2,2}(\Sigma, \text{ad}(P^c)) \rightarrow W^{1,2}(\Sigma, T^*\Sigma\otimes\text{ad}(P))$$ 
										$$L_A(\xi + \textbf{i}\eta) = -d_A \xi + * d_A \eta$$
	contains only constant central sections.
\end{Theorem}

\begin{proof}
We may assume by Lemma \ref{LemmaSSred1} and Lemma \ref{LemmaSSred2} that $Z_0(G)$ is discrete and $\tau = 0$.

Suppose there exists $g \in \mathcal{G}^c(P)$ such that $*F_{gA}= 0$ and $gA$ is irreducible. Then (\ref{Weq1}) shows
						$$w_{0}(gA,\xi) = \int_0^{\infty} ||d_{e^{\textbf{i}t\xi}gA} \xi||_{L^2}^2\,  dt > 0$$
for all $0 \neq \xi \in W^{1,2}(\Sigma,\text{ad}(P))$ and by Proposition \ref{PropStabChar} $(P^c, J_{gA})$ is stable. Since the notion of stability is $\mathcal{G}^c(P)$ invariant, $(P^c, J_A)$ is stable.

Assume conversely that $(P^c, J_A)$ is stable. For every $g \in \mathcal{G}^c(P)$ then $(P^c, J_{gA})$ is stable as well and Proposition \ref{PropStabChar} implies $w_0(gA,\xi) > 0$ for every nonzero $\xi \in W^{1,2}(\Sigma, \text{ad}(P))$. In particular, $gA$ is irreducible and Lemma \ref{KNLemma2} is applicable and shows that $A$ is $\mu_0$-stable.
\end{proof}

\section{The dominant weight theorem}

The dominant weight theorem strengthens the moment weight inequality (Theorem \ref{thmMWI}). It shows that there exists (up to scaling) a unique section $\xi \in \Omega^0(\Sigma, \text{ad}(P))$ which yields equality in the moment weight inequality, whenever the right hand side is positive. In particular, it relates the notion of unstable objects in Definition \ref{DefnStab2} and Definition \ref{DefnStab3}.

\begin{Theorem}[\textbf{The dominant weight theorem}] \label{ThmDomWeight}
Let $G$ be a compact connected Lie group, let $P \rightarrow \Sigma$ be a principal $G$ bundle of central type $\tau \in Z(\mathfrak{g})$ defined by (\ref{eqCT}) and let $A\in\mathcal{A}(P)$ be a smooth $\mu_{\tau}$-unstable connection.
	\begin{enumerate}
			\item There exists an element $\hat{\xi} \in \Omega^0(\Sigma, \text{ad}(P))$ such that
						\begin{align} \label{DWTeq0} \sup_{0\neq \xi \in \Omega^0(\Sigma, \text{ad}(P))} -\frac{w_{\tau}(A,\xi)}{||\xi||} = -\frac{w_{\tau}(A,\hat{\xi})}{||\hat{\xi}||} = \inf_{g\in \mathcal{G}^c(P)} ||*F_{gA} - \tau||. \end{align}
			\item The normalized section $\hat{\xi}/||\hat{\xi}||$ is uniquely determined. Moreover, it is rational in the sense that it generates a closed $\mathbb{C}^*$ subgroup of $\mathcal{G}^c(P)$.
						
			\item	If $A_{\infty}$ is the limit of the Yang-Mills flow (\ref{YMeq2}) starting at $A$, then there exists $u \in \mathcal{G}(P)$ such that $\hat{\xi}= u(*F_{A_{\infty}} - \tau)u^{-1}$ satisfies (\ref{DWTeq0}).
	\end{enumerate}
\end{Theorem}

\begin{proof}
This result is essentially contained in the work of Atiyah and Bott. They determine in (\cite{AtBott:YangMillsEq}, Prop. 8.13 and Prop. 10.13) the infimum of the Yang-Mills functional on the complexified orbit $\mathcal{G}^c(A)$ in terms of the Harder-Narasimhan filtration of $\text{ad}(P^c)$. 

Bruasse and Teleman \cite{BruasseTeleman:2003, Bruasse:2006} show in a more general gauge theoretical setting that the supremum over the normalized weights is attained in a unique direction whenever it is positive. This corresponds to the case where $(P, J_A)$ is unstable and they identify again the dominant weight with the Harder-Narasimhan filtration.

We follow these ideas in our proof below, but simplify the arguments considerably by using the moment weight inequality and the analytic properties of the Yang-Mills flow. The proof will be given on page \pageref{proofDWT}.
\end{proof}

A key ingredient in the proof is the Harder-Narasimhan filtration associated to a holomorphic holomorphic vector bundle. We review this first before we proceed to the proof of the dominant weight theorem.

\subsection{The Harder-Narasimhan filtration}
Let $F$ and $G$ be holomorphic vector bundles over a Riemann surface $\Sigma$ and let $\alpha: F \rightarrow G$ be a holomorphic bundle map. The kernel and cokernel of $\alpha$ are in general not well-defined as holomorphic vector bundles and one may think of them as vector bundles with singularities. These considerations lead naturally to the larger category of coherent analytic sheaves on $\Sigma$ which is closed under taking kernels and cokernels. The next Lemma, however, allows us to get away without considering sheaves.

\begin{Lemma} \label{LemmaNS1}
Let $F$ and $G$ be holomorphic vector bundles over a Riemann surface $\Sigma$ and let $\alpha : F \rightarrow G$ be a nonzero holomorphic bundle map. Then there exists a commutative diagram of holomorphic vector bundles and holomorphic bundle maps
		$$\begin{CD}
				0 @>>> F' @>>> F @>>> F'' @>>> 0 \\
					@.			@.  	 @VV{\alpha}V	 		@VV{\beta}V		\\
				0 @<<< G'' @<<< G @<<< G' @<<< 0
		\end{CD}$$
 with exact rows and $\text{rk}(F'')= \text{rk}(G')$, $\det(\beta) \neq 0$ and $c_1(F'') \leq c_1(G')$.
\end{Lemma}

\begin{proof}
This Lemma is most easily understood in the language of analytic sheaves. Denote by $\mathcal{O}$ the sheaf of germs of holomorphic functions on $\Sigma$. There exists a one to one correspondence between holomorphic vector bundles and locally free $\mathcal{O}$-sheafs on $\Sigma$, which associates to a vector bundle its sheaf of holomorphic sections. The homomorphism $\alpha$ induces a homomorphism between the associated sheaves and the sheaf kernel and sheaf image are clearly torsion free subsheaves. Since the stalks of $\mathcal{O}$ are isomorphic to the principal ideal domain $\mathbb{C}[[z]]$, these sheaves are locally free and correspond to the the vector bundles $F'$ and $G'$.

More concretely, this describes $F'$ and $G'$ it terms of their section as follows. A section $s \in \Omega^0(U,F)$ is a section of $F'$ if and only if $s$ solves $\alpha(s) = 0$. For $x \in \Sigma$ exists a neighborhood $U$ of $x$ over which $G$ is trivial and such that holomorphic sections of $G'$ have the shape $t = \alpha(s)/f $, where $f(z) = (x-z)^r \in \mathcal{O}(U)$ is the greatest common divisor of the components of the images $\alpha(\tilde{s})$ for $\tilde{s} \in \Omega^0(U,F)$.
\end{proof}

Recall that we denote for a complex vector bundle $E \rightarrow \Sigma$ by 
				$$\mu(E) := \frac{c_1(E)}{\text{rk}(E)}$$ 
its slope or normalized Chern-class.

\begin{Cor}\label{CorNS1}
Let $F$ and $G$ be holomorphic vector bundles over $\Sigma$.
		\begin{enumerate}
				\item Suppose $F$ is semistable, $G$ is stable and $\mu(F) = \mu(G)$. Then any nonzero holomorphic bundle map $\alpha: F \rightarrow G$ is surjective.
				\item Suppose $F$ and $G$ are stable and $\mu(F) = \mu(G)$. Then any nonzero holomorphic bundle map $\alpha: F \rightarrow G$ is an isomorphism.
				\item Suppose $F$ and $G$ are semistable and $\mu(F) > \mu(G)$. Then every holomorphic bundle map $\alpha: F \rightarrow G$ vanishes.
		\end{enumerate}
\end{Cor}

\begin{proof}
We prove the first part. Suppose $\alpha: F \rightarrow G$ is neither zero nor surjective. Using the notation of Lemma \ref{LemmaNS1} we see that $G'$ is a proper subbundle and thus
		$$\mu(G) > \mu(G') \geq \mu(F'') \geq \mu(F)$$
contradicting the assumption $\mu(G) = \mu(F)$. In other two parts follow from a similar argument.
\end{proof}

\begin{Lemma}\label{LemHNstable}
Let $E$ be a holomorphic semistable vector bundle. Then there exists a filtration
						$$0 < E_1 < E_2 < \cdots < E_r = E$$
such that each quotient $E_j/E_{j-1}$ is stable and $\mu(E_j/E_{j - 1}) = \mu(E)$.
\end{Lemma}

\begin{proof}
Let $F \subset E$ be a stable subbundle with $\mu(F) = \mu(E)$. Since $E \cong F \oplus (E/F)$ as $C^{\infty}$-bundles, it follows $\mu(E/F) = \mu(E)$. Moreover, any holomorphic subbundle $G \subset E/F$ with $\mu(G) > \mu(E)$ would lift under the projection map $E \rightarrow E/F$ to a holomorphic subbundle $\tilde{G} \subset E$ with $\mu(\tilde{G}) > \mu(E)$ and this contradicts the semistability of $E$. Hence $E/F$ is semistable and the lemma follows by induction.
\end{proof}

The Harder-Narasimhan filtration generalizes Lemma \ref{LemHNstable} to general holomorphic vector bundles.

\begin{Proposition}[\textbf{Harder-Narasimhan filtration}]
Let $E$ be an holomorphic vector bundle. Then there exists a unique holomorphic filtration
				$$0 = E_0 < E_1 < \cdots < E_r = E$$
such that all quotients $E_i/E_{i-1}$ are semistable and the slopes
					$$\mu_j := \frac{c_1(E_j/E_{j-1})}{\text{rk}(E_j/E_{j-1})}$$
satisfy $\mu_1 > \mu_2 > \cdots > \mu_r$.					
\end{Proposition}

\begin{proof}
The degree of any holomorphic subbundles of $E$ is uniformly bounded by Lemma \ref{LemmaNS2} below. Let $E_1 \subset E$ be a semistable subbundle for which $\mu(E_1) =:\mu_1$ is maximal and such that $E_1$ has maximal rank among all such subbundles. We claim that every proper holomorphic subbundle $G' \subset E/E_1$ satisfies $\mu(G') < \mu_1$. Otherwise, the preimage of $G'$ under the projection $E \rightarrow E/E_1$ would be a subbundle $\tilde{G} \subset E$ with $\mu(\tilde{G}) \geq \mu_1$ and of strictly greater rank then $E_1$. This proves the claim and the existence of the Harder-Narasimhan filtration follows by induction.

Let $0 = \tilde{E}_0 < \tilde{E}_1 < \cdots < \tilde{E}_{\ell} = E$ be another filtration of $E$ such that all quotients $\tilde{E}_{j}/\tilde{E}_{j-1}$ are semistable and the slopes $\tilde{\mu}_j := \mu(\tilde{E}_j/\tilde{E}_{j-1})$ are strictly decreasing. In particular, $\tilde{E}_1$ is semistable and the construction above shows 
				$$\mu(E_1) \geq \mu(\tilde{E}_1) = \tilde{\mu}_1 > \tilde{\mu}_2 > \cdots > \tilde{\mu}_{\ell}.$$ 
The last part of Corollary \ref{CorNS1} shows that the projection $E_1 \rightarrow E/\tilde{E}_{\ell - 1}$ must be zero, since $\mu(E_1) > \tilde{\mu}_{\ell}$ and hence $E_1 \subset \tilde{E}_{\ell -1}$. Repeating the argument, it follows by induction that $E_1 \subset \tilde{E}_j$ for all $j \geq 1$. If $\mu(E_1) > \mu(\tilde{E}_1)$, we could go one step further and obtain the contradiction $E_1 \subset \tilde{E}_0 = 0$. This shows $\mu_1 = \tilde{\mu}_1$. Finally, consider the projection
						$$\alpha: E_1 \rightarrow E \rightarrow E/\tilde{E}_1.$$
If it is nonzero, we can apply Lemma \ref{LemmaNS1} with $F = E_1$ and $G= E/\tilde{E}_1$ to obtain the contradiction
						$$\mu_1 = \mu(E_1) \leq \mu(F'') \leq \mu(G') \leq \tilde{\mu}_2 < \tilde{\mu}_1 = \mu_1.$$
This shows $E_1 \subset \tilde{E}_1$ and by maximality of $\text{rk}(E_1)$ equality must hold. The uniqueness of the Harder-Narasimhan filtration follows now by induction.
\end{proof}

\begin{Lemma} \label{LemmaNS2}
Let $(E, h)$ be a hermitian holomorphic vector bundle over $\Sigma$ and denote by $A \in \mathcal{A}(E)$ the associated unitary connection from Lemma \ref{PreLem1}. For a holomorphic subbundle $F \subset E$ the following holds:
		\begin{enumerate}
			\item Let $E = F \oplus G$ be an orthogonal decomposition and identify $G$ with $E/F$. Denote by $A_F$ and $A_G$ the induced connections on $F$ and $G$. Then $A$ has the shape
							$$A = \begin{pmatrix} A_F & \eta \\ -\eta^* & A_G \end{pmatrix}$$
			with $\eta \in \Omega^{0,1}(\Sigma, \text{End}(G,F))$. Moreover, the curvature has the shape
							$$F_A = \begin{pmatrix} F_{A_F}	- \eta \wedge \eta^*	&		d_A\eta \\ -d_A\eta^* & F_{A_G} - \eta^*\wedge\eta \end{pmatrix}.$$
							
			\item	There exists a constant $C > 0$, which does not depend on $F$, such that
							$$ c_1(F) \leq C(1 - ||\eta||_{L^2}^2 ).$$
		\end{enumerate}
\end{Lemma}

\begin{proof}
We leave the first part as an exercise to the reader, see e.g. \cite{GriffithsHarris} Chapter 0.5. For the second part, we calculate
		\begin{align*}
			c_1(F) = \frac{\textbf{i}}{2\pi} \int_{\Sigma} \text{tr}(F_{A_F}) 
						 = \frac{\textbf{i}}{2\pi} \int_{\Sigma} \text{tr}\left(\left.F_{A}\right|_{F}\right) + \text{tr}\left(\eta \wedge \eta^* \right). 
		\end{align*}
In local coordinates write $\eta = \tilde{\eta} d\bar{z}$ and hence $\eta \wedge \eta^* = 2\textbf{i} \tilde{\eta}\tilde{\eta}^* dx \wedge dy$. This yields precisely the $L^2$-norm of $\eta$. Since $F_A$ is uniformly bounded in $L^{\infty}$, the estimate follows.
\end{proof}

We show next that the Harder-Narasimhan filtration is maximal among all holomorphic filtrations in a certain sense. For this we need to introduce some notation. Let 
		$$\mathcal{E}: 0 = E_0 < E_1 < \cdots < E_r = E$$
be a holomorphic filtration of $E$. Denote $n_j := \text{rk}(E_j/E_{j-1})$, $k_j := c_1(E_j/E_{j-1})$ and define the characteristic vector of the filtration $\mathcal{E}$ to be
		\begin{align} \label{eqFV}
				\vec{\mu}(\mathcal{E}) = \left( \frac{k_1}{n_1}, \ldots, \frac{k_1}{n_1}, \ldots, \frac{k_r}{n_r}, \ldots, \frac{k_r}{n_r} \right) \in \mathbb{R}^n
		 \end{align}
where we repeat each entry $k_j/n_j$ exactly $n_j$-times. Moreover define
			$$\ell_{\mathcal{E}}: \{0,\ldots, n\} \rightarrow \mathbb{R}, \qquad \ell_{\mathcal{E}}(m) = \sum_{j=1}^m \left[\vec{\mu}(\mathcal{E}) \right]_j$$
where $\left[\vec{\mu}(\mathcal{E}) \right]_j$ denotes the $j$-th entry of the vector $\vec{\mu}(\mathcal{E})$. The graph of $\ell_{\mathcal{E}}$ interpolates linearly between the points $(0,0)$, $(n_1,k_1)$, $(n_1 + n_2, k_1 + k_2)$, $\ldots$, $(n,k)$. We consider the following ordering on the space of holomorphic filtrations:
			$$\mathcal{E} \geq \mathcal{F} \quad \text{if and only if} \quad \ell_{\mathcal{E}} \geq \ell_{\mathcal{F}}.$$
We call a filtration $\mathcal{E}$ concave, if the function $\ell_\mathcal{E}$ is concave, or equivalently, if the entries of $\vec{\mu}(\mathcal{E})$ are decreasing.

\begin{Proposition}
Let $E$ be a holomorphic vector bundle over $\Sigma$. The Harder-Narasimhan filtration of $E$ is the unqiue maximal concave filtration on $E$.
\end{Proposition}

\begin{proof}
Let
		$$\mathcal{E}_{HN}: \qquad 0 < E_1 < E_2 < \cdots < E_r = E$$
be the Harder-Narasimhan filtration of $E$ and let $F < E$ be a holomorphic subbundle. It suffices to prove that the point $p_F := (\text{rk}(F), c_1(F))$ lies on or below the graph of $\ell_{\mathcal{E}}$. We prove this by induction on $r$.

Suppose $r = 1$. Then $\mathcal{E}$ is semistable and $\mu(E) \geq \mu(F)$. In particular, $\ell_\mathcal{E}$ is a straight line of slope $\mu(E)$ and $p_F$ clearly lies below that line. 

Suppose now $r > 1$. The Harder-Narasimhan filtration of $E/E_1$ is given by
			$$\mathcal{E}_{HN}': \qquad 0 < E_2/E_1 < E_3/E_1 < \cdots < E_r/E_1 = E/E_1$$
and the induction hypothesis applies to $\mathcal{E}_{HN}'$. Consider the commutative diagram from Lemma \ref{LemmaNS1}
			$$\begin{CD}
				0 @>>> F' @>>> F @>>> F'' @>>> 0 \\
					@.			@.  	 @VV{\alpha}V	 		@VV{\beta}V		\\
				0 @<<< G'' @<<< E/E_1 @<<< G' @<<< 0
		\end{CD}$$
with $\alpha: F \rightarrow E \rightarrow E/E_1$. By the induction hypothesis, the point of $(\text{rk}(G'),c_1(G'))$ lies below $\ell_{\mathcal{E}'}$. Since $\text{rk}(F'') = \text{rk}(G')$ and $c_1(F'') \leq c_1(G')$ the same holds with $G'$ replaced by $F''$. This shows
	\begin{align}
			\label{eqHNmax1} c_1(E_1) + c_1(F'') \leq \ell_{\mathcal{E}}(\text{rk}(E_1) + \text{rk}(F'')).
	\end{align}
Since $F'$ gets maped to zero under $\alpha$, we have $F' \subset E_1$ and $\mu(F') \leq  \mu(E_1)$ by semistability of $E_1$. This shows $c_1(F') \leq \ell_{\mathcal{E}}(\text{rk}(F'))$ and with (\ref{eqHNmax1}) follows
		$$c_1(F) = c_1(F') + c_1(F'') \leq \ell_{\mathcal{E}}(\text{rk}(E_1) + \text{rk}(F'')) + \ell_{\mathcal{E}}(\text{rk}(F')) - \ell_{\mathcal{E}}(\text{rk}(E_1)).$$
Since $\ell_{\mathcal{E}}$ is concave and $\text{rk}(F') \leq \text{rk}(E_1)$ we have
		$$\ell_{\mathcal{E}}(\text{rk}(E_1) + \text{rk}(F'')) - \ell_{\mathcal{E}}(\text{rk}(E_1)) \leq  \ell_{\mathcal{E}}(\text{rk}(F') + \text{rk}(F'')) - \ell_{\mathcal{E}}(\text{rk}(F'))$$  
and thus
			$$c_1(F) \leq \ell_{\mathcal{E}}(\text{rk}(F') + \text{rk}(F'')) = \ell_{\mathcal{E}}(\text{rk}(F)).$$
This completes the proof.
\end{proof}

\begin{Cor}\label{HNCor}
Let $E$ be a holomorphic vector bundle over $\Sigma$. Let $\mathcal{E}$ be a concave filtration of $E$ and $\mathcal{E}_{HN}$ the Harder-Narasimhan filtration of $E$. Then follows
			$$||\vec{\mu}(\mathcal{E})||_2 \leq ||\vec{\mu}(\mathcal{E}_{HN})||_2$$
where $||\cdot||_2$ denotes the standard euclidean norm on $\mathbb{R}^n$. Moreover, equality holds if and only if $\mathcal{E} = \mathcal{E}_{HN}$.
\end{Cor}

\begin{proof}
An easy calculation shows that for two concave filtrations with $\mathcal{E}_1 \leq \mathcal{E}_2$ the estimate $||\vec{\mu}(\mathcal{E}_1)||_2 \leq ||\vec{\mu}(\mathcal{E}_2)||_2$ is satisfied. Moreover, equality holds if and only if $\mathcal{E}_1 = \mathcal{E}_2$.
\end{proof}

\subsection{Proof of the dominant weight theorem}
\label{proofDWT}

We proceed now to the proof of Theorem \ref{ThmDomWeight}. We consider first the case $G = U(n)$ and deduce the general case afterwards by choosing a faithful representation $G \hookrightarrow U(n)$.

\subsubsection*{$\mu_{\tau}$-unstable orbits in the unitary case.}

Assume $G = U(n)$ and denote by $E := P\times_G \mathbb{C}^n$ the associated hermitian vector bundle. Note that the constant central type $\tau$ of $P$ is related to the slope of $E$ by the formula	
			\begin{align} \label{TDWeq0} \tau = - 2\pi \textbf{i}\mu(E) \cdot \mathds{1}. \end{align}
If $(E, \bar{\partial}_A)$ is unstable, then Proposition \ref{PropStabChar} implies that there exists a negative weight $w_{\tau}(A,\xi) < 0$ and the moment weight inequality (Theorem \ref{thmMWI}) shows that $A$ is $\mu_{\tau}$-unstable. The following Lemma proves the converse direction.

\begin{Lemma} \label{LemmaSemistableFlow}
Let $A \in \mathcal{A}(E)$ be a unitary connection and suppose $(E,\bar{\partial}_{A})$ is a semistable holomorphic vector bundle. Then the limit $A_{\infty}$ of the Yang-Mills flow $A(t)$ starting at $A$ satisfies
			$$*F_{A_{\infty}} = - 2\pi\textbf{i}\mu(E)\cdot\mathds{1}.$$
\end{Lemma}

\begin{proof}
We show first that the $W^{1,2}$-closure $\overline{\mathcal{G}^c(A)}$ contains a connection $\bar{A}$ with $F_{\bar{A}} = -2\pi \textbf{i}\mu(E) \cdot \mathds{1}$. For this, consider the refined Harder-Narasimhan filtration from Lemma \ref{LemHNstable}
		$$0 < E_1 < E_2 < \cdots < E_r = E$$
with stable quotients $E_j/E_{j-1}$ all having the same slope as $E$. Choose an orthogonal splitting $E = D_1\oplus \cdots \oplus D_r$ such that $E_j = D_1\oplus \cdots \oplus D_j$. With respect to this splitting $\bar{\partial}_{A}$ has the shape
	\begin{align*}
			\bar{\partial}_{A} = 
				\left(\begin{array}{cccc}
									\bar{\partial}_{A_1}  & 			A_{1 2} 					& \ldots & A_{1 r} \\
							 						0    				  & \bar{\partial}_{A_2}		& \ldots & A_{2 r} \\
							 					 \vdots	  		  &      \vdots							&	\ddots & \vdots  \\
							 						0			  			&  				    0						& \ldots & \bar{\partial}_{A_r}
							\end{array}\right).
		\end{align*}
 Define $g_t := \text{diag}(t^{-1},t^{-2}, \ldots, t^{-r})$. Then
		\begin{align*}
			\bar{\partial}_{g_t(A)} = 
				\left(\begin{array}{cccc}
									\bar{\partial}_{A_1}  & 			t A_{1 2} 					& \ldots & t^{r-1}A_{1 r} \\
							 						0    				  & \bar{\partial}_{A_2}		& \ldots & t^{r-2}A_{2 r} \\
							 					 \vdots	  		  &      \vdots							&	\ddots & \vdots  \\
							 						0			  			&  				    0						& \ldots & \bar{\partial}_{A_r}
							\end{array}\right) 
				\rightarrow
				\left(\begin{array}{cccc}
									\bar{\partial}_{A_1}  & 			0				 					& \ldots & 0			 \\
							 						0    				  & \bar{\partial}_{A_2}		& \ldots & 0			 \\
							 					 \vdots	  		  &      \vdots							&	\ddots & \vdots  \\
							 						0			  			&  				    0						& \ldots & \bar{\partial}_{A_r}
							\end{array}\right)
		\end{align*}
as $t \rightarrow 0$. Since $E_j/E_{j-1} \cong (D_j, \bar{\partial}_{A_j})$ are stable holomorphic vector bundles, Theorem \ref{ThmNSR} shows that there exist complex gauge transformations $g_j \in \mathcal{G}^c(D_j)$ such that $\bar{A}_j = g_j(A_j)$ satisfies $*F_{\bar{A}_j} = -2\pi\textbf{i} \mu(D_j)$. Since $\mu(D_j) = \mu(E)$, we conclude that the induced connection $\bar{A} = \bar{A}_1\oplus \cdots \oplus \bar{A}_r$ has curvature $F_{\bar{A}} = -2\pi \textbf{i}\mu(E) \cdot \mathds{1}$.
		
It follows from (\ref{YMeq02}) that $\bar{A}$ minimizes the Yang-Mills functional over $\mathcal{A}_{U(n)}(E)$. The Lemma follows thus from Theorem \ref{MLT} and Theorem \ref{ThmNUgen}.
\end{proof}

\subsubsection*{Proof of Theorem \ref{ThmDomWeight} for $G = U(n)$.}

Let $\xi$ be a section of skew-hermitian endomorphism in $\mathfrak{u}(E) \subset \text{End}(E)$ satisfying $||\xi|| = 1$ and 
			\begin{align} \label{TDWeq1} - w_{\tau}(A,\xi) = \sup_{0 \neq \eta \in \Omega^0(\Sigma, \mathfrak{u}(E))} -\frac{w_{\tau}(A,{\eta})}{||\eta||}. \end{align}
Proposition \ref{PropMaxWeight} shows that $\xi$ determines a holomorphic filtration and orthogonal splitting
		$$\mathcal{E}: \qquad E_1 < E_2 < \cdots < E, \qquad E_j = D_1\oplus \dots \oplus D_j$$
of $(E, \bar{\partial}_A)$. With respect to this orthogonal splitting $\xi$ has the shape
				$$\textbf{i}\xi = \text{diag}(\lambda_1, \lambda_2, \ldots, \lambda_r)$$
with $\lambda_1 < \lambda_2 < \cdots < \lambda_r$ and the weight is given by 
		$$w_{\tau}(A, \lambda) = 2\pi \sum_{j=1}^r \lambda_j \left( c_1(D_j) - \text{rk}(D_j)\mu(E) \right).$$ 
By maximality of the weight $-w_{\tau}(A,\xi)$ we conclude that $\lambda = (\lambda_1,\ldots,\lambda_r)$ is a global minimum of the function
					$$f(x_1,\ldots, x_r) = \sum_{j=1}^r x_j \left(c_1(D_j) -\text{rk}(D_j)\mu(E)\right)$$
on the ellipsoid $\{ \sum_{j=1}^r x_j^2 \text{rk}(D_j) = 1 \}$ under the open condition 
					$$x_1 < x_2 < \cdots < x_r.$$
Since $(E,\bar{\partial}_A)$ is unstable, Proposition \ref{PropStabChar} implies that this minimum is negative and $f$ does not vanish identically. Thus $\nabla f$ vanishes nowhere and $\lambda$ must lie on the ellipsoid. It satisfies there the Lagrange condition 
					$$\left(c_1(D_j) - \text{rk}(D_j) \right)\mu(E)= c \lambda_j \text{rk}(D_j)$$ 
for $j=1,\ldots, r$	and some constant $c \neq 0$. Since $f(\lambda) < 0$ we must have $c < 0$. Since the $\lambda_j$ are increasing this yields
		$$\mu(D_1) > \mu(D_2) > \cdots > \mu(D_r)$$
and $\mathcal{E}$ is a concave filtration of $E$. Solving the Lagrange problem we get 
		$$\lambda_j = \frac{\mu(E) - \mu(D_j)}{\sqrt{\sum_{j=1}^r  \text{rk}(D_j)(\mu(D_j) - \mu(E))^2 }} = \frac{\mu(E) - \mu(D_j)}{\sqrt{||\vec{\mu}(\mathcal{E})||_2^2 - \text{rk}(E)\mu(E)^2}}$$
and
		\begin{align} \label{TDWeq2}-w_{\tau}(A,\xi) = 2\pi \sqrt{||\vec{\mu}(\mathcal{E})||_2^2 - \text{rk}(E)\mu(E)^2}. \end{align}
Now Corollary \ref{HNCor} shows that $\mathcal{E} = \mathcal{E}_{HN}$ must agree with the Harder-Narasimhan filtration of $E$ and $\xi$ is uniquely determined. 

Conversely, we can use the Harder-Narasimhan filtration to define $\xi$ and the argument from above shows that it satisfies (\ref{TDWeq1}). It remains to show it also yields equality in the moment-weight inequality. It follows from the proof of Proposition \ref{PropMaxWeight} that the limit
				$$A_+ := \lim_{t\rightarrow \infty} e^{\textbf{i}t\xi} A$$
exists and splits as $A_+ = A_1 \oplus \cdots \oplus A_r$ with $A_j \in \mathcal{A}(D_j) \cong \mathcal{A}(E_j/E_{j-1})$. The Yang-Mills flow $A_+(t)$ starting at $A_+$ is the product of the Yang-Mills flow on each factor and clearly remains in the closure $\overline{\mathcal{G}^c(A)}$. It follows from Lemma \ref{LemmaSemistableFlow} that the limit $A_{\infty} := \lim_{t\rightarrow \infty}A_+(t)$ of this flow satisfies
				$$F_{A_{\infty}} = -2\pi \textbf{i}			 \begin{pmatrix} \mu(D_1) & 						&					&			 \\
																																					&		\mu(D_2)	&					&				\\
																																					&							&	\ddots	&			 \\
																																					&							&					& \mu(D_r)
																									\end{pmatrix}.$$
Now (\ref{TDWeq0}) and (\ref{TDWeq2}) yield
			$$\inf_{g\in \mathcal{G}^c} ||F_{gA} - \tau|| \leq ||F_{A_{\infty}} -  \tau|| = 2\pi \sqrt{\sum_{j=1}^r \text{rk}(D_j)(\mu(E) - \mu(D_j))^2}= -w(A,\xi).$$
The converse inequality follows from the moment-weight inequality (Theorem \ref{thmMWI}) and this completes the proof in the unitary case.

\subsubsection*{Proof of Theorem \ref{ThmDomWeight} for general compact connected Lie groups $G$.}

Let $G$ be a compact connected Lie group. We show first that one restrict the argument to the case where $Z_0(G)$ is discrete. Recall that the Lie algebra of $G$ decomposes as $\mathfrak{g} = Z(\mathfrak{g})\oplus [\mathfrak{g},\mathfrak{g}]$. The center yields a trivial $Z(\mathfrak{g})$ subbundle $V \subset \text{ad}(P)$ and its orthogonal complement has fiber $[\mathfrak{g},\mathfrak{g}]$ and is canonically isomorphic to $\text{ad}(P/Z_0(G))$. This yields the orthogonal decomposition 
					\begin{align} \label{TDWeq3} \text{ad}(P) \cong V\oplus \text{ad}(P/Z_0(G)). \end{align} 
Let $A \in \mathcal{A}(P)$ and denote by $\bar{A} \in \mathcal{A}(P/Z_0(G))$ the induced connection. Decompose $\xi \in \Omega^0(\Sigma, \text{ad}(P))$ as $\xi = \xi^{z} + \xi^{ss}$ with respect to the decomposition (\ref{TDWeq3}). Then (\ref{eqCT}) and Lemma \ref{WeightsPblLemma} yields
					\begin{align*}  w_{\tau}(A,\xi) = w_{\tau}(A,\xi^{ss}) = w_0(\bar{A}, \xi^{ss}). \end{align*}
Decompose similarly $F_{gA} = F^{z} + F^{ss}$ and note that $F^{ss} = F_{g\bar{A}}$. This yields
				$$||*F_{gA} - \tau||^2 = ||*F^{ss}||^2 + ||*F^z - \tau||^2 \geq ||*F_{g\bar{A}}||^2.$$
As in Lemma \ref{LemmaSSred2} one shows that $g$ can be modified to a gauge transformation $\tilde{g}$ such that $g\bar{A} = \tilde{g}\bar{A}$ and $*F^z = \tau$. Hence
						$$\inf_{g\in\mathcal{G}^c(P)} ||*F_{gA} - \tau|| = \inf_{g \in \mathcal{G}^c(P/Z_0(G))} ||*F_{g\bar{A}}||.$$
This completes the reduction argument.

Now assume that $Z_0(G)$ is discrete and $\tau = 0$. Choose a faithful representation $G \hookrightarrow U(n)$ and identify $G$ with its image in $U(n)$. It follows from Lemma \ref{LemWeightExt} and (\ref{TDWeq0}) that the associated vector bundle $E = P\times_G \mathbb{C}^n$ satisfies $\mu(E) = 0$. For $A \in \mathcal{A}(P)$ Theorem \ref{MLT} yields 
				$$\inf_{g \in \mathcal{G}^c(P)} ||*F_{gA}|| = ||*F_{A_{\infty}}|| = \inf_{g \in \text{GL}(n)} ||*F_{gA}||$$
where we consider $A$ as $G$-connection for the left equality and as $U(n)$-connection for the right equality. It follows from the unitary case that there exists (up to scaling) an unique section $\xi \in \Omega^0(\Sigma, \mathfrak{u}(E))$ satisfying
							$$-\frac{w_0(A,\xi)}{||\xi||} = \inf_{g\in \mathcal{G}^c} ||F_{gA}||.$$
Let $\tilde{\xi}$ be the orthogonal projection of $\xi$ onto $\mathfrak{g}(E) \subset \mathfrak{u}(E)$. Then Lemma \ref{LemWeightExt} shows $w_0(A,\xi) = w_0(A,\tilde{\xi})$ and hence
							$$\inf_{g\in \mathcal{G}^c} ||F_{gA}|| = -\frac{w_0(A,\xi)}{||\xi||} \leq -\frac{w_0(A,\tilde{\xi})}{||\tilde{\xi}||}$$
with equality if and only if $\xi = \tilde{\xi}$. The moment weight inequality (Theorem \ref{thmMWI}) yields the converse inequality and this completes the proof.

\newpage

\bibliographystyle{plain}
\bibliography{references}

\end{document}